\documentclass{amsart}
\usepackage[margin=1in]{geometry}
\usepackage{mine}
\usepackage[utf8]{inputenc}
\renewcommand{\N}{\mathbb{N}}
\DeclareMathOperator{\str}{str}
\DeclareMathOperator{\Kato}{Kato}
\DeclareMathOperator{\disc}{disc}

\renewcommand{\H}{\mathsf H}
\newcommand{\union}{\cup}
\newcommand{\ST}{\mathsf S}
\DeclareMathOperator{\CM}{CM}
\usepackage{mathrsfs}
\newcommand{\s}{\mathsf s}
\renewcommand{\v}{\mathsf v}

\let\rank\relax
\let\corank\relax
\DeclareMathOperator{\rank} {rk}

\DeclareMathOperator{\gr}{gr}
\DeclareMathOperator{\corank}{crk}
\newcommand\sbullet[1][.75]{\mathbin{\vcenter{\hbox{\scalebox{#1}{$\bullet$}}}}}
\usepackage{tikz-cd}
\newcommand{\q}{\mathsf q}
\newcommand{\QT}{\mathsf Q}

\DeclareMathOperator{\Res}{Res}
\DeclareMathOperator{\tors}{tors}
\DeclareMathOperator{\cores}{Cores}

\theoremstyle{plain}
\newtheorem{step}{Step}
\newtheorem{case}{Case}
\newtheorem{construction}[subsubsection]{Construction}
\newtheorem{theo}{Theorem}
\newtheorem{core}[theo]{Corollary}
\theoremstyle{definition}

\newtheorem*{claim}{Claim}

\title{Kolyvagin's Conjecture, bipartite Euler systems, and higher congruences of modular forms}
\author{Naomi Sweeting}
\email{naomisweeting@math.harvard.edu}
\address{Department of Mathematics, Harvard University}
\begin{document}

\begin{abstract}
    Let $E/\mathbb{Q}$ be an elliptic curve and let $K$ be an imaginary quadratic field. Under a certain Heegner hypothesis, Kolyvagin constructed cohomology classes for $E$ using $K$-CM points and conjectured they did not all vanish.  Conditional on this conjecture, he described the Selmer rank of $E$ using his system of classes. We extend work of Wei Zhang to prove new cases of Kolyvagin's conjecture by considering congruences of modular forms modulo large powers of $p $.   Additionally, we prove an analogous result, and give a description of the Selmer rank, in a complementary ``definite'' case (using certain modified $L$-values rather than CM points). 
    Similar methods are also used to improve known results on the Heegner point main conjecture of Perrin-Riou. One consequence of our results is a new converse theorem, that $p$-Selmer rank one implies analytic rank one, when the residual representation has dihedral image.
\end{abstract}

\maketitle
\tableofcontents
\section{Introduction}
Let $f: \T_N \to \O_f$ be a non-CM newform of level $N$, weight 2, and trivial nebentypus, where $\O_f$ is the ring of integers in an algebraic extension of $\Q$. The Birch and Swinnerton-Dyer Conjecture asserts the equality: \begin{equation}r(A/\Q) = \ord_{s=1} L(A, s),\end{equation} where $A=A_f$ is the associated abelian variety to $f$ and $r$ is its Mordell-Weil rank. In pioneering works on this problem, Perrin-Riou \cite{perrin1987fonctions} and Kolyvagin \cite{kolyvagin1991structure, kolyvagin1991sha} studied ranks of elliptic curves over an auxiliary imaginary quadratic field $K$  through the theory of Heegner points on modular curves. We prove, in new cases, conjectures made by both authors.

Fix  a quadratic imaginary field $K $, and a prime $\p\subset\O_f $
of residue characteristic. We write $T_f $
for (a lattice in) the $\p $-adic Galois representation associated to $f $. Assume the following generalized Heegner hypothesis:
\begin{equation}\label{Heegner hypothesis}
    \tag {Heeg}N = N ^ + N ^ - ,\text{ where all }\l | N ^ +\text { are split in } K,\text { all }\l | N ^ -\text { are inert in } K,\text { and } N ^ -\text { is  squarefree,}
\end{equation}
as well as:
\begin{equation}
    \label{good red}\tag{unr} p \nmid 2N\disc (K).
\end{equation}
For purposes of exposition in this introduction, we  also assume:
\begin{equation}
    \label{contains scalar}
    \tag{sclr}\text {The image of the $G_K $ action on }\overline T_f\text { contains a nonzero scalar.}
\end{equation}

To state Kolyvagin's conjecture, assume that the number of prime factors $\nu (N ^ -) $
is even. If $m$ is a squarefree product of primes inert in $K $, one can  use  Heegner points of conductor $m$ 
on the Shimura curve $X_{N^+,N^-}$ to construct classes $$c(m)\in H^1( K, T_f/I_m), $$
where $I_m$ is the ideal of $\O=\O_{f,\p} $
generated by $\l+1$ and $a_\l$ for all $\l | m $. (In the text, $c(m)$ is denoted $\overline c(m, 1)$.) These classes are a mild generalization of the ones constructed by Kolyvagin \cite{kolyvagin1991sha}. 
We are able to prove the following result towards Kolyvagin's conjecture on the nonvanishing of the system $\set{c(m)}$:
\begin{theo}\label{Kolyvagin stupid conjecture Intro}[Corollary \ref{classical vanishing order}]
Assume (\ref{Heegner hypothesis}), (\ref{good red}), and (\ref{contains scalar}) hold for $f,\p, $
and $K $, and $\nu (N ^ -) $
is even. Suppose the following conditions hold:
\[
(\diamondsuit)\hspace{0.5cm} \left\lbrace\parbox{14cm}{$\sbullet$
The modulo $\p$ representation $\overline T_f $
associated to $f $ 
is absolutely irreducible; if $p = 3, $ then $\overline T_f $ is not induced from a character of  $G_{\Q\sqrt{-3}}$. \\ $\sbullet$
If $p$ is inert in $K$ or $a_p$ is not a $\p$-adic unit, then there exists some prime $\l||N$. \\ $\sbullet$ If $a_p$ is not a $\p$-adic unit, then either $\l$ may be chosen above so that $A_f$ has non-split toric reduction at $\l$, or the image of the Galois action on $T_f$ contains a conjugate of $SL_2(\Z_p)$.}\right.
\]
Then there exists a nonzero Kolyvagin class $$0\neq c(m)\in H^1( K, T_f/I_m).$$
\end{theo}

As Kolyvagin observed, Theorem \ref{Kolyvagin stupid conjecture Intro} can be used to give a description of the Selmer ranks $$r ^\pm =\rank_\O\Sel (K, T_f) ^\pm,$$ where superscripts refer to the action of complex conjugation. Indeed, define the vanishing order of the system $\set {c (m)} $ as
\begin{equation}
\nu\coloneqq\min\set {\nu (m)\,:\, c (m)\neq 0}\end{equation}
where as before $\nu $
denotes the number of prime factors.
Then we have:
\begin{core}\label{Selmer rank corollary intro odd}
Under the assumptions of Theorem A, $$\max\set {r ^ +, r ^ -} =\nu +1. $$
Moreover
$r ^ + + r ^ - $
is odd,
and the larger eigenspace has sign $(-1) ^ {\nu +1}\epsilon_f $, where $\epsilon_f $ is the global root number of $f $.
\end{core}
Of course, the latter two assertions follow from the parity conjecture for $f $, already proven
by Nekovar \cite{nekovavr2006selmer}.

Since $c (1)\in\Sel (K, T_f) $
is the Kummer image of the classical Heegner point, the Gross-Zagier formula implies that $L' (f/K, 1)\neq 0 $
if and only if $c (1)\neq 0. $ Hence Corollary B yields a so-called $p $-converse theorem (in fact, under a slightly weaker hypothesis):
\begin{core}\label{converse theorem intro}
Assume that (\ref{Heegner hypothesis}), (\ref{good red}), and Condition $\diamondsuit $ hold for $f $, $\p $, and $K $, and $\nu (N ^ -) $ is even.
Then $$L' (f/K, 1)\neq 0\iff\rank_\O\Sel (K, T_f) = 1\iff\rank_\Z A_f (K) = [\O_f:\Z]. $$
\end{core}
Now suppose instead that $\nu (N ^ -) $
is odd; it turns out  that Kolyvagin's construction, suitably modified, may still be used to relate Selmer ranks and CM points. The Jacquet-Langlands correspondence
associates to $f $
a quaternionic modular form
\begin{equation}
    \xphi_f: X_{N ^ +, N ^ -}\to\O_f,
\end{equation}
where $X_{N ^ +, N ^ -}$
is a double coset space for a definite quaternion algebra, usually called a Shimura set. If $m $
is a squarefree product of primes inert in $K $, there exist analogues of CM points of conductor $m $
on the Shimura set. Using the values of $\xphi_f $
at these points, we construct certain special elements (well-defined up to units)
\begin{equation}\lambda(m)\in \O/I_m\end{equation} ($\lambda (m, 1) $ in the text).
Here the ideal $I_m \subset\O$
is as before.
The elements $\lambda (m) $
encode the same information about the Selmer ranks of $A_f/K $ as Kolyvagin's classes $c (m) $.
\begin{theo}\label{Kolyvagin even rank case}
Suppose that (\ref{Heegner hypothesis}), (\ref{good red}), (\ref{contains scalar}), and Condition $\diamondsuit$
hold for $f,\p, $
and $K $, and that $\nu (N ^ -) $
is odd. Then the vanishing order $$\nu\coloneqq\min\set {\nu (m)\,:\,\lambda (m)\neq 0} $$
is finite and $$\nu=\max\set {r ^ +, r ^ -}. $$
Moreover $(-1) ^ {\nu} =\epsilon_f $ and $r ^ + + r ^ - $
is even.
\end{theo}
As before, the final statement is a consequence of the parity conjecture; we include it only to emphasize that it follows from the non-vanishing of some $\lambda (m) $, in analogy to the indefinite case.
\subsection{Comparison to previous results}
In the indefinite case, the first results towards Kolyvagin's conjecture were obtained by Zhang \cite{zhang2014selmer}, under a number of additional assumptions: that $p\geq 5, $ that the Galois representation associated to $\overline T_f $ is surjective, and additional hypotheses on the residual ramification.  In particular, under the hypotheses of \cite{zhang2014selmer}, there exists a class $c (m) $ whose reduction in $H ^ 1 (K,\overline T_f) $ is nonzero; this is not the case in general.
In the definite case, the classes $\lambda (m) $ are a novel feature of this work and were not considered in \cite{zhang2014selmer}.

The converse theorem we obtain (Corollary \ref{converse theorem intro}) is new in several cases, most notably when the image of the Galois action on $\overline T_f $
is dihedral, or when $p = 3. $ Previous results, under various additional hypotheses,
 were obtained by Zhang as a corollary of his work on Kolyvagin's conjecture, and by Skinner \cite{skinner2020converse} by a purely Iwasawa-theoretic method. For converse theorems in other settings, see Burungale \cite{burungale2020p} for the CM case, Castella-Grossi-Lee-Skinner \cite{castella2020anticyclotomic} for the residually reducible case, Castella-Wan \cite{castella2015iwasawa} for the supersingular case, and Skinner-Zhang \cite{skinner2014indivisibility} for the case of multiplicative reduction.
 
\subsection{Iwasawa theory}

Now suppose again that $\nu(N^-)$ is even. While the Kolyvagin classes are constructed by varying the conductor of CM points on $X_{N^+,N^-}$ over squarefree integers, one may instead $p$-adically interpolate CM points of $p$-power conductor to obtain a class:
\begin{equation}
    \boldsymbol{\kappa}_\infty \in H ^ 1 (K, T_f\otimes\Lambda(\Psi)),
\end{equation}
where $\Lambda =\O\llbracket \Gal (K_\infty/K)\rrbracket$ is the anticyclotomic Iwasawa algebra, given $G_K $-action by the tautological character $\Psi $. (Note that the specialization of $ \boldsymbol{\kappa}_\infty$  at the trivial character is a multiple of  $c(1)$.) The methods used to prove Theorem \ref{Kolyvagin stupid conjecture Intro} also yield the following result towards Perrin-Riou's Heegner point main conjecture.

\begin{theo}[Theorem \ref{Heegner point main conjecture}]\label{Heegner point main conjecture Intro}
Suppose that (\ref{Heegner hypothesis}), (\ref{good red}), and Condition $\diamondsuit$
hold  for $f,\p $, and $K $, and that $\nu (N ^ -) $
is even. Suppose further that
 $a_p $ is a $\p$-adic unit and $p$ splits in $K$. 
Then  there is a pseudo-isomorphism of $\Lambda $-modules: $$\Sel(K_\infty, A_f[\p^\infty])^\vee\approx \Lambda \oplus M\oplus M$$
 for some torsion $\Lambda $-module $M $, and
$$\Char_\Lambda \left(\frac{\Sel(K, T_f\otimes \Lambda)}{\Lambda\cdot\boldsymbol\kappa_\infty}\right)= \Char_\Lambda(M)$$
 as ideals of $\Lambda\otimes\Q_p$. If (\ref{contains scalar}) holds, then the equality is true in $\Lambda. $
\end{theo}
For precise definitions of these Selmer groups and of $\boldsymbol\kappa_\infty $, which is denoted $\boldsymbol\kappa (1) $ in the text, see \S\ref{p adic interpolation}.

Finally, we have the following result on the anticyclotomic main conjecture for $f $
when $\nu (N ^ -) $
is odd. Evaluating the quaternionic modular form $\xphi_f $
on CM points of $p $-power conductor on the Shimura set $X_{N ^ +, N ^ -} $, one constructs the algebraic $p $-adic $L $-function
\begin{equation}
    \boldsymbol\lambda_\infty\in\Lambda,
\end{equation}
denoted $\boldsymbol\lambda (1) $
in the text.
The square of $\boldsymbol\lambda_\infty $
has an interpolation property for twisted $L $-values of $f $.

\begin{theo}[Theorem \ref{ lambda Euler system bound}, Proposition \ref{instructed Euler system overline}]\label{main conjecture Intro}

Suppose that (\ref{Heegner hypothesis}), (\ref{good red}), and Condition $\diamondsuit $
hold for $f,\p $, and $K $, and that $\nu (N ^ -) $
is odd. Suppose further that
 $a_p $ is a $\p$-adic unit and $p$ splits in $K$. 
Then  there is a pseudo-isomorphism of $\Lambda $-modules: $$\Sel(K_\infty, A_f[\p^\infty])^\vee\approx M\oplus M$$
 for some torsion $\Lambda $-module $M $, and
$$(\boldsymbol\lambda _\infty)\subset \Char_\Lambda(M)$$
 as ideals of $\Lambda\otimes\Q_p$. If additionally (\ref{contains scalar}) holds, then the inclusion is
  true in $\Lambda. $
\end{theo}
The opposite inclusion of ideals may be deduced directly from Skinner-Urban's proof of one divisibility in the three-variable main conjecture \cite{skinner2014iwasawa}; indeed, this is an essential ingredient in all of our results, as explained below.
\subsection{Comparison to previous results}
 The technical hypotheses in Zhang's proof of Kolyvagin's conjecture were carried over to Burungale, Castella, and Kim's proof \cite{burungale2019proof} of the lower bound on the Selmer group in the Heegner point main conjecture, where it is also assumed that $p $ is not anomalous. While the methods used in this paper build on those of \cite{ burungale2019proof},  Castella and Wan \cite{castella2020iwasawa} have given a completely independent proof of a three-variable main conjecture when $\nu (N ^ -) $ is even. Their result also requires some hypotheses on residual ramification avoided here, and that $N $ be squarefree.

For upper bounds on the Selmer group in Theorem \ref{Heegner point main conjecture Intro} and Theorem \ref{main conjecture Intro}, various technical assumptions on the residual representation and on the image of the Galois action were used in prior works by Howard \cite{howard2004heegner, howard2006bipartite} and
Chida-Hsieh \cite {chida2015anticyclotomic}.

\subsection{Overview of the proofs}
To prove
Theorems \ref{Kolyvagin stupid conjecture Intro} and \ref{Kolyvagin even rank case}, 
we extend Kolyvagin's construction
to a larger system of classes \begin{equation}\label{Euler system}c (m, Q_1)\in H ^ 1 (K, T_f/\p ^ M),\;\;\lambda (m, Q_2)\in\O/\p ^ M,\end{equation}
where $M $ is a fixed integer, and $m, Q_1, Q_2 $
are squarefree product of auxiliary primes satisfying certain congruence conditions, such that $\nu (N ^ -Q_1) $
is even and $\nu (N ^ -Q_2) $
is odd. 
The classes (\ref{Euler system}) form a bipartite Euler system in the sense of Howard \cite{howard2006bipartite}
for each fixed $m $ and a Kolyvagin system for each fixed $Q_1. $
   If $\nu (N ^ -) $ itself is even, then
the classes $c (m, 1) $ agree with Kolyvagin's original construction. The Euler system relations are of the form:
\begin{equation}
    \loc_q c (m, Q_1)\sim\lambda (m, Q_1q)\sim\partial_{q'} c (m, Q_1qq'),
\end{equation}
where  $q, q' $
are two additional auxiliary primes not dividing $Q_1 $; 
and
\begin{equation}
    \loc_\l^\pm c (m, Q_1)\sim\partial_\l ^\mp c (m\l, Q_1),
\end{equation}
where $\l $
is an additional auxiliary prime not dividing $m $.
(Here $\loc_q $, $\partial_{q'} $, $\loc_\l ^\pm $, $\partial^\pm_\l$
are certain  localization maps landing in subspaces of the local cohomology free of rank one over $\O/\p ^ M. $) The classes $c (m, Q_1) $
were introduced by Zhang, although the $\lambda (m, Q_2) $
are only implicit in \cite{zhang2014selmer}.

If $c (m, Q_1)\neq 0 $,
 then one can use the Kolyvagin system relation to find an auxiliary $\l $  --- either prime or equal to $1$---
 such that $\partial_q c (m\l,Q_1)\neq 0$.
 By the bipartite Euler system relation, this implies $\lambda (m\l, Q_1/q)\neq 0.$
 On the other hand if $\lambda (m, Q_2)\neq 0 $
 and $q | Q_2, $
 then $c (m, Q_2/q)\neq 0. $
 Combining these two observations, we reduce the non-vanishing of some class $c (m, 1) $
 or $\lambda (m, 1) $ --- depending on the parity of $\nu (N ^ -) $ --- to
 exhibiting a single $Q_2 $
 such that $\lambda (1, Q_2)\neq 0. $

Now, if there exists a newform $g$  of level $NQ_2 $ with a congruence to $f $ modulo $\p^M $, then $\lambda (1, Q_2) $ is essentially the reduction of the algebraic $L $-value $L ^ {\text {alg }} (g/K, 1) $ modulo $\p ^ M$, which  is related to the length of the Selmer group of $g$
by the Iwasawa main conjecture \cite{skinner2014iwasawa,wan2016iwasawa}.  To complete the proof, it therefore suffices to choose a  suitable $Q_2 $
and construct such a $g$ with a small Selmer group. We remark that our results can only be obtained by working modulo $\p ^ M $
for a large $M $, since in general it will not be possible  to choose $g $
such that $L ^ {\text {alg}} (g/K, 1) $
is a $\p $-adic unit; in \cite{zhang2014selmer},
$M = 1 $
is fixed throughout, and the need to show that the $L $-value is a unit is responsible for most of the additional hypotheses.

To construct $g $, we use the deformation-theoretic techniques developed by Ramakrishna \cite{ramakrishna2002deforming}. Standard level-raising methods work by producing a modulo $\p $ eigenform of the desired level, and then using that all modulo $\p $ eigenforms lift to characteristic zero, but this is not the case modulo $\p ^ M $. Instead, we deform the representation $T_f/\p ^ M $ to a $\p $-adic Galois representation of a suitable auxiliary level, and then apply modularity lifting to ensure the resulting representation is modular. The auxiliary level $Q_2 $ must be chosen to control two Selmer groups: the adjoint Selmer group governing the deformation problem, and the 
Selmer group $\Sel (K, A_g [\p ^\infty]) $
that is related to the $L $-value. 

We now make some remarks on the construction of the Euler system. The elements $c (m, Q_1) $ (resp.
 $\lambda (m, Q_2) $)
are constructed from CM points of conductor $m $
on the Shimura curve $X_{N ^ +, N ^ - Q_1}$
(resp. Shimura set $X_{N ^ +, N ^ - Q_2} $). Similar Euler system constructions have been made by many authors, e.g. in \cite{chida2015anticyclotomic, bertolini2005iwasawa} as well as in \cite{zhang2014selmer}, but all have relied on certain hypotheses
 ensuring an integral multiplicity one property for the space of  algebraic modular forms on $X_{N ^ +, N ^ - Q_i} $, which we do not impose here. 
Instead, we obtain a  control on the failure of multiplicity one, using the work of Helm \cite{helm2007maps} on maps between Jacobians of modular curves and Shimura curves. The construction of the Euler system is intimately related to level-raising, and so our method can also be viewed as improving results on level-raising of $f$ to algebraic eigenforms modulo $\p^M$ new at multiple auxiliary primes, which had previously been restricted to the multiplicity one case.

 The proof of Theorem \ref{Heegner point main conjecture Intro} is similar to that of Theorem \ref{Kolyvagin stupid conjecture Intro}:
 the $p $-adically interpolated Heegner class $\boldsymbol\kappa_\infty $ is viewed as the bottom layer of an Euler system $\set {\boldsymbol\kappa (Q_1),\boldsymbol\lambda (Q_2)} $. (The squarefree conductor $m $ 
 no longer plays a role.) If $g $, as above, is a newform of level $NQ_2 $ with a congruence to $f $, then $\boldsymbol\lambda (Q_2) $
 is congruent to Bertolini and Darmon's anticyclotomic $p $-adic $L $-function of $g $ \cite{bertolini2005iwasawa}. Using this and an Euler system argument, we reduce the lower bound on the Selmer group in the Heegner point main conjecture to the anticyclotomic main conjecture for $g $, which was proven in \cite{skinner2014iwasawa}.
Finally, the upper bound on the Selmer group in the Heegner point main conjecture, as well as Theorem \ref{main conjecture Intro}, follow by standard arguments from the construction of the Euler system. 

In the text, the arguments described above are phrased in the language of
  ultrapatching, which amounts to a formalism for letting $M$  tend to infinity; this also forces each prime factor of $m $, $Q_1 $, $Q_2 $
  to tend to infinity in order to satisfy the congruence conditions. (The number of prime factors of $m $, $Q_1, $ and $Q_2 $
  remain bounded.) This method was inspired by \cite{scholze2015p}, where ultrapatching was applied to the Taylor-Wiles construction. Our setting is different in that
 we patch Galois cohomology groups and Selmer groups rather than geometric \'etale cohomology groups.
 The benefit of ultrapatching is that it allows us to consider the Euler system classes as characteristic zero objects in patched Selmer  groups, significantly streamlining the Euler system arguments. For instance, with patching, we are able to make precise the heuristic that the non-vanishing of each Euler system class $c(m,Q_1)$ or $\lambda(m, Q_2)$ is equivalent to the $(m , Q_i)$-transverse Selmer group being rank one or zero, respectively, cf. Lemma \ref{ KolyvaginTheorem lemma}.
\subsection*{Structure of the paper}
 In \S\ref{ultra section}, we review basic properties of ultrafilters and introduce patched cohomology and Selmer groups. In \S\ref{Howard section}, we present a simplified version of the theory of bipartite Euler systems that appeared in \cite{howard2006bipartite}, using patched cohomology. In \S\ref{geometry sec}, we recall the geometric inputs that will be used to construct bipartite Euler systems: the work of Helm on maps between modular curves  and Shimura curves, and the  behavior of Heegner points on Shimura curves under reduction and specialization. In \S\ref{ construction section}, we prove the  modulo $\p^M$ level-raising result and present a general framework for constructing bipartite Euler systems out of CM points, which we then specialize in \S\ref{construction per to} for our applications. In \S\ref{deform sec}, we give the deformation-theoretic input to construct the newform $g$ (in fact a sequence $g_n $ satisfying increasingly deep congruence conditions). Finally, we prove the main results in  \S\ref{proof section}. An additional calculation in cyclotomic Iwasawa theory is required for Kolyvagin's conjecture when $p $ is non-ordinary or inert in $K $; this is done in the appendix.
 \subsection*{Acknowledgments}
 I am  grateful to Mark Kisin for suggesting this problem and for his ongoing encouragement. 
 Special thanks are additionally due to Francesc Castella, who first alerted me to the relation between Kolyvagin's conjecture and the Heegner point main conjecture. It is also a pleasure to
thank many other people with whom I had helpful conversations and correspondence over the course of this project: Robert Pollack, Xin Wan, Christopher Skinner, Alexander Petrov, Sam Marks, Aaron Landesman, and  Alexander Smith. This work was supported by NSF Grant \#DGE1745303.

\section{Ultrafilters and patching}\label{ultra section}

\subsection{Ultraproducts}
The facts recalled in this subsection are discussed in more detail in \cite{manning}.
\subsubsection {}
A (non-principal) ultrafilter $\mathfrak F $
for the natural numbers $\N =\set {0, 1,\ldots} $
is a collection of subsets of $\N $
satisfying the following properties:
\begin {enumerate}
\item Every set $S\in\mathfrak F $
is infinite.
\item For every $S\subset\N $, either $S\in\mathfrak F $
or $\N - S
\in\mathfrak F $.
\item If $S_1\subset S_2\subset\N $
and $S_1\in\mathfrak F $, then $S_2\in\mathfrak F $.
\item If $S_1, S_2\in\mathfrak F $, then $S_1\intersection S_2\in\mathfrak F $.
\end {enumerate}
Throughout this paper, we fix once and for all a non-principal ultrafilter $\mathfrak F $
on $\N $, which is possible assuming the axiom of choice.
We will say that a statement $P $
holds for $\mathfrak F $-many $n\in\N $
if the set $S $
of $n $
for which $P $
holds lies in $\mathfrak F $. 
\begin{prop}\label{ultrafilter finite}
Suppose that $\mathcal {C} $
is a finite set and $S\subset\N $
lies in $\mathfrak F $. Then for any function  $t:S\to \mathcal C $, there is a unique $c\in\mathcal C $
such that $t (n) = c $
for $\mathfrak F $-many $n $.
\end{prop}
\begin{proof}
The function $t $
defines a finite partition of $\N $:
$$\N =(\N - S)\sqcup\bigsqcup_{c\in\mathcal C} t ^ {-1} (c). $$
An easy induction argument shows that, for any partition of $\N $
into a finite number sets, exactly one of the sets lies in $\mathfrak F $.
Since $\N - S\not\in\mathfrak F $, the result follows.
\end{proof}
\subsubsection {}\label{equivalence}
If $\mathcal M =\set {M_n} _{n\in\N}$
is a sequence of sets indexed by $\N $, then $\mathfrak F $
defines an equivalence relation $\sim $
on $\prod M_n $:
$$(m_n)_{n\in\N}\sim (m'_n)_{n\in\N}\iff\set {n\,:\, m_n = m'_n}\in\mathfrak F. $$
The quotient $\prod M_n/\sim $
is called the \textbf {ultraproduct} of the sequence $\mathcal M $
and is denoted $\mathcal U (\mathcal M) $.
 The ultraproduct is functorial: let $\mathcal M' =\set {M_n'} $
 be another sequence of sets and suppose given, for $\mathfrak F $-many $n $,  maps $\phi_n: M_n\to M'_n$. Then there is a natural map $\phi^ {\mathcal{U}}:\mathcal U (\mathcal M)\to\mathcal U (\mathcal M') $.
 Similarly, if $R $
 is a fixed (topological) ring and each $M_n $
 is a (continuous) $R $-module, then $\mathcal U (\mathcal M) $
 may be naturally endowed with the structure of a (continuous) $R $-module; in particular, if each $M_n $
 is an  abelian group, then $\mathcal U(\mathcal M)$
 has a natural abelian group structure as well.
 
The following basic properties are proven in \cite [\S I.1] {manning} using the functoriality of the ultraproduct.
\begin{prop}\label{basic properties}
Suppose that each $M_n $
is a finite set, and that $\#M_n <C $
for some constant $C $
and for $\mathfrak F $-many $n $. Then:
\begin {enumerate}
\item $\mathcal U (\mathcal M) $
is finite and $\#\mathcal U (\mathcal M) =\#M_n $
for $\mathfrak F $-many $n $.
\item Suppose that each $M_n $
is additionally endowed with the structure of a (continuous) $R $-module, and let $A $
be another (continuous) $R $-module.
Given a family of isomorphisms $\phi_n: M_n\xrightarrow\sim A $
for $\mathfrak F $-many $n $, there is an induced isomorphism $$\phi ^ {\mathcal U}:\mathcal U (\mathcal M)\xrightarrow\sim A. $$

\end{enumerate}
\begin {rmk}
If $R $
is (topologically) finitely generated, then there
are only finitely many isomorphism classes
of (continuous) $R $-modules
of a fixed cardinality. Hence, if $\mathcal M $
is a sequence of finite $R $-modules of bounded cardinality, Proposition \ref {ultrafilter finite}
and Proposition \ref {basic properties}
together imply that $\mathcal U (\mathcal M)$
is non-canonically isomorphic to $\mathfrak F $-many $M_n $.
\end {rmk}
\end{prop}
\begin {prop}\label{exactness}
Let $\mathcal C $
be the category of sequences of (continuous) $R $-modules of uniformly bounded cardinality. Then $\mathcal U $
is exact as a functor from $\mathcal C $
to the category of (continuous) $R $-modules.
\end {prop}
\begin {proof}
We wish to show that $\mathcal U $
preserves finite limits and
colimits. Since any given finite limit or colimit
in $\mathcal C $
involves only sequences of $N $-torsion $R $-modules, for some integer $N $, the limit or colimit may be computed in the category of $\Z/N \Z$-modules, in which case \cite [Proposition I.2.2] {manning}
applies.
\end {proof}
 \subsection{Ultraprimes}
 \subsubsection{}
 Fix a number field $L $
 and let $M_L $
 be its set of  places;
 for each $v\in M_L $,  
 fix as well
  an embedding $\overline L\hookrightarrow\overline{L_v} $. If $\mathcal M_L $
 is the constant sequence of sets $\set {M_L}_{n\in\N} $, then we define the set of ultraprimes of $L $
 as $$\mathsf M_L =\mathcal U (\mathcal M_L). $$
 By definition, an ultraprime $\mathsf v\in\mathsf M_L $
 is an equivalence class of sequences $(v_n)_{n\in\N} $, where each $v_n $
 is a place  of $L $; note that $\Gal (L/\Q) $
 acts by set automorphisms on $\mathsf M_L $, compatibly with the natural projection $\mathsf M_E\to\mathsf M_L $
 for a finite extension $E/L $.
 The map $v\mapsto (v, v,\ldots) $
 induces an embedding $M_L\hookrightarrow\mathsf M_L $, written $v\mapsto\underline v $, and we say an ultraprime is constant if it lies in the image of this embedding. 
 
  \begin{prop}
    Let $\mathsf v $
    be a non-constant ultraprime.
  Then there
  exists a unique Frobenius element $\Frob_\v\in\Gal (\overline L/L) $
  with the following property: for each finite Galois extension $L\subset E\subset\overline L $,
  and for any representative
  $(v_n) $
  of $\mathsf v $,
  there are $\mathfrak F $-many $n $
  such that
  $v_n $
  is unramified in $E/L $
  and the Frobenius of $v_n $
  in $\Gal (E/L) $
  is the natural image of $\Frob_\v$.
  \end{prop}
  
  \begin{proof}
  Let $(v_n)_{n\in\N} $
  be a representative of $\mathsf v $, and fix
  for the time being a finite extension $E/L $.
  If $v_n $
  is archimedian or ramified in $E $
  for $\mathfrak F $-many $n $, then Proposition \ref{ultrafilter finite}
  implies that $\mathsf v $
  is constant. Thus the map that sends $n $
   to the Frobenius of $v_n $
   in $\Gal (E/L) $
  is defined for $\mathfrak F $-many $n $;
  by Proposition \ref{ultrafilter finite},
  it sends $\mathfrak F $-many $n $
  to a (unique) common value $g_E\in\Gal (E/L) $.
  Note that $g_E $
  does not depend on the representative $(v_n) $.
  By the uniqueness of $g_E $, the association $E\mapsto g_E $
  is compatible with restriction to subextensions $E'\subset E $, hence defines an element of the absolute Galois group.
  \end{proof}
  \subsubsection{}
  Let $\mathsf v $
  be an ultraprime. We define its abstract Galois group $\mathsf G_{\mathsf v} $
  as $\Gal (\overline {L_v}/L_v) $
  if $ \mathsf v =\underline v $
  is constant, and as the semi direct product
  $$\widehat\Z (1)\rtimes\langle\Frob_\v \rangle$$
  otherwise. Here, $\langle\Frob_\v\rangle$
  denotes the free profinite group on one generator, acting on $\widehat\Z (1) $
  by $\Frob_\v$. We  define the inertia group $\mathsf I_{\mathsf v} \subset\mathsf G_{\mathsf v}$ of $\mathsf v $
  to be the usual inertia group if $\mathsf v $
  is constant, and the normal subgroup $\widehat\Z (1)\subset\mathsf G_{\mathsf v} $
  otherwise. 
  \subsection{Local cohomology}

  \subsubsection{}For any (continuous) Galois module $A $
  defined over $L $,
  and for any $\mathsf v\in\mathsf M_L $, there is a natural action of $\mathsf G_{\mathsf v} $
  on $A $
  (factoring through the quotient $\mathsf G_{\mathsf v}\to\Frob_\v$
  if $\mathsf v $
  is nonconstant).
  We define local cohomology groups  by:
  \begin{align*}
  &\mathsf H ^ i (L_{\mathsf v}, A)\coloneqq H ^ i_{cts} (\mathsf G_{\mathsf v}, A),\\
  &\mathsf H ^ i (L_{\mathsf v} ^ {nr}, A)\coloneqq H ^ i_{cts} (\mathsf I_{\mathsf v}, A),
  \;\;\; i\geq 0. 
  \end {align*}
  
 Note that the local cohomology  commutes with direct limits and countable inverse limits of finite, discrete Galois modules; the former is essentially by definition
 of continuous cohomology and the latter is by \cite[Corollary 2.6.7]{neukirch2013cohomology} applied to $\mathsf G_{\mathsf v},\mathsf I_{\mathsf v}. $
  \begin {prop}\label{local cohomology}
  Let $\mathsf v\in\mathsf M_L $
  be an ultraprime represented by a sequence $(v_n)_{n\in\N} $.
  If $A $
  is a finite, discrete Galois module over $L $, then for $\mathfrak F $-many $n $
  there are natural isomorphisms (compatible with the restriction maps and with the cup product):
  \begin {align*}
  &  H ^ i (L_{ v_n}, A)\simeq\mathsf H ^ i (L_{\mathsf v}, A),\\
& H ^ i (L_{v_n} ^ {nr}, A)\simeq \mathsf H ^ i (L ^ {nr}_{\mathsf v} , A),\;\;\; i\geq 0.\end {align*}

\end {prop}
\begin{proof}
 If $\mathsf v $
  is the constant ultraprime $\underline v $, then $v_n = v  $
  for $\mathfrak F $-many $n $, and the desired isomorphisms are given by the identity maps; so suppose $\mathsf v $
  is nonconstant. 
  For $\mathfrak F $-many $n $, the action of the decomposition group $G_{v_n} $ at $v_n $
  on $A $
  is unramified and the Frobenius of $v_n $
  acts by $\Frob_\v $.
  Let $\l_n $
  be the prime of $\Q $
  lying under $v_n $;  since $L/\Q $
  is a finite extension and $A $
  is a finite Galois module, for $\mathfrak F $-many
  $n $
  we have $\l_n\nmid | A | $.
  Restricting to these $n $, the inflation map induces isomorphisms:
  $$H ^ i (G_{v_n} ^ t, A)\simeq H ^ i(L_{v_n}, A),\;\;H ^ i (I_{v_n} ^ t, A)\simeq H ^ i (L_{v_n} ^ {nr}, A), $$
  where $G_{v_n} ^ t $
  and $I_{v_n} ^t$
  denote the tame quotients.
  The tame Galois group $G_{v_n} ^ t $
  is  identified with the semi direct product:
  $$I_{v_n} ^ t\rtimes\langle\Frob_{v_n}\rangle \simeq\widehat\Z ^ {(\l_n)} (1)\rtimes\langle\Frob_{v_n}\rangle. $$
  Since $\Frob_{v_n} $
  and $\Frob_\v$
  may act differently
  on the Tate twist, $G_{v_n} ^ t $
  and $\mathsf G_{\mathsf v} $
  cannot be compared directly; we wish to show that the cohomologies are nonetheless canonically isomorphic for $\mathfrak F $-many $n $. 
  
  Let $G = I\rtimes\langle F\rangle $
  be an abstract group, where $I $
  is abelian and profinite, and $\langle F\rangle$
  denotes the free profinite group on one generator, acting on $I$ by an automorphism. 
  If $A $
   is a $\Z [F] $-module,
   then 
    the Galois cohomology groups $H ^ i (G, A) $
  and $H ^ i (I, A) $
  depend only on $A $
  and $\Hom (I, A) $
  as $\Z [F] $-modules; in particular, the cohomology groups  for $G$
  are canonically isomorphic to the cohomology groups for its quotient $I/| A |\rtimes\langle F\rangle $, and similarly for $I $
  and $I/| A | $.
  Applying this to $G_{v_n} ^ t $
  and $\mathsf G_{\mathsf v} $
  completes the proof, since $\Frob_{v_n} $
  and $\Frob_\v$
  have the same action on the finite Tate twist $\Z/| A | (1)$
  for $\mathfrak F $-many $n $.
  
\end{proof}
\subsection{Patched cohomology}
 
  \subsubsection{}
  Let $\mathsf S\subset\mathsf M_L $
  be a finite set of ultraprimes $\set {\mathsf
s_1,\mathsf s_2,\ldots,\mathsf s_r} $. A \textbf{representative} of $\ST$ is a sequence of sets $S^n\subset M_L$ such that $S^n = \set{s_1^n, \cdots,  s_r^n}$ for some sequences $(s_i ^ {n})_{n\in\N} $ representing $\mathsf s_i$.
If $A$ is a $\Gal(\overline L/L)$ module, we say $A$ is unramified outside $\ST\subset \mathsf M_L$ if it is unramified outside $\ST \cap M_L$. 

\begin{definition}
 Let  $A$ be a topological $\Gal (\overline L/L) $-module
 unramified outside a finite set  $\ST\subset \mathsf M_L$, represented by a sequence $S^n\subset M_L$. If $A$ is profinite, then we define the $i $th unramified-outside-$\ST$ patched cohomology, for all $i\geq 0, $ by:
$$\H^i(L^\ST/S, A) = \lim_{\substack {\longleftarrow\\A\twoheadrightarrow  A'}}
\mathcal{U }\left (\set {H ^ i (L ^ {S ^ n}/L, A')}\right)_{n\in\N}, $$
where the inverse limit runs over continuous finite quotients of $A $.
If $A$
is ind-finite, then its unramified-outside-$\ST $
patched cohomology is defined as:
$$\H ^ i (L ^\ST/L, A ) =\lim_{\substack{\longrightarrow\\A'\subset A}}\mathcal U\left (\set {H ^ i (L ^ {S ^ n}/L, A')}_{n\in\N}\right), $$
where the direct limit runs over finite submodules. If $A $
is either  profinite or ind-finite, then the totally patched cohomology is defined as $$\H ^ i(L, A) =\lim_{\substack {\longrightarrow\\\mathsf S\subset\mathsf M_L}}\H ^ i (L ^ {\ST}/L, A), $$
where the direct limit is over finite subsets and the transition maps are induced by the functoriality of the ultraproduct.
\end{definition}
\begin {rmk}
\begin{enumerate}
\item
    To see that these cohomology groups are well-defined, first note that they are independent of the choice of $S ^ n $
    since any two representatives of a finite set $\ST\subset\mathsf M_L $
    agree for $\mathfrak F $-many $n $. Moreover,
    if $A $
    is both profinite and ind-finite, then it is finite,
    and it is clear that either definition gives the same cohomology groups.
    \item There is a canonical isomorphism $\H ^ 0 (L ^\ST/L, A) = H ^ 0 (L, A) $ for all finite $\ST\subset\mathsf M_L $ and all profinite or ind-finite  $A $.
        \item The assignments $$A\mapsto\H ^ i (L ^\ST/L, A),\;\;\; A\mapsto\H ^ i (L, A) $$
    are functorial in $A $. If $A $
    is an $R $-module for some ring $R $, then each patched cohomology group $\H ^ i (L ^\ST/L, A), $
    $\H ^ i (L, A) $
    has a natural $R $-module structure.
    \item In practice, we will want our profinite Galois modules to be \textbf{countably profinite}, i.e. to have a presentation as a countable inverse limit of finite, discrete topological Galois modules. The significance of this technical hypothesis is that countable inverse limits of finite abelian groups are exact. For example, see  \cite[Corollary 2.7.6]{neukirch2013cohomology}.
    \item Suppose $A$ is ind-finite or countably profinite. If every ultraprime in $\mathsf S $
    is constant, and $S\subset M_L $
    is the corresponding finite set of places, then $\H ^ i (L ^\ST/L, A) $
    is canonically isomorphic to $H ^ i (L ^ S/L, A) $.
    \item  Suppose $A$ is ind-finite or countably profinite. For each ultraprime $\mathsf v $, there are  natural localization maps $$\Res_\v:\H ^ i (L, A)\to\H ^ i (L_\v, A) $$
    deduced from Proposition \ref{local cohomology} (and from \cite[Corollary 2.7.6]{neukirch2013cohomology} applied to $\mathsf G_\v$ in the profinite case).
    \item If the Galois action on $A $
    is the restriction of an action of $G_K $, where $L/K $
    is a Galois extension, then $\Gal (L/K) $ acts naturally on $\H ^ i (L, A) $, again by functoriality of ultraproducts; this is compatible with the localization maps in the obvious way.
    \end{enumerate}
    \end{rmk}

\begin {lemma}\label{finiteness}
For any finite set of primes $S\subset M_L $, and any finite Galois module $A $ over $L $, the cardinality of $H^i(L^S/L, A)$ is uniformly bounded, with a bound depending only on $A$, $L$, and $|S|$. In particular, if
$\ST\subset\mathsf M_L $
is finite, then the patched cohomology groups $\H ^ i (L ^\ST/L, A) $ are finite for each finite Galois module $A $
and each $i\geq 0$.
\end {lemma}
\begin {proof}
The first claim is easily seen from
 \cite[Theorem 4.10]{milne2006arithmetic}; the second follows by
 Proposition \ref{basic properties}.
 \end {proof}
\begin{prop}\label{unramified kernel}
If $A $
is either countably profinite or ind-finite, then, for all $i $, the natural map induces an isomorphism
$$\H ^ i (L^\ST/L, A)\simeq\ker\left (\H ^ i (L, A)\to\prod_{\mathsf v\in \mathsf M_L -\ST}\H ^ i (L_\v ^ {nr}, A)\right). $$
\end{prop}
\begin{proof}
 It suffices to show that, for all finite sets $\mathsf T\subset\mathsf M_L -\ST $,
 $$\H ^ i (L^\ST/L, A)\simeq\ker\left (\H ^ i (L^ {\ST\cup \mathsf T}, A)\to\prod_{\mathsf t\in\mathsf T}\H ^ i (L_{\mathsf t} ^ {nr}, A)\right). $$
 This holds when $A $
 is finite by Lemma \ref{finiteness} and Proposition \ref{exactness}; the general case follows by taking limits.
\end{proof}
\begin {lemma}\label {LES}
 Let $$0\to A\to B\to C\to 0 $$
be an exact sequence of either countably profinite or ind-finite Galois modules unramified outside $\ST.$
Then there is an induced long exact sequence
beginning:
\begin {align*}
0 &\to\mathsf H ^ 0 (L ^ {\mathsf S}/L, A)\to\mathsf H ^ 0 (L ^ {\mathsf S}/L, B)\to\mathsf H ^ 0(L ^ {\mathsf S}/L, C)\to\\
&\to\mathsf H ^ 1 (L ^ {\mathsf S}/L, A)\to\mathsf H ^ 1 (L ^ {\mathsf S}/L, B)\to\cdots
\end {align*}
\end {lemma}
\begin {proof}
If $A $, $B $, and $C $
are all finite, then this follows from Proposition \ref {exactness} and Lemma \ref{finiteness}.

 Now suppose that $A $, $B $, and $C $
are all profinite.
Let $I $, $J $, and $K $
be directed sets indexing the finite  quotients $A\twoheadrightarrow A_i $, $B\twoheadrightarrow B_j $, and $C\twoheadrightarrow C_k $, respectively. We define morphisms of directed sets $t: J\to I $
and $s: J\to K $
by $$A_{t (j)} =\im (A\to B_j),\;\; C_{s (j)} = B_j/A_{t (j)}. $$
Because the subgroup and quotient topologies
on $A $
and $C $
agree with the profinite topologies, the images of $t $
and $s $
are cofinal in $I $
and $K $, respectively.
We therefore have: $$\mathsf H ^ \ast(L ^ {\mathsf S}/L, A) =\lim_{\substack {\leftarrow\\j\in J}}\mathsf H ^ \ast (L ^ {\mathsf S}/L, A_{t (j)}),\;\;\mathsf H ^ \ast (L ^ {\mathsf S}/L, C) =\lim_{\substack {\leftarrow\\j\in J}}\mathsf H ^ \ast (L^ {\mathsf S}/L, C_{s (j)}). $$
For each $j $, we have a long exact sequence
associated to the short exact sequence of finite Galois modules
$$0\to A_{t (j)}\to B_j\to C_{s (j)}\to 0;$$
by Lemma \ref{finiteness}, each term in the long exact sequence is finite.
Since countable inverse limits of finite abelian groups
are exact, taking limits completes the proof.
The ind-finite case is completely analogous.\end {proof}
 \subsection{Selmer structures and patched Selmer groups}
\begin{definition}
Let $A $
be a countably profinite or ind-finite $\Z_p [G_L]$-module.
A \textbf{generalized Selmer structure}
$(\mathcal F,\mathsf S) $
for $A $ 
consists of:
\begin {itemize}
\item a finite set $\mathsf S\subset\mathsf M_L $
containing all Archimedian places, all places over $p $, and all ramified places for $A $;
\item for each $\mathsf v\in\mathsf M_L $, a closed $\Z_p $-submodule (the \textbf{local condition}) $$\mathsf H ^ 1_{\mathcal F} (L_{\mathsf v}, A)\subset\mathsf H ^ 1 (L_{\mathsf v}, A) $$
such that $$\H ^ 1_{\mathcal F} (L_\v, A) =\mathsf H ^ 1_{\unr}  (L_{\mathsf s}, A)\coloneqq\ker\left (\mathsf H ^ 1 (L_{\mathsf s}, A)\to\mathsf H ^ 1 (L ^ {nr}_{\mathsf s}, A)\right) $$
for all $\v\not\in\ST $.
\end {itemize}
If $A $
is an $R$-module for some ring $R $ and $G_L $
acts on $A $
by $R $-module automorphisms, a Selmer structure for $A $ \textbf{over} $R $ is a Selmer structure such that every local condition is an $R $-submodule.
\end{definition}
\subsubsection{}

If $B\subset A $
is any closed Galois-stable submodule, then a Selmer structure $(\mathcal F,\mathsf S) $
for $A $
induces Selmer structures on $B $
and $A /B$
defined in the usual way:
\begin {align*}\mathsf H ^ 1_{\mathcal F} (L_{\mathsf s}, B) =\ker\left (\mathsf H ^ 1 (L_{\mathsf s}, B)\to\frac {\mathsf H ^ 1 (L_{\mathsf s}, A)} {\mathsf H_{\mathcal F} ^ 1 (L_{\mathsf s}, A)}\right),\\\mathsf H  ^ 1_{\mathcal F} (L_{\mathsf s}, A/B) =\im\left (\mathsf H ^ 1_{\mathcal F} (L_{\mathsf s}, A)\to\mathsf H ^ 1 (L_{\mathsf s}, A/B)\right).
\end {align*}

\subsubsection{}
To a generalized Selmer structure we 
associate the \textbf {patched Selmer group}, defined by the exact sequence:
\begin{equation}0\to\Sel_{\mathcal F} (A)\to\mathsf H ^ 1 (L^\ST/L, A)\to\prod_{\mathsf s\in\mathsf S}\frac {\mathsf H ^ 1 (L_{\mathsf s}, A)} {\mathsf H ^ 1_{\mathcal F} (L_{\mathsf s}, A)},\end{equation}
or equivalently (by Proposition \ref{unramified kernel}):
\begin{equation}
    0\to\Sel_{\mathcal F} (A)\to\mathsf H ^ 1 (L, A)\to\prod_{\mathsf s\in\mathsf S}\frac {\mathsf H ^ 1 (L_{\mathsf s}, A)} {\mathsf H ^ 1_{\mathcal F} (L_{\mathsf s}, A)}\times\prod_{\mathsf s\not\in\ST}\H ^ 1 (L_{\mathsf s} ^ {nr}, A).
\end{equation}
(Note that the Selmer group attached to a Selmer structure does not depend on the choice of set $\ST$  but only on the local conditions; we will therefore
sometimes omit $\ST $ from the notation when there is no risk of confusion.)
\subsubsection{}
If $B\subset A $
is Galois-stable, and $B, A/B $
are equipped with the induced Selmer structures,
then by definition there are natural maps:
$$\Sel_{\mathcal F} (B)\to\Sel_{\mathcal F} (A)\to\Sel_{\mathcal F} (A/B). $$

\begin {prop}\label{Selmer proposition investment}
Let $(\mathcal F,\mathsf S) $
be a generalized Selmer structure for $A $. If $A $
is countably profinite and  each continuous finite  quotient  $A \twoheadrightarrow A'$ is equipped
with the Selmer structure induced by $\mathcal F $, then:
$$\lim_{\longleftarrow}\Sel_{\mathcal F} (A')\simeq\Sel_{\mathcal F} (A). $$
If instead $A $
is ind-finite and each finite submodule $A'\subset A $
is given its induced Selmer structure, then:
$$\lim_{\longrightarrow} \Sel_{\mathcal F} (A')\simeq\Sel_{\mathcal F} (A). $$
\end{prop}
\begin {proof}
We show the countably profinite case; the ind-finite case is similar. By definition, $\Sel_{\mathcal F} (A) $
is the kernel of $$\lim_{\longleftarrow}\mathsf H ^ 1 (L ^ {\mathsf S}/L, A')\to\prod_{\mathsf s\in\ST}\frac{\mathsf  H ^ 1 (L_{\mathsf s}, A)} {\mathsf H ^ 1_{\mathcal F} (L_{\mathsf s}, A)},$$
whereas\begin{align*}
    \lim_{\longleftarrow}\Sel_{\mathcal F} (A') & =\lim_{\longleftarrow}\ker\left (\mathsf H ^ 1 (L ^ {\mathsf S}/L, A')\to\prod_{\mathsf s\in\ST}\frac {\H ^ 1 (L_{\mathsf s}, A')} {\H ^ 1_{\mathcal F} (L_{\mathsf s}, A')}\right)
     =\ker\left (\H ^ 1 (L ^\ST/L, A)\to\lim_{\longleftarrow}\frac {\H ^ 1 (L_{\mathsf s}, A')} {\H ^ 1_{\mathcal F} (L_{\mathsf s}, A')}\right).
\end{align*} 
Since $\H ^ 1_{\mathcal F} (L_{\mathsf s}, A) $
is closed, we have $$\lim_{\longleftarrow}\H ^ 1_{\mathcal F} (L_{\mathsf s}, A') =\H ^ 1_{\mathcal F} (L_{\mathsf s}, A), $$
which implies the result.
\end {proof}
\subsubsection{}
Given two Selmer structures $(\mathcal F,\ST) $
and $(\mathcal G,\mathsf T) $
for $A $, we may define Selmer structures $(\mathcal F + \mathcal G,\ST\union\mathsf T)$ and $(\mathcal F \intersection\mathcal G, \ST\union\mathsf T)$ by the local conditions:
\begin{equation}
   \H^1_{ \mathcal  F+\mathcal{G}}(L_\v,A) =\H ^ 1_\mathcal F (L_\v, A)+\H ^ 1_\mathcal G (L_\v, A),\;\;\;\H ^ 1_{\mathcal  F\intersection\mathcal G} (L_\v, A) =\H ^ 1_\mathcal F (L_\v, A)\intersection\H ^ 1_\mathcal G (L_\v, A).
\end{equation}
\subsection{Dual Selmer groups}
\subsubsection{}\label {define pairings}
Fix an
ultraprime $\v\in\mathsf M_L $.
 If $A $
is a countably profinite $\Z_p $-Galois module and $A ^\ast $
denotes the Cartier dual, then the cup product  induces pairings:
\begin {equation}\label{profinite duality}
\langle\cdot,\cdot\rangle_{\mathsf v}:\mathsf H ^ i (L_{\mathsf v}, A)\times\mathsf H^ {2 - i} (L_{\mathsf v}, A ^\ast)\to\Q_p/\Z_p,\;\;\; i = 0, 1, 2.\end {equation}
\begin{prop}\label{Perfect pairing}
The pairing (\ref {profinite duality})
is perfect if $\v$ is non-Archimedean. Moreover, the induced pairing $$\H ^ 1 (L ^\ST/L, A)\times\H ^ 1 (L ^\ST/L, A ^\ast)\to\prod_{\s\in\ST}\H ^ 1 (L_\s, A)\times\H ^ 1 (L_\s, A ^\ast)\xrightarrow {\Sigma\langle\cdot,\cdot\rangle_\s}\Q_p/\Z_p $$
is identically zero.
\end{prop}
\begin {proof}
For the perfectness of (\ref{profinite duality}), the usual proof of Poitou-Tate duality applies equally well to $\mathsf G_{\mathsf v} $; alternatively, one may take limits using Proposition \ref{local cohomology}. 
The second claim is clear when $A $
is finite by functoriality of the ultraproduct, and the general case follows by taking limits.
\end {proof}
\subsubsection{}
Suppose that $A $
is either countably profinite or countably ind-finite, i.e. the Pontryagin dual of a countably profinite Galois module.
If $(\mathcal F,\ST) $
is any Selmer structure for $A $, then we define the dual Selmer structure $(\mathcal F ^\ast,\ST) $
for $A ^\ast $
by:
$$\H ^ 1_{\mathcal F ^\ast} (L_\s, A ^\ast) =\H ^ 1_{\mathcal F} (L_\s, A) ^\perp.  $$
Here $\perp $
denotes the orthogonal complement under either the pairing of (\ref {define pairings}), or the usual modified Tate pairing of \cite[Theorem 2.17]{darmon2000fermat} at Archimedian places.
 We observe that the dual Selmer structure to $(\mathcal F ^\ast,\mathsf S) $
is again $(\mathcal F,\mathsf S) $. 
When$A $
is finite, the dual Selmer groups are related by the
 Greenberg-Wiles formula:
\begin {prop}\label {Greenberg Wiles}
Let $(\mathcal F,\ST) $
be a Selmer structure for a finite $\Z_p[G_L]$-module $A $.
We have: $$\frac {\#\Sel_{\mathcal F} (A)} {\#\Sel_{\mathcal F ^\ast} (A ^\ast)} =\frac {\#\H ^ 0 (L ^\ST/L, A)} {\#\H ^ 0 (L ^\ST/L, A ^\ast)}\prod_{\mathsf s\in\ST}\frac {\#\H ^ 1_{\mathcal F} (L_{\mathsf s}, A)} {\#\H ^ 0 (L_{\mathsf s}, A)}. $$
\end {prop}
\begin {proof}
This follows from \cite[Theorem 2.19]{darmon2000fermat} by the exactness of ultraproducts
and Proposition \ref{basic properties}(1).
\end {proof}
\subsection{Selmer groups over discrete valuation rings}\label{DVR Selmer section}
\subsubsection{}
Let $R $
be a discrete valuation ring with uniformizer $\pi$ which is a finite, flat extension of $\Z_p $, and suppose that $A = T $
is a free $R $-module of finite rank, with $G_L $
action through $R $-module automorphisms. In particular, $T $
is countably profinite. Suppose $\ST\subset\mathsf M_L $
is a finite set  containing all Archimedian places and all places over $p $, such that $T $
is unramified outside $\ST $. 
If $T ^\dagger =\Hom_R (T,R (1)) $
is the dual, then the cup product induces a local Tate pairing \begin{equation}\label{ imperfect pairing}
    \langle\cdot,\cdot\rangle_\v:\H ^ 1 (L_\v, T)\times\H ^ 1 (L_\v,T ^\dagger)\to R.
\end{equation}
\begin{prop}\label{ dualityover DVR}
The kernels on the left and right of (\ref{ imperfect pairing}) are the $R$-torsion submodules; moreover, the induced pairing $$\H ^ 1 (L ^\ST/L, T)\times\H ^ 1 (L ^\ST/L, T ^\dagger)\to\prod_{\v\in\ST}\H ^ 1 (L_\v, T)\times\H ^ 1 (L_\v,T ^\dagger)\xrightarrow {\sum\langle\cdot,\cdot\rangle_\v} R $$
is identically zero.
\end{prop}
\begin{proof}
 This follows from Proposition \ref{Perfect pairing}.
\end{proof}
Given a Selmer structure $(\mathcal F ,\ST)$
for $T $ over $R $, taking the orthogonal complement of each local condition under (\ref { imperfect pairing}) yields a Selmer structure $(\mathcal F ^\dagger,  \ST) $
for $T ^\dagger $. Note that $\mathcal F ^ {\dagger\dagger}\neq\mathcal F $
in general, but we always have $\mathcal F ^ {\dagger\dagger\dagger} =\mathcal F ^\dagger $.
\begin{prop}\label{two different duals}
Let $(\mathcal F,\ST) $
be a Selmer structure such that $$\frac{\H ^ 1 (L_\v, T)}{\H ^ 1_{\mathcal F} (L_\v, T) }$$
is torsion-free. Then, for all $j $ and all $\v\in \mathsf M_L $,
$$\H ^ 1_{\mathcal F ^\ast} (L_\v, T ^\ast [\pi ^ j]) =\H ^ 1_{\mathcal F ^\dagger} (L_\v, T ^\dagger/\pi ^ j) $$
under an identification $T ^\ast [\pi ^ j]\simeq T^\dagger/\pi ^ j $, and in particular
$$\Sel_{\mathcal F ^\ast} (T ^\ast [\pi ^ j]) =\Sel_{\mathcal F ^\dagger} (T^\dagger/\pi ^ j).$$
\end{prop}
\begin{proof}
 Although this fact is presumably standard,
 we give a proof for lack of a reference.
 For ease of notation, we abbreviate $\H ^ i (T ^\dagger) =\H ^ i (L_\v, T ^\dagger) $, etc. and $R_j= R/\pi ^ j $. The choice of uniformizer induces an identification $T ^\dagger\otimes_R (R [1/\pi]/R)\simeq T ^\ast $ and an embedding $T ^\dagger/\pi ^ j\hookrightarrow T ^\ast $;
 let $\H ^ 1 _{\mathcal F^\ast}(T ^\dagger/\pi ^ j)$
 be the induced local condition from this embedding.
 Consider the following commutative diagram with exact rows:
 \begin{center}
\begin{tikzcd}[column sep=small]
0\arrow [r] &\H ^ 0 (T ^\ast)_{/\div}\arrow [r]\arrow [d, "\alpha"] &\H ^ 1_{\mathcal F^\dagger} (T ^\dagger)\arrow [r]\arrow [d, "\beta"] &\Hom (\H ^ 1(T), R)\arrow [r]\arrow [d, "\gamma"] &\Hom (\H ^ 1_{\mathcal F} (T), R)\arrow [d, "\delta"]\arrow [r] & 0
\\
0 \arrow [r]&\H^ 0 (T ^\ast)/\pi ^ j\arrow [r] &
\H ^ 1 _{\mathcal F^\ast}(T ^\dagger/\pi ^ j)\arrow[r] &\Hom (\H ^ 1 (T), R_j)\arrow [r] &\Hom (\H ^ 1 _{\mathcal F}(T), R_ j)\arrow[r]& 0
\end{tikzcd}
\end{center}
Here, the first horizontal map on each row is the Kummer map, and the subscript $/\div$
 refers to the quotient by the maximal divisible submodule. 
 By the hypothesis on $\H ^ 1_{\mathcal F} (T) $, the maps $\coker \gamma\to\coker\delta $ and $\ker\gamma\to \ker\delta $
 are  injective and surjective, respectively. 
Also, $\alpha $
 is clearly surjective. 
 Breaking the diagram into two and applying the snake lemma, it follows that $\beta $
 is surjective.
\end{proof}
\begin{prop}\label{Greenberg Wiles over DVR}
Let $(\mathcal F,\ST) $
be a Selmer structure for $T $ over $R $. Then:
\begin{equation*}
        \rank_{R}\Sel_{\mathcal F} (T) -\rank_{R}\Sel_{\mathcal F ^\dagger} (T^\dagger) =
        \rank_{R} H ^ 0 (L, T) -\rank_{R} H ^ 0 (L, T ^\dagger) +\sum_{\mathsf s\in\ST}\left (\rank_{R}\H ^ 1_{\mathcal F} (L_\mathsf s, T) -\rank_{R}\H ^ 0 (L_\mathsf s, T )\right).
\end{equation*}
\end{prop}
\begin{proof}
 Without loss of generality, we may assume that $$\frac {\H ^ 1 (L_{\mathsf v}, T)} {\H ^ 1_{\mathcal F} (L_{\mathsf v}, T)} $$
 is torsion-free for all $\v$. By Propositions  \ref{Greenberg Wiles}
 and \ref{two different duals}, we then have for each $j $:
 \begin{equation*}
     \begin{split}
         \lg\Sel_{\mathcal F} (T/\pi ^ j)-\lg\Sel_{\mathcal F ^\dagger}(T ^\dagger/\pi ^ j) =\lg\H ^ 0(L, T/\pi ^ j)-\lg\H ^ 0 (L,T ^\dagger/\pi ^ j)\\+\sum_{\v\in\ST}\left (\lg\H ^ 1_\mathcal F (L_\v, T/\pi ^ j)-\lg\H ^ 0 (L_\v, T/\pi ^ j)\right).
     \end{split}
 \end{equation*} 
 Since $\Sel_{\mathcal F} (T) $
 is a  finitely generated $R $-module, it follows from \cite[Lemma 3.7.1]{mazur2004kolyvagin} that
 $$\lg\Sel_{\mathcal F} (T/\pi ^ j) = (\rank_{R}\Sel_{\mathcal F} (T))\cdot \lg R/\pi ^ j +O (1) $$
 as  $j $
 varies, and likewise for $\Sel_{\mathcal F ^\dagger} (T^\dagger) $ and each term on the right-hand side; the proposition follows.
\end{proof}
\section {Bipartite Euler systems}\label{Howard section}
\subsection{Admissible primes}
\subsubsection{}\label{modular forms notation}
Let $f $
be a modular form of weight two, trivial character, and level $N $, and let $\p\subset\O_f $
be a prime ideal of the ring of integers of its field of coefficients. We assume the rational prime $p$ lying under $\p $ is odd,
and write $\O $
for the completion of $\O_f $
at $\p $.
Fix  a Galois-stable $\O $-lattice $T_f$ in the $\p $-adic Galois representation associated to $f $, and let $\overline T_f $
be the residual representation $T_f/\p $; we also write $W_f $
for $T_f\otimes\Q_p/\Z_p $. Also let $K/\Q $
be an imaginary quadratic field. We assume throughout this section that $\overline T_f $
is absolutely irreducible as a $G_K $-module. We will sometimes use the condition:
\begin{equation}
    \tag{sclr}\text {The image of the $G_K $ action on }\overline T_f\text { contains a nonzero scalar.}
\end{equation}
\begin{definition}\label{admissible def}
A nonconstant ultraprime $\mathsf q\in \mathsf M_\Q$
is said to be \textbf{admissible} with sign $\epsilon_{\mathsf q}=\pm 1 $ for $f $
if $\Frob_{\mathsf q} $
has nonzero image in $\Gal (K/\Q) $, $\chi (\Frob_{\mathsf q})\not\equiv 1\pmod p $, and there is a rank-one direct summand $\Fil ^ +_{\q,\epsilon_\q} T_f\subset T_f $
on which $\Frob_{\mathsf q}$ acts as $\chi (\Frob_{\mathsf q})\epsilon_{\mathsf q}$. (Equivalently, $\chi (\Frob_{\mathsf q})\not\equiv 1\pmod p $
and $T_f $
admits a basis of eigenvectors for $\Frob_{\mathsf q} $
with eigenvalues $\epsilon_{\mathsf q}$ and $\chi (\Frob_{\mathsf q})\epsilon_{\mathsf q}$.)
\end{definition}
For example, if $\Frob_{\mathsf q}\in G_\Q $
is a complex conjugation, then $\mathsf q $
is admissible with either choice of $\epsilon_{\mathsf q} $. We abusively write $\mathsf q $
for the unique ultraprime in $\mathsf M_K $
lying over $\mathsf q\in\mathsf M_\Q $, whose Frobenius is $\Frob_{\mathsf q} ^ 2.$
\begin{definition}
If $\mathsf q $
is admissible with sign $\epsilon_{\mathsf q} $
for $f $, then we define the \textbf{ordinary} local condition (with sign $\epsilon_{\mathsf q} $)
as: $$\H ^ 1_{\ord,\epsilon_\q} (K_{\mathsf q}, T_f) =\im\left (\H ^ 1 (K_{\mathsf q}, \Fil ^ +_{\q,\epsilon_\q} T_f)\to\H ^ 1 (K_{\mathsf q}, T_f)\right). $$
The subscript $\epsilon_q $
will often be omitted (from this and future notation) when there is no risk of confusion.
\end{definition}
\subsubsection{}
Note that the ordinary local condition is self-annihilating under the local Tate pairing $$\H^1(K_\q,T_f)\times\H ^ 1 (K_\q, T_f)\to\O $$
induced by (\ref{ imperfect pairing}) and the Weil pairing.  
\subsubsection{}
For any finite set $\ST\subset\mathsf M_K $ such that $T_f$ is unramified outside $\ST$,
and any admissible $\q\not\in\ST $
with sign $\epsilon_\q $, define a localization map \begin{equation}\label{ localization definition}
    \loc_{\q,\epsilon_\q}:\H ^ 1 (K ^\ST/K, T_f)\to\H ^ 1_{\unr} (K_\q, T_f)\to\frac{\H ^ 1_{\unr} (K_\q, T_f)}{\H ^ 1_{\unr} (K_\q, T_f)\intersection\H ^ 1_{\ord,\epsilon_\q}(K_\q, T_f)}\approx \O.
\end{equation} 
Define as well a residue map \begin{equation}\label{ residue definition}
        \partial_{\q,\epsilon_\q}:\H ^ 1 (K, T_f) \to \H ^ 1 (K_\q, T_f)\to\H ^ 1_{\ord,\epsilon_\q} (K_\q, T_f)\to\frac{\H ^ 1_{\ord,\epsilon_\q} (K_\q, T_f)}{\H ^ 1_{\unr} (K_\q, T_f)\intersection\H ^ 1_{\ord,\epsilon_\q}(K_\q, T_f)}\approx \O,
\end{equation}
 where the second map is given by the projection $T_f\twoheadrightarrow (\Frob_\q -\epsilon_\q) T_f \simeq \Fil ^ +_{\q,\epsilon_\q} T_f $. 
 The maps $\loc_{\q,\epsilon_q} $ and $\partial_{\q,\epsilon_\q} $
 may be extended in the obvious way to the patched cohomology for $W_f$ and all  $T_f/\p ^ j $.

\subsection{Euler systems for anticyclotomic twists}
\subsubsection{}
Let $R $
be a complete flat Noetherian local $\O $-algebra with finite residue field, equipped with an anticyclotomic character $\phi: G_K\to R ^\times $
which is trivial modulo the maximal ideal of $R $.
We write $T_\phi$
for the anticyclotomic twist $T_f\otimes_\O R (\phi) $, which is a countably profinite Galois module.
If $\mathsf q $
is admissible with sign $\epsilon_\q $, then $\phi (\Frob_{\mathsf q}^ 2) = 1, $
so $$\H ^ 1 (K_\q, T_\phi) = \H ^ 1 (K_\q, T_f)\otimes_\O R. $$
We extend the ordinary local condition of the previous subsection by linearity to define $\H ^ 1_{\ord,\epsilon_\q} (K_\q, T_\phi)$, and  likewise the maps  $ \loc_{\q,\epsilon_\q}, \partial_{\q,\epsilon_\q} $.

 
 \subsubsection{}\label{where admissible set defined}
Suppose given a finite set $\ST\subset\mathsf M_K $
and a generalized Selmer structure $(\mathcal F,\ST) $
for $T_\phi$. Let $\mathsf N=\mathsf N_\ST $ be the set of pairs $\set{\QT,\epsilon_\QT}$
where  $\mathsf Q\subset\mathsf M_K -\ST$
is a finite set of ultraprimes
and $\epsilon_\QT:\QT\to\set{\pm 1}$ is a function such that $\q$ is admissible with sign $\epsilon_\QT(\q)$ for all $\q\in \QT$. (We will drop the subscript $\ST$ when it is clear from context, or when $\ST$ contains only constant ultraprimes.) Given a pair $\set{\QT,\epsilon_\QT}\in \mathsf N $,  define a generalized Selmer structure $(\mathcal F (\mathsf Q, \epsilon_\QT),\ST\cup \mathsf Q) $
for $T_\phi $
by the local conditions:
\begin{equation}
    \H ^ 1_{\mathcal{F}(\QT,\epsilon_\QT)} (K_\v, T_\phi) = \begin{cases}\H ^ 1_{\mathcal{F} } (K_\v, T_\phi), &\v\not\in\QT\\\H ^ 1_{\ord,\epsilon_\QT(\q)} (K_\q, T_\phi), &\v=\q\in\QT.\end{cases}
\end{equation}

For  $\delta\in\Z/2\Z $, let $\mathsf N  ^\delta\subset\mathsf N  $
be the collection of pairs $\set{\QT,\epsilon_\QT}\in\mathsf N  $
such that $|\QT|\equiv\delta\pmod 2. $
Also, given two pairs $\set{\QT,\epsilon_\QT}\in\mathsf N ^\delta$
and $\set{\QT',\epsilon_{\QT'}}\in\mathsf N ^ {\delta'}  $
such that $\QT\intersection\QT' =\emptyset, $
write $$\set{\QT\QT',\epsilon_{\QT\QT'}}\in\mathsf N ^ {\delta +\delta'} $$
for the pair formed in the obvious way from $\QT\cup \QT' $
and the sign functions $\epsilon_\QT,\epsilon_{\QT'} $. The pair $\set{\emptyset,\emptyset}\in \mathsf N  $
will be abbreviated as 1.
\begin{definition}\label{Euler definition}
A \textbf{bipartite system} $(\kappa,\lambda) $
for $(T_\phi,\mathcal F,\ST) $
of parity $\delta \in\Z/2\Z$
consists of the following data:
\begin {enumerate}
\item for each pair $\set {\QT,\epsilon_\QT}\in\mathsf N ^ {\delta}$, a principal submodule $$(\kappa (\mathsf Q,\epsilon_\QT))\subset\Sel_{\mathcal F(\QT)}(T_\phi);$$
\item for each pair $\set {\QT,\epsilon_\QT}\in\mathsf N ^{\delta+1}$, a principal ideal $$(\lambda (\mathsf Q,\epsilon_\QT))\subset R.$$
\end{enumerate}
A \textbf{bipartite Euler system}
is a bipartite system satisfying the ``reciprocity laws'':
\begin{enumerate}
    \item For each $\set {\QT\q,\epsilon_{\QT\q}}\in\mathsf N ^{\delta +1} $, $$\loc_{\q} ((\kappa (\QT)) = (\lambda (\QT\q))\subset R. $$
    \item For each $\set {\QT\q,\epsilon_{\QT\q}}\in\mathsf N ^{\delta} $, $$\partial_\q ((\kappa (\QT\q)) = (\lambda (\QT))\subset R. $$
\end{enumerate}
We say $(\kappa ,\lambda ) $ is \textbf{nontrivial} if there exists some $\set {\QT,\epsilon_\QT}\in\mathsf N  $
such that either $\lambda (\QT,\epsilon_\QT)\neq 0 $
or $\kappa (\QT,\epsilon_\QT)\neq 0 $
depending on the parity of $|\QT|+\delta$.
\end{definition}
\subsection{Euler systems over discrete valuation rings}
\subsubsection{}
Suppose that $R $
is a discrete valuation ring with uniformizer $\pi $,
and let $W_\phi = T_\phi\otimes\Q_p/\Z_p $. Exactly as in \cite{howard2004heegner}, there is a perfect pairing $T_\phi\times T_\phi\to R(1)$, $G_K $-equivariant up to a twist,
which induces local pairings: \begin{align*}
    \H ^ 1 (K_\v, T_\phi)\times\H ^ 1 (K_{\overline\v}, W_\phi)&\to R\otimes\Q_p/\Z_p,\\
    \H ^ 1 (K_\v, T_\phi/\pi ^ j)\times\H ^ 1 (K_{\overline\v}, W_\phi[\pi^ j]) &\to R/\pi ^ j,\\
    \H ^ 1 (K_\v, T_\phi)\times\H ^ 1(K_{\overline\v}, T_\phi) &\to R.
\end{align*}
Here $\overline\v\in\mathsf M_K $
is the complex conjugate of $\v $; the first two pairings are perfect.
 A Selmer structure $(\mathcal F,\ST) $ for $T_\phi $
induces a Selmer structure for $W_\phi $, denoted the same way, by taking orthogonal complement local conditions.
\begin{definition}
We say $(\mathcal{F},\ST) $
is \textbf{self-dual}
if, for all $\v\in\mathsf M_K $,
$\H ^ 1_{\mathcal  F}(K_\v, T_\phi)$
 and $\H ^ 1_{\mathcal F} (K_{\overline\v}, T_\phi) $
 are exact annihilators under the local pairing.
 \end{definition}
 
\begin{prop}\label{Ray duality}
Suppose that $(\mathcal F,\ST) $
is a self-dual Selmer structure for $T_\phi $.
Then, for each $\set {\QT,\epsilon_\QT}\in\mathsf N $:
\begin{enumerate}
    \item  $(\mathcal{F}(\QT),\ST) $
    is self-dual and
    $$\Sel_{\mathcal F (\QT)} (W_\phi) \approx (R\otimes \Q_p/\Z_p)^{r_\QT}\oplus M_\QT\oplus M_\QT $$
    for some
    torsion $R $-module $M_\QT $
   and an integer $r _\QT $.
    \item   $r_\QT =\rank_R\Sel_{\mathcal F (\QT)} (T_\phi)$.
    \item For any $\set{\QT\q,\epsilon_{\QT\q}}\in\mathsf N $, 
    one of the following holds:
    \begin{enumerate}
        \item $\loc_\q (\Sel_{\mathcal F (\QT)} (T_\phi)) = 0, $ $\partial_\q (\Sel_{\mathcal F (\QT\q)} (T_\phi))\neq 0,$
        $r_{\QT\q} = r_\QT +1, $
        and 
        there exists an exact sequence of $R$-modules: $$0\to M_{\QT\q}\to M_\QT \to \loc_\q (\Sel_{\mathcal F (\QT)}) (W_\phi)\to 0. $$
        Moreover, $$\lg \loc_\q (\Sel_{\mathcal F (\QT)}) (W_\phi)=\lg\coker\partial_\q (\Sel_{\mathcal F (\QT\q)}  (T_\phi)). $$
        \item $\loc_\q (\Sel_{\mathcal F (\QT)} (T_\phi)) \neq 0, $ $\partial_\q (\Sel_{\mathcal F (\QT\q)} (T_\phi))=0,$
        $r_{\QT\q} = r_\QT -1, $
        and there exists an exact sequence of $R$-modules:
        $$0\to M_\QT\to M_{\QT\q}\to \partial_\q(\Sel_{\mathcal F(\QT\q)}(W_\phi)) \to 0.$$
        Moreover, $$\lg \partial_\q(\Sel_{\mathcal F(\QT\q)}(W_\phi)) = \lg\coker\loc_\q (\Sel_{\mathcal F (\QT)} (T_\phi)) .$$
    \end{enumerate}
\end{enumerate}
\end{prop}
\begin{proof}
\begin{enumerate}
    \item The self-duality claim is clear since $H^1_{\ord}(K_\q, T_f)$ is self-dual.
    Now, for all $j\geq 0,$ \begin{equation}
        \Sel_{\mathcal F (\QT)} (W_\phi[\pi ^ j])\simeq\Sel_{\mathcal F (\QT)} (W_\phi) [\pi^ j] \tag{$\ast$}
    \end{equation} 
   by Lemma \ref{LES} and the definition of the  induced Selmer structure on $W_\phi[\pi^j]$. (Note $H^0(K, \overline T_f) = 0$ since we have assumed $\overline T_f$ is an absolutely irreducible $G_K$-module.) Since $\Sel_{\mathcal F (\QT)} (W_\phi) $
   is co-finitely generated, we may conclude by \cite[Theorem 1.4.2]{howard2004heegner} (or its proof).
   \item 
   As explained in
\cite{howard2004heegner,howard2006bipartite}, the cohomological pairings deduced from $T_\phi\times T_\phi\to R(1)$ behave ``exactly like'' the Tate pairing; in particular, by the self-duality of the local conditions and Proposition \ref{two different duals}, $\Sel_{\mathcal F (\QT)} (T_\phi/\pi ^ j) =\Sel_{\mathcal F (\QT)} (W_\phi [\pi ^ j]) $
and the result follows as in the proof of Proposition \ref{Greenberg Wiles over DVR}.
   \item 
   Consider the Selmer structures $\mathcal F ^\q (\QT)=\mathcal{F}(\QT) +\mathcal F (\QT\q) $ and $\mathcal F_\q (\QT) = \mathcal{F}(\QT) \intersection\mathcal F (\QT\q) $.
   By Proposition \ref{Greenberg Wiles over DVR},
   $$\rank_{R}\Sel_{\mathcal F ^\q (\QT)} (T_\phi) =\rank_{R}\Sel_{\mathcal F_\q (\QT)} (T_\phi) +1. $$
   Moreover, because $\mathcal F (\QT) $
   is self-dual, Proposition \ref{ dualityover DVR}
   implies that the image of
   \begin{align*}\frac{\Sel_{\mathcal F ^\q (\QT)} (T_\phi)}{\Sel_{\mathcal F_\q (\QT)}(T_\phi)}&\hookrightarrow\frac {\H ^ 1_{\mathcal F^\q(\QT)} (K_\q, T_\phi)} {\H ^ 1_{\mathcal F_\q(\QT)} (K_\q, T_\phi) } 
= \frac {\H ^ 1_{\unr} (K_\q, T_\phi) }{\H ^ 1_{\mathcal F_\q(\QT)} (K_\q, T_\phi)}\oplus \frac{\H ^ 1_{\ord} (K_\q, T_\phi)} {\H ^ 1_{\mathcal F_\q(\QT)} (K_\q, T_\phi)}\approx R^2\end{align*}
   is self-annihilating under the induced local pairing, hence is contained either in the ordinary
   or unramified part. 
   
   For the relation between $M_\QT$ and $M_{\QT\q}$, we suppose we are in case (a), because the two arguments are identical. Using the perfect pairing between $W_\phi$ and $T_\phi$, we see by Proposition \ref{Perfect pairing} that $\loc_\q (\Sel_{\mathcal F(\QT)}(W_\phi)) \oplus \partial_\q (\Sel_{\mathcal F(\QT)}(W_\phi))$ is the exact annihilator of $\partial_\q (\Sel_{\mathcal F(\QT\q)}(T_\phi))$ under the perfect induced local pairing $$\frac{\H^1_{\mathcal F^\q(\QT)}(K_\q, T_\phi)}{\H ^ 1_{\mathcal F_\q(\QT)} (K_\q, T_\phi)}\times \frac{\H^1_{\mathcal{F}^\q (K_\q,\QT)} (W_\phi)}{\H ^ 1_{\mathcal F_\q(\QT)} (K_\q, W_\phi)}\to R\otimes\Q_p/\Z_p.$$
   This implies that $\partial_\q (\Sel_{\mathcal F (\QT\q)}) (W_\phi)$ is divisible and $$\lg \loc_\q (\Sel_{\mathcal F (\QT)}) (W_\phi)=\lg\coker\partial_\q (\Sel_{\mathcal F (\QT\q)}  (T_\phi)).$$ Furthermore, the direct sum decomposition of $\Sel_{\mathcal F (\QT\q)} (W_\phi) $
   may be chosen so that $\partial_\q(M_{\QT\q}\oplus M_{\QT\q}) = 0,$ and in particular $$M_{\QT\q}\oplus M_{\QT\q}\simeq (M_\QT\oplus M_\QT) \intersection \Sel_{\mathcal F_\q(\QT)}\subset (\Sel_{\mathcal F (\QT)}) (W_\phi).$$
   Since $\loc_\q (M_\QT\oplus M_\QT)$ must generate the image of $\loc_\q(\Sel_{\mathcal F (\QT)}) (W_\phi)$, the desired exact sequence follows.
\end{enumerate}
\end{proof}
The following result will allow us to control the alternative in Proposition \ref{Ray duality}(3).
\begin{thm}\label{localization nontrivial}
Let $c\in\H ^ 1 (K^ {\mathsf T}/K, T_\phi) $
 be any nonzero element, where $\mathsf T\supset\ST $
 is a finite set. Then there are infinitely many admissible ultraprimes $\q \not\in\mathsf T$, with associated signs $\epsilon_\q $, such that $\loc_\q c\neq 0. $
\end{thm}
The proof is via a series of lemmas.
\begin{lemma}\label{bounding down cohomology}
There is an integer $j $ 
such that, for all $n\geq 0, $ $$\pi ^ j H ^ 1 (K (T_\phi)/K, T_\phi/\pi ^ n) = 0. $$ If (\ref{contains scalar}) holds, then we may take $j = 0. $
\end{lemma}
\begin{proof}
Let $G =\Gal (K (T_\phi)/K) $, and let $Z\subset G $
be its center; since $T_f $
is absolutely irreducible over $K $, $Z$
acts on $T_\phi $
by scalars. We claim:
\begin{equation}\label{nontrivial scalar}
    Z\neq\set{1}.
\end{equation}
Assuming (\ref{nontrivial scalar}), the lemma follows from the inflation-restriction exact sequence
$$H^ 1 (G/Z, H ^ 0 (Z, T_\phi/\pi ^ n))\hookrightarrow H ^ 1 (G, T_\phi/\pi ^ n)\to H ^ 1 (Z, T_\phi/\pi ^ n). $$ 
Let us now prove (\ref{nontrivial scalar}).
Let $G' =\Gal (K (T_f)/K) $, and let $L/K $
be the Galois subfield of $K (T_f) $
cut out by the center $Z' =Z (G')\subset G' $. By a result of Momose \cite{ribet1985adic}, $Z' $
is nontrivial. Let $E/K$ be the Galois extension determined by the kernel of $\phi $; then  it suffices to show that $EL/L $
and $K (T_f)/L $
are linearly disjoint. Both $EL$ and $K (T_f) $
are Galois over $\Q $, so $G_\Q $
acts on $\Gal (EL/L) $
and $\Gal (K (T_f)/L) $
by conjugation. If $\tau\in G_\Q$
is a complex conjugation, then $\tau $
acts trivially on $\Gal (K (T_f)/L) $
but nontrivially on $\Gal (EL/L)$, so the two groups have no nontrivial common quotient compatible with the $G_\Q$-action; hence $EL\intersection K (T_f) = L $.
\end{proof}

\begin{lemma}\label{just one localization}
Suppose given a  cocycle $$c\in H ^ 1 (K,  T_\phi/\pi^n) $$
such that $\pi ^ jc\neq 0, $
where $j $ is as in Lemma \ref{bounding down cohomology}. Then, for any integer $N\geq n$, there exists a sign $\epsilon=\pm 1 $ and infinitely many rational primes $q$
such that:
\begin{enumerate}
\item $q $
is inert in $K $ and unramified in the splitting field $\Q (T_f,c) $.
    \item $\Frob_q\in\Gal (\Q (T_f)/\Q) $
    has distinct eigenvalues $\pm 1$ on $T_f\otimes R/\pi ^ N $ (where $R $
    has trivial Galois action).
    \item For any cocycle representative, $c (\Frob_q ^ 2) $
    has nonzero component in the $\epsilon $
    eigenspace for $\Frob_q $.
\end{enumerate}
\end{lemma}
\begin{proof}
Abbreviate $L = K (T_\phi/\pi^N)$, and
let $\xphi\in\Hom_{G_K} (G_{L}, T_\phi/\pi ^ n)$
be the image of $c $
under restriction; by hypothesis $\xphi\neq 0.$
Without loss of generality, we may suppose that the image of $\xphi $
is contained in $T_\phi/\pi ^ n [\pi] \simeq T_\phi/\pi$, which, since $\phi $
is residually trivial,  is an extension of scalars $\overline T_f \otimes_{\O/\p} k$.
Now, $$\Hom_{G_K} (G_{L},\overline T_f\otimes k) $$
has a natural action of $\Gal (K/\Q) $, and we may assume without loss of generality that $\xphi$
lies in the $\epsilon$
eigenspace for some $\epsilon\in\set{\pm 1}.$ 
Fix a complex conjugation $\tau\in G_\Q $. Since $\overline T_f $
is absolutely irreducible over $G_K $, there exists $g\in G_{L} $
such that $\xphi (g) $
has nonzero component in the $\epsilon$
eigenspace of $\tau$. Then $$\xphi (\tau g\tau g) =\epsilon\tau\xphi (g) +\xphi (g) $$
has nonzero component in the $\epsilon $
eigenspace as well. Any $q $
with Frobenius $\tau g $
in $L (\xphi) $
satisfies the desired conditions.
\end{proof}
\begin{proof}[Proof of Theorem \ref{localization nontrivial}]
Since $H ^ 0 (K, T_\phi/\pi) = 0, $
Lemma \ref{LES} implies that
$$\H ^ 1 (K^ {\mathsf T}/K, T_\phi)[\pi] = 0. $$
Thus there exists some $n $
such that the image $\overline c$ of $c$ in $\H ^ 1 (K^ {\mathsf T}/K, T_\phi/\pi ^ n) $
satisfies $\pi ^ j\overline c\neq 0, $
for some $j $
as in Lemma \ref{bounding down cohomology}. By definition, $\overline c $
is represented by a sequence of classes $c_m\in H ^ 1 (K^ {T_m}/K, T_\phi/\pi ^ n) $
such that $\pi ^ jc_m\neq 0 $
for $\mathfrak F$-many $m $, where $\set{T_m}_{m\in\N} $
represents $\mathsf T $.
For each $m $, apply Lemma \ref{just one localization} (with $N=m$)
to obtain a prime $q_m \not\in T_m$
and a sign $\epsilon_m $. If $\q\in\mathsf M_\Q $
is the equivalence class of the sequence $\set {q_m}_{m\in\N} $,
and $\epsilon \in\mathcal U (\set{\pm 1}_{m\in\N})\simeq\set {\pm 1}$
is the equivalence class of the sequence $\set {\epsilon_m}_{m\in\N}, $
then the pair $\set {\q,\epsilon} $
has
the desired properties. Since there are infinitely many choices for each $q_m $, there are also infinitely many choices for $\q $.
\end{proof}
\begin{cor}\label{the Selmer annihilation corollary}
For any $\set{\QT,\epsilon_\QT}\in \mathsf N$,  there exists some $\set{\QT\QT',\epsilon_{\QT\QT'}}\in\mathsf N $
such that $r_{\QT\QT'} = 0. $
\end{cor}
\begin{proof}
 This is an obvious induction argument using Theorem \ref{localization nontrivial} and Proposition \ref{Ray duality}.
\end{proof}
Combining Proposition \ref{Ray duality} and Theorem \ref{localization nontrivial} allows us to prove the main result of this subsection.
\begin{thm}\label{DVR Euler system bound}
Suppose that $(\mathcal F,\ST) $
is self-dual and that $(\kappa,\lambda) $
is a nontrivial bipartite Euler system with sign $\delta $ for $(T_\phi,\mathcal F,\ST) $.
Then there exists a  nonzero fractional ideal $I $
    of $R $
    such that:
    \begin{enumerate}
        \item For all $\set{\QT,\epsilon_\QT}\in\mathsf N ^\delta $, $r_\QT $
    is odd, $r_\QT = 1 $
    if and only if $\kappa (\QT) \neq 0,$
    and in that case $$\Char_R (M_\QT) \cdot I =\Char_R \left(\frac{\Sel_{\mathcal F (\QT)} (T_\phi)} {(\kappa (\QT))}\right). $$
    \item For all $\set{\QT,\epsilon_\QT}\in\mathsf N ^ {\delta +1}$, $r_\QT $
    is even, $r_\QT = 0 $
    if and only if $\lambda (\QT)\neq 0, $
    and in that case $$\Char_R (M_\QT)\cdot I = (\lambda (\QT)). $$
\end{enumerate}
In particular, $$\delta =\rank_R\Sel_{\mathcal F} (T_\phi) +1\pmod 2. $$
\end{thm}
\begin{proof}
The proof will be in several steps.
\begin{step}
If $\lambda (\QT)\neq 0
$
for some $\set{\QT,\epsilon_\QT}\in\mathsf N ^ {\delta +1} $, then $r_\QT = 0. $
\end{step}
\begin{proof}
If $0\neq c\in\Sel_{\mathcal F (\QT)} (T_\phi) $, then by Theorem \ref{localization nontrivial}, there exists an admissible ultraprime $\q $
with sign $\epsilon_\q $
such that $\loc_\q c\neq 0. $ By Proposition \ref{Ray duality}, $\partial_\q (\kappa (\QT\q)) = 0, $
which contradicts the reciprocity laws.
\end{proof}
\begin {step} If $\kappa (\QT)\neq 0 $
for some $\set {\QT,\epsilon_\QT}\in\mathsf N ^\delta $, then $r_\QT = 1. $
\end {step}
\begin {proof}
Choose an admissible ultraprime $\q $
with sign $\epsilon_\q $
such that $\loc_\q\kappa (\QT)\neq 0. $
Then by the reciprocity laws, $\lambda (\QT\q)\neq 0, $
so by Step 1 $r_{\QT\q} = 0. $ Proposition \ref{Ray duality} implies $r_{\QT} = 1. $\end {proof}
\begin {step}
 For all $\set {\QT,\epsilon_\QT}\in\mathsf N $, $r_\QT\equiv\delta + |\QT | +1\pmod 2. $
\end{step}
\begin {proof}
If $\set{\QT\QT',\epsilon_{\QT\QT'}}\in\mathsf N $, then by Proposition \ref{Ray duality} $$r_{\QT}-r_{\QT\QT'}\equiv |\QT' |\pmod 2. $$
So Steps 1 and 2 imply Step 3.
\end {proof}
\begin {step} Suppose  $r_\QT = 0 $
for some $\set {\QT,\epsilon_\QT} $. Then, for all admissible ultraprimes $\q \not\in\QT\cup \ST$
with sign $\epsilon_\q $, $r_{\QT\q} = 1 $
and $$\Char_R (M_{\QT\q})\cdot(\lambda (\QT)) =\Char_R (M_\QT)\cdot\Char_R\left (\frac{\Sel_{\mathcal F (\QT\q)}(T_\phi)}{(\kappa (\QT\q))}\right). $$
\end {step}
\begin {proof}
By Step 3, $\lambda(\QT) $
and $\kappa (\QT\q) $
are well-defined. Then Step 4 follows from Proposition \ref{Ray duality}, since $$\Char_R\left (\frac{\Sel_{\mathcal F (\QT\q)}(T_\phi)}{(\kappa (\QT\q))}\right)\cdot (\partial_\q(\Sel_{\mathcal F (\QT\q)} (T_\phi)) = (\lambda (\QT))\subset R. $$
\end{proof}
The exact same reasoning implies:
\begin{step}
Suppose that $r_{\QT} = 1$
and  $\q\not\in\QT\cup\ST $
is an admissible ultraprime with sign $\epsilon_\q $
such that $r_{\QT\q}= 0.$
Then $$\Char_R (M_{\QT\q})\cdot\left (\frac {\Sel_{\mathcal F (\QT)} (T_\phi)} {(\kappa (\QT))}\right)=\Char_R (M_\QT)\cdot (\lambda (\QT\q)). $$
\end{step}
Now consider the graph $\xchi $ \cite{howard2006bipartite} whose vertices are the elements of $\mathsf N $, and where the edges are between vertices of the form $\set{\QT,\epsilon_\QT}$
and $\set {\QT\q,\epsilon_{\QT\q}} $, for some admissible ultraprime $\q $
with sign $\epsilon_\q $. We say $\set {\QT,\epsilon_QT} $
is a \textbf{core vertex} if $r_\QT\leq 1. $
The \textbf{core subgraph} $\xchi _0 $
of $\xchi $
is the full subgraph on core vertices.
Applying Steps 4 and 5, it suffices to show that $\xchi_0 $
is  path-connected to complete the proof of the theorem.
\begin{step}
If $v=\set {\QT,\epsilon_\QT} $
and $v'=\set {\QT\QT',\epsilon_{\QT\QT'}} $
are core vertices, then they are connected by a path in $\xchi_0$.
\end{step}
\begin{proof}
 We proceed by induction on $|\QT' | $, where the base case is trivial. If $r_{\QT\QT'/\q}\leq 1 $
 for any $\q\in\QT' $, then we may apply the inductive hypothesis, so assume otherwise. By Proposition \ref{Ray duality}, 
 $r_{\QT\QT'} = 1 $ and
 $\partial_\q (\Sel_{\mathcal F (\QT\QT')} (T_\phi)) = 0 $
 for all $\q\in\QT' $. Hence $$\Sel_{\mathcal F (\QT\QT')} (T_\phi)\subset\Sel_{\mathcal F (\QT)}(T_\phi).$$ Then, by Theorem \ref{localization nontrivial} and Proposition \ref{Ray duality}, there exists an admissible ultraprime $\q \not\in\QT\cup\QT'\cup\ST$
 with sign $\epsilon_\q $
 such that $r_{\QT\q} = r_{\QT\QT'\q} = 0. $
 If $\q'\in\QT' $
 is any factor, then $\set{\QT\QT'\q/\q',\epsilon_{\QT\QT'\q/\q'}}\in\mathsf N $
 is a core vertex, which is connected to $v' $
 in $\xchi_0 $. By the inductive hypothesis, $\set{\QT\QT'\q/\q',\epsilon_{\QT\QT'\q/\q'}}$ is also connected to the core vertex $\set{\QT\q,\epsilon_{\QT\q}} $, hence to $v$,
 by a path in $\xchi_0. $
 This completes the inductive step.
\end{proof}
\begin{step}
If $v=\set{\QT,\epsilon_\QT}$ is a core vertex and $\mathsf T\subset\mathsf M_\Q $
is any finite set, then there exists a core vertex $v'=\set {\QT',\epsilon_{\QT'}} $
such that $v $
and $v' $
are connected by a path in $\xchi_0 $
and $\QT'\intersection\mathsf T =\emptyset $.
\end{step}
\begin{proof}
By iterating, it suffices to assume that $\QT\intersection\mathsf T $
consists of exactly one ultraprime $\q\in\QT $. If $r_{\QT/\q}\leq 1, $
then the conclusion is obvious, so suppose otherwise. As in the proof of Step 6, choose an admissible ultraprime $\q'\not\in\QT\cup\ST\cup \mathsf T $
 with associated sign $\epsilon_{\q'} $
such that $r_{\QT\q'} = 0, $
which implies $r_{\QT\q'/\q} = 1$. The core vertex $v'=\set {\QT\q'/\q,\epsilon_{\QT\q'/\q}} $
has the desired properties.
\end{proof}
Finally, we have:
\begin{step}
The core subgraph $\xchi_0 $
is path-connected.
\end{step}
\begin{proof}
Let $\set{\QT_1,\epsilon_{\QT_1}}$
and $\set {\QT_2,\epsilon_{\QT_2}} $
be two core vertices. Without loss of generality, by Step 7, we may assume $\QT_1\intersection\QT_2 =\emptyset $. (This step is necessary because the sign functions $\epsilon_{\QT_1} $
and $\epsilon_{\QT_2} $
need not agree on $\QT_1\intersection\QT_2. $) Consider $\set {\QT_1\QT_2,\epsilon_{\QT_1\QT_2}}\in\mathsf N $. This may not be a core vertex, but, by repeatedly applying Theorem \ref{localization nontrivial} and Proposition \ref{Ray duality}, there exists $\set {\QT_3,\epsilon_{\QT_3}}\in\mathsf N $
such that $\set {\QT_1\QT_2\QT_3,\epsilon_{\QT_1\QT_2\QT_3}} $
is a core vertex. We may then conclude by Step 6.
\end{proof}
\end{proof}
\begin{prop}\label{ uniformDVR bound}
Under the hypotheses of Theorem \ref{DVR Euler system bound}, there exists a constant $C$ depending on $|\ST |$, $T_f$, and the ramification index of $R/\O $, but not on $\phi$, such that $I\pi^C\subset R$. If (\ref{contains scalar}) holds,
then we may take $C = 0. $
\end{prop}
\begin{proof}
 By Theorem \ref{DVR Euler system bound}, it suffices to show that there exists a constant with the desired dependencies and a pair $\set {\QT,\epsilon_\QT}\in\mathsf N ^\delta $
 such that $\pi^C\in\Char_R(M_\QT)$. We first note that the constant $j$ in Lemma \ref{just one localization} depends only on $T_f$ and the ramification index of $R/\O$, and can be taken to be 0 under (\ref{contains scalar}).
 
 Moreover, if $k$ is the residue field of $R$, then $d=\dim_k\H^1(K^{\ST}/K,W_\phi[\pi])$
 is also bounded with the desired uniformity. We now construct a sequence $\set{\QT_i,\epsilon_{\QT_i}}$ recursively (starting from $\QT_1 =1 $) by the following rules:
 \begin{itemize}
    \item If $r_{\QT_i} >0, $
    then choose any  $\q_{i+1}\not\in \QT_i$ with sign $\epsilon_{\q_{i+1}}$ such that $$\lg\coker(\loc_\q \Sel_{\mathcal F (\QT_i)} (T_\phi))\leq j.$$
     \item If $r_{\QT_i}=0$ and the exponent of  $\Sel_{\mathcal F (\QT_i)} (W_\phi)\neq 0 $
     is $n_i> i\cdot j$,
     then choose any $\q_{i+1}\not\in \QT_i$ with sign $\epsilon_{\q_{i+1}}$ such that the exponent of $\loc_\q (\Sel_{\mathcal F (\QT_i)} (W_\phi))$ is at least $n_i - j$.
 \end{itemize}
 These choices are possible by Lemma \ref{just one localization}.
In either of the above two cases,
set $$\set{\QT_{i +1} ,\epsilon_{\QT_{i+ 1}}} =\set {\QT_i\q_{i +1},\epsilon_{\QT_i\q_{i+1}}}; $$
     if neither holds, then end the construction.
 For each $i$, let $r_{\QT_i}'$ be the minimal number of generators of the $R $-module $\pi ^ {i\cdot j} M_{\QT_i}. $ 
 In the first case of the construction, $r_{\QT_{i+1}}'\leq r_{\QT_i}'$; in the second case, $r_{\QT_{i+1}}' < r_{\QT_i}'$ (by Proposition \ref{Ray duality}(3b,a) respectively). 
 After $r_{1}\leq d$ steps, we alternate between the two cases of the construction, taking at most $2r_{1}'\leq 2d$ more steps. 
Hence for some $i\leq 3d $, $r'_{\QT_i} = 0$ and $r_{\QT_i} = 0$, and the construction halts. 
 For this $i$, $$\lg M_{\QT_i} \leq ij \dim_k \Sel_{\mathcal F(\QT_i)} (W_\phi)[\pi] \leq 3d j (d + 3d),$$
  the last inequality by the reasoning of \cite[Corollary 2.2.10]{howard2006bipartite}. (Less precisely, we could deduce the bound $3dj(d + 6d)$ directly from Proposition \ref{Ray duality}(3).)
  
  Since $d$ and $j$ have bounds of the desired sort, the claim follows.
\end{proof}
\subsection{Euler systems over $\Lambda $}
Let $\Lambda $
be the anticyclotomic Iwasawa algebra $\O\llbracket T\rrbracket $ with canonical character $$\Psi: G_K\to\Lambda ^\times. $$ For each height-one prime $\P\subset\Lambda $, let $S_\P $
be the integral closure of $\Lambda/\P $
in its field of fractions, so that $\Psi $
induces a character $G_K\to\Lambda ^\times\to S_\P ^\times $. We write $T_\P $
for the twist $T\otimes_\O S_\P (\Psi) $
and $\mathbf T$ for the interpolated twist $T\otimes_\O\Lambda (\Psi) $. Also let $\mathbf W =\mathbf T ^\ast $
be the Cartier dual with $\Lambda$
action twisted by the canonical involution $\iota $, so that for each $\P $
there is a natural map $$W_\P\to\mathbf W $$
of $\Lambda[G_K] $-modules (see, e.g., \cite{howard2004heegner}). The following definition is motivated by \cite[Lemma 5.3.13]{mazur2004kolyvagin} and its applications in \cite{howard2004heegner,howard2006bipartite}.
\begin{definition}
An \textbf {interpolated self-dual} Selmer structure $$(\ST,\mathcal F_\Lambda,\mathcal F_\P,\Sigma_\Lambda) $$ for $\mathbf T $
consists of the following data:
\begin{itemize}
\item A finite set $\ST\subset\mathsf M_K $.
    \item For each height-one prime $\P\subset\Lambda $, a  self-dual Selmer structure $(\mathcal F_\P,\ST) $
    for $T_\P $.
    \item A finite set $\Sigma_\Lambda $
    of height-one primes $\P\subset\Lambda $.
    \item A Selmer structure $(\mathcal F_\Lambda,\ST) $
    for $\mathbf T $
    such that, for all $\v\in\mathsf M_K $ and all $\P\subset\Lambda $, there are well-defined maps induced in the obvious way: \begin{equation}\label{Lambda interpolating Selmer}
    \begin{split}
        \H ^ 1_{\mathcal F_\Lambda} (K_\v,\mathbf T/\P) &\to \H ^ 1_{\mathcal F_\P} (K_\v, T_\P),\\
        \H^ 1_{\mathcal F_\P} (K_\v, W_\P) &\to\H ^ 1_{\mathcal F_\Lambda} (K_\v,\mathbf W [\P]).
    \end{split}      
    \end{equation}
    Moreover, for all $\P\not\in\Sigma_\Lambda $, the maps (\ref{Lambda interpolating Selmer}) have finite kernel and cokernel with order bounded by a constant depending only on $[S_\P:\Lambda/\P] $
    as $\P $ varies.
\end{itemize}
\end{definition}
\begin{prop}[\cite{mazur2004kolyvagin}, Proposition 5.3.14]\label{interpolating Selmer groups}
Suppose $(\ST,\mathcal F_\P,\mathcal F_\Lambda, \Sigma_\Lambda) $
is an interpolated self-dual Selmer structure for $\mathbf T $. Then for all $\P\subset\Lambda $, there are well-defined maps induced in the obvious way:
\begin{align*}
    \Sel_{\mathcal F_\Lambda} (\mathbf{T})/\P &\to\Sel_{\mathcal F_\P} (T_\p)\\
    \Sel_{\mathcal F_\P} (W_\P) &\to\Sel_{\mathcal F_\Lambda} (\mathbf W) [\P].
\end{align*}
For all $\P\not\in\Sigma_\Lambda $, these maps have finite kernel and cokernel with a bound depending on $\mathcal F$ and on $[S_\P:\Lambda/\P]$, but not on $\P$ itself.
\end{prop}
\subsubsection{}
Recall that, for any finitely generated $\Lambda $-module $M$, there exists a unique $\Lambda $-module $N$ of the form $\Lambda ^r \oplus\bigoplus\Lambda/\P_i ^ {e_i} $ such that $M $
admits a map to $N $
with finite kernel and cokernel; we denote this relationship by $M\sim N $. The characteristic ideal $\Char_\Lambda (M) $
is zero if $r\geq 1, $
and equal to $\prod\P_i ^ {e_i} $
otherwise.   The following easy lemma is implicit in \cite[p. 66]{mazur2004kolyvagin}.
\begin{lemma}\label{approaching in lambda}
Let $\P\subset\Lambda $ be a height-one prime. Then there exists an integer $d$ and a sequence
 of height-one primes $\P_m $
 such that, for all finitely generated torsion $\Lambda $-modules $M $, $$\lg_{\O} (M/{\P_m}) = md \ord_\P\Char_\Lambda (M) + O (1) $$
 as $m $
 varies (holding $M $ fixed). Moreover $[S_{\P_m}:\Lambda/\P_m]$ is constant for large enough $m$, and if $\P\neq (\p), $
 then the rings $\Lambda/\P_m $
 are abstractly isomorphic.
\end{lemma}
\begin{proof}
If $\P\neq (\p) $
is generated by a distinguished polynomial $f\in\Lambda $, and $\pi $
is a uniformizer for $\O $, then we may take $\P_m = f+\pi^m$ (for sufficiently large $m $)
and $d = [S_\P:\O] $. If $\P = (\p)= (\pi) $, then we may take $\P_m = T ^ m +\pi $
and $d = 1. $
\end{proof}

\begin{prop}\label{self dual upshot}
Suppose that $(\ST,\mathcal F_\P,\mathcal F_\Lambda, \Sigma_\Lambda) $
is an interpolated self-dual Selmer structure for $\mathbf T $. Then for all $\set {\mathsf Q,\epsilon_\QT}\in\mathsf N  $:
\begin{enumerate}
    \item $(\ST\cup \QT,\mathcal F_\P (\QT),\mathcal F_\Lambda(\QT),\Sigma_\Lambda)$
    is an interpolated self-dual Selmer structure for $\mathbf T $
    and $$\Sel_{\mathcal F_\Lambda(\QT)}(\mathbf W) ^\vee\sim\Lambda ^ {r_\QT}\oplus M_\QT\oplus M_\QT $$
    for some torsion $\Lambda $-module $M_\QT $
    and an integer $r_\QT $.
    \item $r_\QT =\rank_\Lambda\Sel_{\mathcal F_\Lambda(\QT)} (\mathbf T) $.
\end{enumerate}
\end{prop}
\begin{proof}
 At places $\q\in\QT $,
 $$\H ^ 1_{\mathcal F_\Lambda (\QT)} (K_\q,\mathbf T) =\H ^ 1_{\ord} (K_\q, T_f)\otimes\Lambda$$ 
 and  $$\H ^ 1_{\mathcal F_\P(\QT)} (K_\q,\mathbf T) =\H ^ 1_{\ord} (K_\q, T_f)\otimes S_\P ,$$ so we clearly have local maps with kernel and cokernels
 bounded as desired (and similarly for $\mathbf W_f$ and $W_\P$); so indeed
$(\ST\cup \QT,\mathcal F_\P (\QT),\mathcal F_\Lambda(\QT),\Sigma_\Lambda) $ is an interpolated self-dual Selmer structure. The rest of the claims are deduced from Proposition \ref{interpolating Selmer groups} and Proposition \ref{Ray duality}(1,2) exactly as in \cite[Theorem 2.2.10]{howard2004heegner}.
\end{proof}

\begin{thm}\label{ lambda Euler system bound}
Suppose that $(\ST,\mathcal F_\P,\mathcal F_\Lambda, \Sigma_\Lambda) $
is an interpolated self-dual Selmer structure for $\mathbf T $ and  $\set {\boldsymbol\kappa,\boldsymbol\lambda} $
is a nontrivial bipartite Euler system with parity $\delta $ for the triple $(\mathbf T,\mathcal F_\Lambda,\ST) $. Then there exists a nonzero fractional ideal $I\subset\Lambda\otimes \Q_p $
such that:
\begin{enumerate}
 \item For all $\set{\QT,\epsilon_\QT}\in\mathsf N ^\delta $, $r_\QT $
    is odd, $r_\QT = 1 $
    if and only if $\boldsymbol\kappa (\QT) \neq 0,$
    and in that case $$\Char_\Lambda (M_\QT) \cdot I =\Char_\Lambda \left(\frac{\Sel_{\mathcal F_\Lambda (\QT)} (\mathbf T)} {(\boldsymbol\kappa (\QT))}\right). $$
    \item For all $\set{\QT,\epsilon_\QT}\in\mathsf N ^ {\delta +1}$, $r_\QT $
    is even, $r_\QT = 0 $
    if and only if $\boldsymbol\lambda (\QT)\neq 0, $
    and in that case $$\Char_\Lambda (M_\QT)\cdot I = (\boldsymbol\lambda (\QT)). $$
\end{enumerate}
In particular, $$\delta =\rank_R\Sel_{\mathcal F} (\mathbf T) +1\pmod 2. $$
If (\ref{contains scalar}) holds, then $I\subset \Lambda$. 
\end{thm}
\begin{proof}
 Let $\P\subset\Lambda $
 be any height-one prime; via the natural maps $\Sel_{\mathcal F_\Lambda}(\mathbf T)\to\Sel_{\mathcal F_\P} (T_\P) $
 and $\Lambda\to S_\P $, the Euler system $(\boldsymbol\kappa,\boldsymbol\lambda) $
 defines an Euler system $(\kappa_\P,\lambda_\P) $ of  parity $\delta$
 for the triple $(T_\P,\mathcal F_\P,\ST) $. In particular, Theorem \ref{DVR Euler system bound} applies.
 
 If $\boldsymbol\kappa (\QT)\neq 0, $
 then by Proposition \ref{interpolating Selmer groups} $\kappa_\P (\QT)\neq 0 $
 for all but finitely many $\P $ (since $\Sel_{\mathcal F_\Lambda (\QT)} (\mathbf T) $
 is  torsion-free), and similarly for $\boldsymbol\lambda (\QT) $.
 Because $$\rank_\Lambda\Sel_{\mathcal F_\Lambda} (\mathbf T) \leq\rank_{S_\P}\Sel_{\mathcal F_\P} (T_\P) $$
with equality for all but finitely many $\P $,
 the claims about $r_\QT $
 follow from Theorem \ref{DVR Euler system bound}.
 
 For any $\P $
 and $\set{\QT,\epsilon_{\QT}}\in\mathsf N ^ {\delta +1} $
such that $\boldsymbol{\lambda}(\QT)\neq 0$, by Proposition \ref{self dual upshot} and Lemma \ref{approaching in lambda}
we have
\begin{align*}
   e_\P(\QT)&\coloneqq\ord_\P (\boldsymbol\lambda (\QT)) -\ord_\P\Char_\Lambda (M_{\QT})\\ &=\lim_{m\to\infty}\frac {\lg_\O (S_{\P_m}/\lambda_{\P_m}(\QT)) -\lg_\O M_{\QT,\P_m}} {m d}.
\end{align*}
Applying Theorem \ref{DVR Euler system bound}, this quantity does not depend on $\set{\QT,\epsilon_\QT} $ (as long as $\boldsymbol\lambda (\QT)\neq 0 $); it is also clearly zero for almost all $\P$, so that $\prod_\P \P^{e_\P}$ defines a fractional ideal $I$ of $\Lambda $ satisfying (2). 
The same calculation shows that
$I$ satisfies (1) as well, and the integrality properties follow from Proposition \ref{ uniformDVR bound}. 
\end{proof}
\section{Geometry of modular Jacobians}\label{geometry sec}
\subsection{Multiplicity one}
\subsubsection{}
Let $N_1 $
and $N_2 $
be coprime positive integers. 
Consider the Hecke algebra $\T =\T_{N_1, N_2} $
 generated over $\Z $
by operators $T_\l $
for all primes $\l\nmid N =N_1N_2 $
and $U_\l $
for all $\l | N, $
acting on the modular forms of weight two and level $\Gamma_0 (N)$
which are new  at all factors $\l | N_2. $
If $I $
is the kernel of the projection $\T_{N_1N_2, 1}\to\T $, then we set \begin{equation}
    J_{\min}^{N_1,N_2}\coloneqq J_0 (N)/IJ_0 (N), 
\end{equation}
an abelian variety with a (faithful) action of $\T $. If $N_1, N_2 $
are clear from context, we will omit the superscript.

For any abelian variety $A $
with an action of $\T $, and any maximal ideal $\m\subset\T $, the $\m $-adic Tate module is defined  to be the localization \begin{equation}
    T_\m A\coloneqq T_pA\otimes_\T\T_\m, 
\end{equation}
where $p $
is the residue characteristic of $\m $. (Note that this is dual to the notation of \cite{helm2007maps}.) For any $\m $
which is non-Eisenstein with odd residue characteristic $p\nmid N $, it follows from \cite{tilouine1997hecke} that $T_\m J_{\min}$
is free of rank two over $\T_\m$; by \cite[Corollary 4.7]{helm2007maps}, the natural map then induces an isomorphism \begin{equation}\label{all endomorphisms}
    \T_\m\xrightarrow{\sim}\End_{\T}(J_{\min})_\m.
\end{equation}
\subsubsection{}
Now suppose that $A $
is an abelian variety with a Hecke-equivariant isogeny to $J_{\min} $. For any $\l | | N_2, $
let $\mathcal A_{/\Z_\l} $
be the N\'{e}ron model of $A $. The connected component $\mathcal A ^ 0_{\F_\l} $
of the special fiber of $\mathcal A $
is a torus, and we write $\xchi_\l (A) =\Hom (\mathcal A ^ 0_{\F_\l},\mathbb G_m) $
for its character group. The association $A\mapsto\xchi_\l (A) $
is contravariantly functorial.
\begin{prop}[Helm]\label{Tate module}
Let $\m\subset \T $
be non-Eisenstein of residue characteristic $p\nmid 2N $. Then the natural maps induce $\T_\m $-module isomorphisms:
$$T_\m J_{\min}\otimes\Hom(J_{\min}, A)_\m\xrightarrow\sim T_\m A, $$
$$\xchi_\l (J_{\min} ^\vee)\otimes\Hom (J_{\min}, A)_\m\xrightarrow\sim\xchi_\l (A^\vee),$$
 $$\Hom (A, J_{\min})_\m\xrightarrow\sim\Hom_{\T_\m}\left (\Hom (J_{\min},
 A)_\m,\End(J_{\min})_\m\right).$$
 Here, all Hom-sets are understood to be $\T$-equivariant morphisms, and tensor products are taken modulo $\Z $-torsion. 
\end{prop}
\begin{proof}
This follows by duality from \cite[Corollary 4.1, Theorem 4.11, Proposition 4.14]{helm2007maps}. 
\end{proof}
We record the following elementary lemma for later use.

\begin{lemma}\label {Xlbound2} 
Let $\xchi =\xchi_\l (J_{\min} ^\vee)_\m $
for some $\l| N_2 $
and $\m\subset\T $, where  $\m $
is non-Eisenstein of odd residue characteristic $p $. If the associated residual representation $\overline\rho_\m $
is ramified at $\l $, then $\xchi $
is free of rank one over $\T_\m $. In general, there exist $\T_\m $-module maps $$\xphi_i:\xchi\to\T_\m,\;\;\psi_i:\T_\m\to\xchi,\;\; i = 1, 2 $$
such that $$\xphi_i\circ\psi_i =\psi_i\circ\xphi_i=t_i\in\T_\m\subset\End (\xchi) $$
and $$t_1+ t_2 = \l -1\in\T_\m. $$
\end{lemma}
\begin{proof}
If $\l -1
$
is a $p $-adic unit, or if $\overline\rho_\m $
is ramified, then this  follows from \cite[Lemma 6.5]{helm2007maps}. In general, we have\begin{equation}\xchi = \Hom((\mathcal{J}^\vee_{\min})_{\F_\l}^0[\m^\infty], \mu_{p^\infty})\end{equation} so that $\xchi$ may be identified with a $\T_\m [G_{\Q_\l}] $-module quotient  \[\pi: T_\m J_{\min} \to \xchi;\]
the Galois action on $\xchi $
is unramified and Frobenius acts as $U_\l $, which is a constant $\pm$1
because the residue characteristic of $\m $
is $p >2. $

Because $T_\m J_{\min} $
is free of rank two over $\T_\m $, it may be equipped with a basis $\set {e_1, e_2} $, and moreover an alternating $\T_\m $-module pairing
\begin{equation}
    \langle\cdot,\cdot\rangle: T_\m J_{\min}\times T_\m J_{\min}\to\T_\m
\end{equation}
such that\begin{equation}\label{alternating identity}
   y = \langle e_1, y\rangle e_2 -\langle e_2, y\rangle e_1
\end{equation}
for all $y\in T_\m J_{\min}. $ Define maps\begin{align*}
    \xphi_i: T_\m J_{\min}&\to\T_\m,\;\;i = 1, 2\\
    \xphi_1: y &\mapsto\langle y, (F - U_\l) e_2\rangle\\
    \xphi_2:y &\mapsto\langle y, (F - U_\l) e_1\rangle,
\end{align*}
where $F\in G_{\Q_\l} $
is any lift of Frobenius. We first claim that the maps $\xphi_i $
factor through $\pi $. Since $\T_\m $
is $p $-torsion-free, it suffices to check this after inverting $p $. On $T_\m J_{\min}\otimes\Q_p $, $F $
acts with distinct eigenvalues $U_\l $
and $\l U_\l $, and $\pi\otimes\Q_p: T_\m J_{\min}\otimes\Q_p\to\xchi\otimes\Q_p $
coincides with the projection onto the $U_\l $-eigenspace.
Since $\langle \cdot,\cdot\rangle $
is alternating and $\T_\m $-linear,
it follows that each $\xphi_i$
does indeed descend to a $\T_\m $-module map $\xchi\to\T_\m $. Now define maps 
\begin{align*}
    \psi_i:\T_\m&\to\xchi,\;\;i = 1, 2\\
    \psi_1: 1 &\mapsto U_\l \pi (e_1)\\
    \psi_2:1&\mapsto - U_\l \pi (e_2).
\end{align*}
We claim that $\psi_i $
and $\xphi_i $
satisfy the conclusion of the lemma. One readily calculates:
\begin {align*}
\xphi_1\circ\psi_1 (1) & =U_\l\langle e_1, (F - U_\l) e_2\rangle\\
\psi_1\circ\xphi_1 (e_1) & = U_\l\langle e_1, (F - U_\l) e_2\rangle\pi (e_1)\\
\psi_1\circ\xphi_1 (e_2) &= U_\l\langle e_2, (F - U_\l) e_2\rangle\pi (e_1)\\
& = U_\l\langle e_1, (F - U_\l) e_2\rangle\pi (e_2) - U_\l (F - U_\l)\pi (e_2) \\
& =U_\l\langle e_1, (F - U_\l) e_2\rangle\pi (e_2),
\end {align*}
where in the last two steps we have used (\ref{alternating identity})
and the fact that $F= U_\l $
on $\xchi $.
Similarly, $$\xphi_2\circ\psi_2 =\psi_2\circ\xphi_2 =- U_\l\langle  e_2, (F - U_\l) e_1\rangle, $$
and $$U_\l\langle e_1, (F - U_\l) e_2\rangle- U_\l\langle  e_2, (F - U_\l) e_1\rangle=\tr_{T_\m J_{\min}} U_\l (F - U_\l) =\l -1. $$
\end{proof}
\subsection{Shimura curves}
\subsubsection{}\label{Shimura curves}
If $\nu (N_2) $
is even, then there exists a Shimura curve $X_{N_1, N_2},$
with $\Gamma_0 (N_1) $
level structure, associated to the indefinite quaternion algebra over $\Q $
of discriminant $N_2. $
Let $$J^ {N_1, N_2}\coloneqq J (X_{N_1, N_2}),$$
an abelian variety with a natural action of $\T $
by correspondences (induced by Picard functoriality). When $N_1 $ and $N_2 $ are understood, we abbreviate $J = J^{N_1, N_2} $. There is a noncanonical Hecke-equivariant isogeny $J\to J_{\min} $. Consider the following technical hypothesis on the residual representation $\overline\rho_\m: G_\Q\to GL_2 (\T/\m) $ associated to $\m $:
\begin{equation}
  \tag {$\ast $} 
\begin{split}
  &\text{If } p= 3 \text{ and } \overline\rho_\m \text{ is induced from a character of } G_{\Q (\sqrt{-3})},  \;\exists\; \l|| N_2 \\&\text{such that either } \l\equiv -1\pmod 3 \text{ or } \overline\rho_\m \text{ is ramified at } \l.
\end{split}
\end{equation}
\begin{thm}[Helm]\label{Helm theorem}
Let $\m \subset\T$
be a non-Eisenstein maximal ideal of residue characteristic $p\nmid 2N $
satisfying ($\ast$). Then there is an isomorphism of $\T_\m $-modules: $$\Hom (J_{\min}, J)\simeq\otimes_{\l | N_2}\xchi_\l (J_{\min} ^\vee)_\m, $$
modulo $\Z $-torsion on the right-hand side.
\end{thm}
\begin{proof}
This is essentially \cite[Theorem 8.7]{helm2007maps}; to complete the case $p = 3, $
by \cite[Remark 8.12]{helm2007maps} one only needs a level-raising input that is provided by \cite{diamond1994lifting}.
\end{proof}
\subsection{Shimura sets}\label{Shimura set}
Now suppose that $\nu (N_2) $
is odd, and consider the finite double coset space (often called a \textbf{Shimura set}):
\begin{equation}
    X_{N_1, N_2}\coloneqq R (\mathbb A_\Q) ^\times\backslash B (\mathbb A_\Q) ^\times/B (\Q) ^\times,
\end{equation}
where $B $
is a definite quaternion algebra over $\Q $
ramified at $N_2 $
and $\infty $, and $R $
is an Eichler order in $B $
of level $\Gamma_0 (N_1). $
When $N_1 $
and $N_2 $
are clear from context, the subscripts may be omitted.
\subsubsection{}
The $\Z $-module $\Z [X] ^ 0 $
of formal degree-zero divisors in $X $
has two natural actions of $\T =\T_{N_1, N_2} $
by correspondences: an ``Albanese'' action induced by viewing an element of $\Z [X]^0 $
as a formal sum of points in a double coset  space, and a ``Picard'' action induced by identifying $\Z [X] =\Hom_{\text {Set}} (X,\Z) $. We will consider $\Z [X]^0 $
as a $\T $-module through the latter action.
The analogue of Theorem \ref{Helm theorem} is:
\begin{thm}\label{Helm theorem odd}
Let $\m \subset\T$
be a non-Eisenstein maximal ideal of residue characteristic $p\nmid 2N $
satisfying ($\ast$). Then there is an isomorphism of $\T_\m $-modules: $$\Z[X]^0\simeq\otimes_{\l | N_2}\xchi_\l (J_{\min} ^\vee)_\m, $$
modulo $\Z $-torsion on the right-hand side.
\end{thm}
\begin{proof}
Choose any prime $q | N_2, $
 so that $\nu (N_2/q) $
 is even. Let $\T' =\T_{N_1q, N_2/q},$
 and write $\m $ as well
 for the maximal ideal of $\T' $
 induced by the map $\T'\to\T $.

 Applying Theorem \ref{Helm theorem}
 to the pair $N_1q,N_2/q $, we obtain an isomorphism of $\T'_{\m} $-modules (modulo $\Z $-torsion)
\begin{equation}
    \Hom (J_{\min}^{N_1q,N_2/q}, J^{N_1q, N_2/q})_\m\simeq\otimes_{\l | N_2/q}\xchi_\l (J_{\min} ^ {N_1q,N_2/q,\vee}).
\end{equation}
 By \cite[Corollary 5.3, Lemma 8.2]{helm2007maps}, this implies an isomorphism of $\T_\m $-modules \begin{equation}
     \Hom (J_{\min}^ {N_1, N_2}, J^{N_1, N_2/q}_{ q\new})_\m\simeq\otimes_{\l | N_2/q}\xchi_\l (J_{\min} ^{N_1,N_2,\vee}),
 \end{equation} 
  where $J^{N_1q, N_2/q}_{q\new} $
  is the $q $-new quotient of $J^{N_1q, N_2/q} $.  Then, by Proposition \ref{Tate module}, we have \begin{equation}
      \xchi_q (J^{qN_1, N_2/q,\vee}_{ q\new} )_\m\simeq\xchi_q (J_{\min}^{N_1,N_2,\vee})_\m\otimes_{\l | N_2/q}\xchi_\l (J_{\min} ^{N_1,N_2,\vee})_\m.
  \end{equation} 
  By \cite[Proposition 5.3]{bertolini2005iwasawa}, $\xchi_q (J^{qN_1, N_2/q,\vee}) $
  is canonically identified with $\Z [X_{N_1, N_2}] ^ 0.$
  It remains to show that the inclusion $J^{N_1q, N_2/q,\vee}_{q\new} \hookrightarrow J^{N_1q, N_2/q,\vee} $
  induces an isomorphism on character groups at $q $. Indeed, since $J^{N_1q, N_2/q, \vee}_{q\new}  $
  has purely toric reduction at $q $, there is a surjection of character groups $\xchi_q (J^{N_1q, N_2/q,\vee})\twoheadrightarrow\xchi_q(J^{N_1q, N_2/q, \vee}_{q\new}) $,  which is clearly an isomorphism after tensoring both sides with $\Q $, hence also before.

 \end{proof}
 \subsection{CM points}
 \subsubsection{}\label{embeddings}
 Let us now fix an imaginary quadratic field $K/\Q $, and, once and for all, an embedding $K\hookrightarrow GL_2(\Q)$ such that $K\intersection M_2(\Z) = \O_K$.
 
 For any integer $N = N ^ + N ^ - $
 such that every prime factor of $N ^ + $
 is (unramified and) split in $K $, and $N ^ - $
 is a squarefree product of primes (unramified and) inert in $K $,  let $B =B_{N^-}$
 be the quaternion algebra over $\Q $
 ramified exactly at the factors of $N ^ - $ (and possibly $\infty $), and $R\subset B $
 an Eichler order of level $N ^ + $. For each such $B$, we fix an optimal embedding $K\hookrightarrow B $
 such that $K\intersection R =\O_K $. Then we may define the space of $K$-CM points: $$\mathscr C_{N^ +, N ^ -} = K ^\times\backslash B (\mathbb A_f)/\widehat R^\times.$$
 The Galois group $G_K ^ {\ab}$ acts on $\mathscr C_{N ^ +, N ^ -} $ through the reciprocity map  \begin{equation}\label{reciprocity map}\rec : G_K ^ {\ab}\xrightarrow {\sim} K ^\times\backslash\widehat K ^\times, \end{equation}
 i.e. $\sigma[b] = [\rec(\sigma)b].$

\subsubsection{} If $S\subset M_\Q $
 is a finite  set of primes, consider the set of $S$-CM points:
 $$\mathscr C ^ S= K ^\times\backslash \widehat K ^\times GL_2 (\mathbb A_S)/\widehat O^\times_K GL_2 (\widehat \Z_S). $$
 We again have an action of $G_K^{\ab}$, and by \cite[Proposition 2.5]{nekovar2007euler} the field of definition $K[y]$ of any $y\in\mathscr C ^ S $
 is contained in the compositum $K[S]$ of all ring class fields unramified outside $S$. 
 If $S$ is  
 disjoint from the factors of $N$, 
 then there is a unique map $GL_2(\mathbb A_S)\hookrightarrow B (\mathbb A_S)$ so that the composite $$K \otimes_\Q \mathbb A_S\hookrightarrow GL_2(\mathbb A_S)\hookrightarrow B(\mathbb A_S)$$ agrees with the embedding deduced from $K\hookrightarrow B$. This embedding identifies $\mathscr C^S$ with a Galois-stable subset of $\mathscr C_{N^+,N^-}$. 
 
 \subsubsection{}\label{cm on curve}
 If $\nu (N ^ -) $ is \textbf{even},
 then the Shimura curve $X=X_{N ^ +, N ^ -}
 $
 of (\ref{Shimura curves}) admits a complex uniformization:
 \begin{equation}\label{ complex uniformization}
     X(\C) = B ^\times\backslash\mathcal H ^ {\pm}\times B (\mathbb A_f) ^\times/\widehat R ^\times,
 \end{equation}
 where $\mathcal H ^\pm=\C\setminus \mathbb R $. 
We have an injective map:
 \begin{equation}
 \begin{split}
    \CM_{N ^ +, N ^ -} :  \mathscr C_{N ^ +, N ^ -}&\to X (\C)\\
    y&\mapsto [(h_0, y)],
 \end{split}
 \end{equation}
 where $h_0 $
 is the unique fixed point on $\mathcal H ^ + $
 of the action of $K $
 through
the embedding $K\hookrightarrow B $.
 By Shimura's reciprocity law, the image of $\CM_{N^+,N^-} $
is contained in $X(K^{\ab}), $
and $\CM_{N^+,N^-}$ is Galois-equivariant for the action on $\mathscr C_{N^+,N^-}$ defined above. 
\subsubsection{}
If instead $\nu (N ^ -) $
is \textbf {odd}, then there is a natural projection $$\CM_{N ^ +, N ^ -}:\mathscr C_{N ^ +, N ^ -}\to X_{N ^ +, N ^ -}, $$
where $X_{N ^ +, N ^ -} $ is the Shimura set of (\ref{Shimura set}).
\subsubsection{}
With notation as above in (\ref{embeddings}), let $q \nmid N$
be a prime inert in $K $
and not in $S $; note that the unique prime of $K $ above $q $
splits completely in $K  [S]/K $. We fix a prime $\mathfrak q$ of $K[S]$ lying above $q$ according to a fixed embedding $G_{\Q_q}\hookrightarrow G_\Q$. If $\nu(N^-)$ is \textbf{even}, then all points of $\CM_{N^+,N^-}(\mathscr  C ^S) $
have supersingular reduction modulo $\mathfrak q $, and the set of supersingular points $X(\mathbb F_{q ^ 2}) ^ {ss} $
may be identified with the Shimura set $X_q =X_{N ^ +, N ^ - q} $
associated to the definite quaternion algebra $B_q $
ramified at $N ^ - q \infty $.  Consider the following hypothesis on the residual representation $\overline\rho_\m$ associated to a non-Eisenstein maximal ideal $\m\subset \T_{N_1,N_2}$ of residue characteristic $p$ (for any Hecke algebra $\T_{N_1,N_2}$):
 \begin{equation}\label{Taylor Wiles}
     \tag{TW} \text{if }p=3,\text{ then }\overline \rho_\m\text{ is absolutely irreducible over } \Q{\sqrt{-3}}.
 \end{equation}
 Note that this is strictly stronger than condition ($\ast$) above.
\begin{prop}\label{reduction cd}
With identifications chosen compatibly, there is a commutative diagram:
\begin {center}
\begin {tikzcd}
\Z\left[ \mathscr C ^ S \right] ^ 0\arrow [r]\arrow [d,equal] &\Z\left[\CM_{N^+,N^-}(\mathscr C ^ S)\right]^0 \arrow [r]\arrow [d,"\Red_{\mathfrak q}"] & J (K[S])\arrow [d, "\Red_{\mathfrak q}"]\\
\Z \left[\mathscr C ^ S\right] ^ 0\arrow [r]&\Z [X_q] ^ 0\arrow [r] &J (\mathbb F_{q ^ 2})
\end {tikzcd}
\end {center}
 Moreover, the map $\Z [X_q] ^ 0\to J (\mathbb F_{q ^ 2}) $
is compatible with the action of $\T_{N ^ + q, N ^ - Q} $
where $U_q $
acts on $J (\mathbb F_{q ^ 2}) $
through $\Frob_q $, and is surjective after localizing at a non-Eisenstein maximal ideal $\m\subset\T_{N ^ + q, N ^ - Q} $
satisfying condition (\ref{Taylor Wiles}).
\end{prop}
\begin{proof}
 See \cite[Lemma 5.4.3]{zhang2001gross} for the commutativity of the diagram and
 \cite{ribet1990modular}  for the $U_q $ action. The surjectivity is an application of Ihara's Lemma which can be deduced from the argument in \cite[Proposition 9.2]{bertolini2005iwasawa}: we add auxiliary level of the form $\Gamma_1(\l)$, where $\l\nmid N$
 is a prime such that $\l -1, T_\l-\l-1\not\in\m$. That such a prime exists follows from condition (\ref{Taylor Wiles}) by \cite [Lemma 3] {diamond1994lifting}.
\end{proof}
\subsubsection{}
Now suppose instead that $\nu (N ^ -) $ is \textbf{odd}, and let $q\nmid N ^ - $
be a prime inert in $K $ and not lying in $S $. The Shimura curve $X_q = X_{N ^ +, N ^ - q} $
has a canonical, semistable integral model over $\Z_q $, whose irreducible components are  identified with two copies of $X = X_{N ^ +, N ^ -} $. We denote this set by $X ^\pm $, where the positive
copy is the one containing the reduction of the point $[(h_0, 1)] $
in the uniformization of (\ref{ complex uniformization}). We define a map $\mathscr C^S\to X^\pm$ by the composition $\mathscr C^S\xrightarrow{\CM_{N^+,N^-}} X\simeq X^+\subset X^\pm$. 

\subsubsection{}The N\'{e}ron model $\mathcal J_q $
of the Jacobian $J_q = J ^ {N ^ +, N ^ - q} $
 has purely toric reduction, and we write
$\xchi $ and $\Phi$
 for the character group and the group of connected components, respectively, of its special fiber.
 Recall the rigid-analytic uniformization of $J_q $, which gives rise to an exact sequence:
 \begin{equation}\label{rigid analytic}
     0\to\xchi\to\xchi ^\dagger\otimes\overline\Q_q\to J_q (\overline\Q_q)\to 0.
 \end{equation}
 Here, $\xchi ^\dagger =\Hom (\xchi,\Z) $, and the maps are Hecke-equivariant if $\xchi $
 is given Hecke action through Albanese functoriality, and the actions on $\xchi^\dagger$
 and $J_q (\overline\Q_q) $
 are induced by Picard functoriality.
 Importantly for our later applications, (\ref{rigid analytic}) is compatible  with the
 Galois action of $G_{\Q_q} $ \cite{bosch1996component}, where
 the action on $\xchi $
 is unramified and Frobenius acts through $U_q$ \cite[Proposition 3.8]{ribet1990modular}. The rigid analytic uniformization is related to the monodromy pairing $j:\xchi\to\xchi^\dagger$ of Grothendieck \cite{grothendieck1972modeles} by the commutative diagram with exact rows:
 \begin {center}
\begin {tikzcd}
0\arrow [r] &\xchi\arrow [d, equal]\arrow [r]&\xchi ^\dagger\otimes\Q_{q ^ 2}\arrow [r] \arrow [d, "\ord"]& J_q (\Q_{q ^ 2})\arrow [r]\arrow [d, "\Sp_q"] & 0\\0\arrow [r] &\xchi\arrow [r, "j"] &\xchi ^\dagger\arrow [r] &\Phi\arrow [r] & 0
\end {tikzcd}
\end {center}
In particular, the specialization map is well-defined on $J_q (K [S]) $.
 \begin{prop}\label{specialization cd}
 With identifications chosen compatibly,  there is a commutative diagram:
 \begin {center}
\begin {tikzcd}
\Z\left[ \mathscr C ^ S \right] ^ 0\arrow [r]\arrow [d,equal] &\Z\left[\CM_{N ^ +, N ^ - q} (\mathscr C ^ S)\right]^0 \arrow [r]\arrow [d,"\Sp_{\mathfrak q}"] & J_q(K[S])\arrow [d,"\Sp_{\mathfrak q}"]\\
\Z \left[\mathscr C ^ S \right] ^ 0\arrow [r]&\Z [X^\pm] ^ 0\arrow [r] &\Phi
\end {tikzcd}
\end {center}
 Moreover, 
the map $\Z[X^\pm]^0\to \Phi$
is compatible with the action of $\T_q =\T_{N ^ + q, N ^ -} $, where $U_q$ acts on $\Z [X ^\pm] ^ 0 $ by the matrix $$\begin{pmatrix} T_q & q\\ -1 & 0 \end{pmatrix}.$$
After localizing at a non-Eisenstein  maximal ideal $\m\subset \T_{q}$, the map $\Z[X^\pm]^0\to\Phi$ induces an isomorphism $$\Z [X ^\pm] ^ 0_\m\otimes_{\T_q}\T_q/(U_q ^ 2-1)\xrightarrow\sim \Phi_\m. $$
 \end{prop}
 \begin{proof}
  For the commutativity of the diagram, see \cite[Lemma 5.4.6]{zhang2001gross}; for the rest, see  
  \cite[Proposition 5.5]{bertolini2005iwasawa}.
 \end{proof}
\section {CM classes in cohomology}\label{ construction section}
\subsection{Level raising}\label{ level raising section}
\subsubsection{}\label{construction notation}
 Fix notation as in (\ref{modular forms notation}), and suppose that $N = N ^ + N ^ - $ is coprime to $pD_K$,
 where all factors of $N ^ + $
 are split in $K $
 and $N ^ - $
 is a squarefree product of primes inert in $K $. (In particular, $H^0(G_K, \overline T_f) = 0.$) We shall denote by $\pi $ a uniformizer of $\O $. From now on, we additionally assume that the maximal ideal $\m $ associated to $T_f $
 satisfies (\ref{Taylor Wiles}) above, i.e.:
 \begin{equation}
     \tag{TW} \text{if }p=3,\text{ then }\overline T_f\text{ is absolutely irreducible over } \Q{\sqrt{-3}}.
 \end{equation}
  \subsubsection{}
  We say a prime $q \nmid N$ is weakly admissible with sign $\epsilon_q=\pm 1$ if $q$ is inert in $K $, $a_q\equiv\epsilon_q(q+1)\pmod \p,$
 and   $q\not\equiv 1\pmod p. $
A weakly admissible pair $\set{Q,\epsilon_Q}$ is an ordered pair of a squarefree number $Q $ and a function  $\epsilon_Q: \set{q | Q}\to\set{\pm 1}$ such that $q$ is weakly admissible with sign $\epsilon_Q (q) $
 for all $q | Q $.
If $\set {Q,\epsilon_Q} $
is a weakly admissible pair, then for all $q | Q $, there is a unique root $u_q\in \O$
of the polynomial $y ^ 2 - ya_q+q$ such that $u_q\equiv\epsilon_Q (q)\pmod\p $. We may view $\O$ as a $\T_{N ^ + Q, N ^ -} $-algebra by letting $U_q$ act through $u_q$; let $\m_Q^{\epsilon_Q}$ be the associated maximal ideal (we will usually drop the superscript). 

If $Q = Q' Q'' $, then we abbreviate $\T^{Q'}_{Q''} = \T_{N ^ + Q', N^-Q''} $, omitting any superscript or subscript which is equal to 1.
\subsubsection{}

 In light of the structural similarity of Theorems \ref{Helm theorem}  and \ref{Helm theorem odd}, let \begin{equation}
     M_{Q} =\begin{cases}\Hom (J_{\min} ^ {N^+, N^-Q}, J ^ {N^+,N^-Q}), &\nu (N^-Q)\text { even},\\
     \Z [X_{N ^ +, N ^ - Q}] ^ 0, &\nu (N ^ -Q)\text { odd}.\end {cases}
 \end{equation}
 It is well-known that $M_Q$ is a faithful $\T_Q$-module, and indeed $M_Q\otimes\Q $
 is free of rank one over $\T_Q\otimes\Q $.
 \begin{lemma}\label{approximation Lemma}
 Suppose $\set {Q,\epsilon_Q} $
 is a weakly admissible pair, and
 let $$C= \sum_{\substack{\l | N ^ - \\ \overline T_f\text{ unram at }\l}}\ord_\pi(\l-1).$$
 Then there exists an $\O$-module map $$M_Q\otimes_{\T ^ Q}\O\to\T_Q\otimes_{\T^Q}\O$$
 with kernel and cokernel annihilated by $\pi ^ C $; in particular, 
 $\pi ^ {C} (M_Q\otimes\O) $
 is principal of length at least $\lg(\T_Q\otimes\O) - 2C. $
 \end{lemma}
 \begin{proof}
 We may assume that $\m_Q\subset \T ^ Q $
 descends to  $\T_Q $. Now, by Theorems  \ref{Helm theorem}  and \ref{Helm theorem odd},
 we have $$ M_{Q,\m_{Q}}\simeq\otimes_{\l |N ^ - Q}\xchi_\l \left(J_{\min} ^ {N ^ +, N ^ - Q,\vee}\right)_{\m_{Q}}, $$
  modulo $\Z $-torsion on the right.
  Lemma \ref{Xlbound2} implies that there exist a collection of $\T_Q$-module maps $$\xphi_i: M_{Q,\m_{Q}} \to\T_{Q,\m_Q}, \;\;\psi_i:\T_{Q,\m_Q}\to M_{Q_n,\m_{Q_n}},\;\; i=1,\ldots, r$$
  such that $$\xphi_i\circ\psi_i =\psi_i\circ\xphi_i = t_i\in\T_{Q,\m_Q}\subset\End (M_{Q_n,\m_{Q_n}}) $$
  and $$t_1+\ldots + t_r = \prod_{\substack{\l | N ^ -\\ \overline T_f \text{ unram at }\l}}(\l-1)\in\T_{Q,\m_Q}.$$
  Since $\O $ is principal, we may choose some $i $
  such that the image of $t_i $
  in $\T_Q\otimes\O $
  divides $\pi ^ C $. Then $\xphi_i $
  and $\psi_i $
  induce $\O $-module maps
  $$M_Q\otimes\O\to\T_Q\otimes\O\to M_Q \otimes\O$$
  whose composition in either direction is multiplication by a divisor of $\pi ^ C $, which implies the result.
\end{proof}
\begin{thm}\label{level raising}
If $\set{Qq,\epsilon_{Qq}}$ is a weakly admissible pair, then $$\lg (\T_{Qq}\otimes_{\T^{Qq}}\O)\geq \lg\left(\frac{\T_Q\otimes_{\T^{Q}}\O}{a_q-\epsilon_q(q+1)}\right) - C,$$
where $C $ is the number of Lemma \ref{approximation Lemma}.
\end{thm}
\begin{proof}
The proof depends on the parity of $\nu (N ^ - Q) $.
 \begin{case}
 $\nu (N ^ -Q)  $
 is even.
 \end{case}
Let us abbreviate  $J^Q = J^{N^+,N^-Q}$ and $J ^ {Q}_{\min}= J ^ {N ^ +, N ^ - Q}_{\min} $.
Consider the composite $$M_{Qq}\to J^{Q} (\mathbb F_{q ^ 2})\to H ^ 1 \left(\F_{q}^2, T_{\m_{Q}}  J ^ {Q}\right)_{\m_{Qq}}\simeq M_{Q}\otimes \frac{T_{\m_{Q}} J^{Q}_{\min}}{(U_{q}-\epsilon_q)}$$
 induced from the diagram of Proposition \ref{reduction cd}, the Kummer map, and Proposition \ref{Tate module}. These are surjective maps of $\T_{Q}^q$-modules, where $U_{q}$
 acts on the three latter modules through $\Frob_{q}$. Since $T_{\m_{Q}}J_{\min}^{Q}$ is free of rank two  over $\T_{ Q,\m_{Q}}$ and $\Frob_{q} $
 acts with the characteristic polynomial $\Frob_{q} ^ 2- T_{q} \Frob_{q} + q $ (whose roots are distinct modulo $\m_Q$), we may fix an identification 
 $$\frac{T_{\m_{Q}} J^{Q}_{\min}}{U_{q}-\epsilon_q}\simeq \frac{\T_{Q,\m_{Q}}}{T_{q} - \epsilon_{q}(q+1)},$$
  considered as a $\T^q_{Q,\m_{Qq}}$-module  again through $U_{q} $
  acting by $\epsilon_q $.
  Tensoring with $\O $,
   we obtain a surjective map $$ M_{Qq}\to M_Q\otimes_{\T ^ {Qq}}\frac{\O}{a_q -\epsilon_q (q +1)},$$
    hence (by Lemma \ref{approximation Lemma}) a map of $\T^{Qq}$-modules $M_{Qq}\to \frac{\T_Q\otimes \O}{a_q -\epsilon_q (q +1)}$
   with cokernel annihilated by $\pi ^ C $.
   Since the action  of $\T ^ {Qq} $ on $M_{Qq} $
   factors through $\T_{Qq} $, we obtain, by taking eigenvalues, a surjection $\T_{Qq}\to\O/\pi ^ j $
   for some $$j\geq\lg\left(\frac{\T_Q\otimes_{\T^{Q}}\O}{a_q-\epsilon_q(q+1)}\right) - C.$$
 
\begin{case}
 $\nu (N ^ -Q)  $
 is odd. 
\end{case}
By Proposition \ref{specialization cd}, the action of $\T^q_{Q, \m_{Qq}}$
on $$M_{Q,\m_{Q}} \otimes_{\T_Q}\left(\T_Q^2/\im \begin{pmatrix}  T_{q}-\epsilon_q &q\\-1 & -\epsilon_q\end{pmatrix}\right) ,$$
 with $U_{q}$ acting by $\epsilon_q$,
factors through $\T_{Qq,\m_{Qq}} $. Hence the action of $\T^q_{Q}$
on $$A=M_Q\otimes_{\T ^ Q}\frac{\O}{a_q - \epsilon_q(q+1)}$$
likewise factors through $\T_{Qq} $  (again with
$U_{q}$ acting by $\epsilon_q$). The  conclusion of
Lemma \ref{approximation Lemma} 
implies that $A$ has a $\T^q_Q$-module map to $$\frac{\T_Q\otimes_{\T ^ Q}\O}{a_q-\epsilon_q (q)} $$
with cokernel annihilated by $\pi ^ C $, from which the result follows.
 \end{proof}
 \begin{rmk}\label{level raising remark}
 If $\set {\QT,\epsilon_\QT}\in \mathsf N $, then for $\mathfrak F $-many $n $
 there is a corresponding weakly admissible pair $\set {Q_n,\epsilon_{Q_n}} $, where $Q_n $
 is a sequence representing $\QT $. To be precise, if $\QT =\set {\q_1,\ldots,\q_r }$,
 we choose sequences $q_i^n$ representing each $q_i$; for $\mathfrak F$-many $n$, the product $Q_n =q_1 ^ n\cdots q_r ^ n $, equipped with sign function $\epsilon_{Q_n} (q_i ^ n) =\epsilon_\QT (\q_i) $, forms a weakly admissible pair $\set{Q_n,\epsilon_{Q_n}} $.
 It follows from the definition of $\mathsf N $
 and from the theorem that, for any $j\geq 0, $
 there exist $\mathfrak F $-many $n $
 such that $\T^{Q_n}\to \O\to \O/\pi ^ j $
 factors through $\T_{Q_n} $. We say that a sequence of weakly admissible pairs $\set {Q_n,\epsilon_{Q_n}} $ (defined for $\mathfrak F $-many $n $) represents the pair $\set {\QT,\epsilon_\QT} $ if it is obtained from this construction for some choice of representatives $q_i ^ n $.
 \end{rmk}
 \subsection{The CM class construction}
 \subsubsection{}
 Let $S\subset M_\Q $
 be a finite set of primes, and fix an element
 $$y\in\Z[\mathscr C^{S})]^0;$$
  let $K [y]\subset K [S] $
  be its field of definition and $G_y =\Gal (K [y]/K) $
  the corresponding Galois group.
 A weakly admissible prime $q $
 with sign $\epsilon_q $
 is called $j $-admissible if $a_q\equiv\epsilon_q (q +1)\pmod {\pi ^ j} $; in this case, $T_j\coloneqq T_f/\pi ^ j $
 has a unique subspace $\Fil_{q,\epsilon_q}^+T_j$, free of rank one over $\O_j\coloneqq\O/\pi ^ j $, on which $\Frob_q $
 acts as $q\epsilon_q $. We will omit the subscript $\epsilon_q$ when there is no risk of confusion.
 A weakly admissible set $\set {Q,\epsilon_Q} $
 is called $j $-admissible if $\lg (\T_Q\otimes_{\T^Q}\O)\geq j+2C$; note that each $q|Q$ is then necessarily $j $-admissible (since $U_q=\epsilon_q$ in $\T_Q\otimes \O$). Let $N_j $
 be the collection of $j $-admissible sets.
 
For any $j $-admissible prime with sign $\epsilon_q $, we define the ordinary subspace:
\begin {equation}
H ^ 1 _{\ord,\epsilon_q}(K_q, T_j)=\im\left (H^ 1 (K_q,  \Fil^+_{q,\epsilon_q}T_j)\to H ^ 1 (K_q, T_j)\right).
\end{equation}
Using  the map obtained from Schapiro's Lemma (e.g. \cite[\S3.1.2]{skinner2014iwasawa})
 \begin{equation}\label{shapiro}
     \Res_q:H^1(K[y], T_j) \to \Hom_{\text{Set}}(G_y,H^1(K_q, T_j)),
     \end{equation}
      we also have maps:\begin{align*}
    \partial_{q,\epsilon_q}: H^1(K [y], T_j) &\to
     \Hom_{\text{Set}}(G_y,H^1(I_q,\Fil_q^+T_j))\approx \O_j[G_y],\\
     \loc_{q,\epsilon_q}: H^1(K[y]^\Sigma/K[y], T_j) &\to\Hom_{\text{Set}}(G_y,T_j/\Fil_q^+T_j),\approx \O_j[G_y],\;\;q\not\in \Sigma,
 \end{align*}
defined as in (\ref{ localization definition},\ref{ residue definition}).
 \begin{construction}\label{construction and reciprocity}
If $\Sigma\subset M_\Q $
is the set of places dividing $Np\infty $, then for all  $\set {Q,\epsilon_Q}\in N_{j} $, there exist principal sub-$\O_j$-modules: \begin{align*}
    (\kappa_j (y,Q,\epsilon_Q))&\subset H ^ 1 (K [y] ^ {\Sigma\cup Q}/K [y], T_j), &\nu(N^-Q) \text{ even},\\
    (\lambda_j (y, Q,\epsilon_Q)) &\subset\O_j[G_y],&\nu(N^-Q)\text{ odd},
\end{align*}
compatible under the natural reduction maps for $j'\leq j$, and
satisfying the following properties.
\begin{enumerate}
    \item If $\set{Qq, \epsilon_{Q}}\in N_{j}$
    where  $\nu(N^-Q)$ is even, 
    then for all $q|Q$ and all $g\in G_y$, $$\Res_q(\kappa_j(y,Q,\epsilon_Q))(g)\subset H ^ 1_{\ord,\epsilon_Q (q)} (K_\q, T_j). $$
    \item If $\set{Qq,\epsilon_{Qq}},\set{Q,\epsilon_Q}\in  N_j$  where $\epsilon_Q =\epsilon_{Qq} |_Q $ and $\nu(N^-Qq)$ is even, then
  $$\partial_{q,\epsilon_{Qq} (q)}(\kappa_j(y,Qq,\epsilon_{Qq})) = (\lambda _j(y,Q,\epsilon_{Q}))\subset\O_j [G_y]. $$
  \item If $\set{Qq,\epsilon_{Qq}},\set{Q,\epsilon_Q}\in  N_j$  where $\epsilon_Q =\epsilon_{Qq} |_Q $ and $\nu(N^-Qq)$ is odd, then  $$\loc_{q,\epsilon_{Qq}(q)} (\kappa_j (y,Q,\epsilon_Q)) = (\lambda_j (y,Qq,\epsilon_{Qq}))\subset\O_j[G_y]. $$
\end{enumerate}
\end{construction}
\begin{proof}
The specifications $\epsilon_{Q}$ will be dropped to ease notation.
 Suppose first that $\nu(N^-Q)$ is odd. By Lemma \ref{approximation Lemma},  there is a unique  map (up to scalars) $M_Q\to\O_j $ 
 of $\T ^ Q $-modules that factors through multiplication by $\pi ^ C $ and is surjective after $\O $-linearization.
 We define $\lambda_j(y, Q)(g)$ to be the image of $gy\in \Z[\mathscr  C^S]^0$ by  the composite $\Z[\mathscr C^S]^0\to M_Q\to  \O_j$.
 
Now suppose that  $\nu(N^-Q)$ is even. Adopting the abbreviations of Case 1 of Theorem \ref{level raising},  $CM_{N ^ +, N ^ -Q}(y)$ is a formal divisor on $X_{Q} = X_{N ^ +, N ^ - Q}$ for each $n$, and its image in $J^{Q}$ is defined over $K[y]$. Let $$d(y, Q)\in H^1(K[y]^{\Sigma\cup Q}/K[y],T_{\m_{Q}}J^{Q})$$
be the Kummer image. By Lemma \ref{approximation Lemma} and Proposition \ref{Tate module}, there is a unique (up to scalars) map of $\T^Q[G_\Q]$-modules $T_{\m_Q}J^Q\to T_j$ that is surjective after $\O $-linearization. We define $\kappa_j(y, Q)$ to be the image of $d (y, Q) $
under the induced map $$H^1(K[y]^{\Sigma\cup Q}/K[y],T_{\m_{Q}}J^{Q})\to H^1(K[y]^{\Sigma\cup Q}/K[y],T_j). $$

Note that, for each $g\in G_y $,
$\Res_q\kappa_j (y,Q) (g)$ is the local Kummer image of $g\cdot \CM_{N ^ +, N ^ -} (y)=\CM_{N ^ +, N ^ -} (gy)  $ by Shimura's reciprocity law.
\begin{enumerate}
    \item This follows from the rigid analytic uniformization (\ref{rigid analytic}) by the argument of \cite[p. 15]{gross2012local}. Indeed, the argument there shows that the image of the Kummer map $$J^{Qq}(\Q_{q^2})\to H^1(\Q_{q^2}, T_{\m_{Q}} J^{Q})$$
    agrees with the image of the map $$H^1(\Q_{q^2}, \xchi(J^{Q}) (1))\to H^1(\Q_{q_n^2}, T_{\m_{Q}} J ^ {Q}),$$ where  $\xchi(J^{Q})(1) \to  (T_{\m_{Q}} J^{Q})[\Frob_{q} - U_{q} q]$ is a canonical isomorphism. Since $q\nmid p - 1,$ the surjection $$T_{\m_Q} J ^ Q\twoheadrightarrow T_j$$
    induces $$\xchi(J^{Q})(1) \twoheadrightarrow \Fil ^ +_{q} T_j, $$ and the claim follows.
    \item 
    We claim that  the 
    composite $$J^{Qq}(\Q_{q^ 2})\to H^1(\Q_{q^ 2}, T_{\m_{Q}}J^{Q})\to H^1(I_{q}, T_j)$$
    factors through $$\Sp_{q}: J^{Qq}(\Q_{q^ 2})\to \Phi_{Q,\m_{Q}}.$$
    Indeed, the target of the composite map has Frobenius eigenvalue $U_{q}$, and, because $p\nmid q - 1$, a diagram chase using (\ref{rigid analytic}) shows that the pro-$p$ part of the kernel of $\Sp_{q}$ has Frobenius eigenvalue $-U_{q}$ (if it is nontrivial at all). 
    
    Using the commutativity of the diagram in Proposition \ref{specialization cd} (and the fixed embedding $K[S]\hookrightarrow\Q_{q ^ 2} $),
    $\partial_q(\kappa_j(y,\QT))(g)$ is therefore the image of $gy$ under the composite of the canonical map $$\Z[\mathscr C ^ S]^0 \to M_{Qq}$$
     with \emph{some} surjective map of Hecke modules $M_{Qq} \to \O_j$, which factors through multiplication by $\pi^C$ by the choice of map $T_{\m_Q}J^Q\to T_j$. 
    We may conclude by  Lemma \ref{approximation Lemma}.
    \item The proof is very similar to (2), invoking instead Proposition \ref{reduction cd}.
\end{enumerate}
\end{proof}
\section{From CM classes to bipartite Euler systems}\label{construction per to}
\subsection{$p$-adic interpolation}\label{p adic interpolation}
\subsubsection{}
Suppose for this subsection that:
\begin{equation}
    \tag{spl} p \text { splits in } K
\end{equation}
and
\begin{equation}
    \tag{ord} a_p\not\in \p.
\end{equation}
\subsubsection{}Fix an auxiliary prime $\l_0\nmid Np$ such that $a_{\l_0} -\l_0-1\in \O ^\times $ (one exists by the irreducibility of $\overline  T_f $), and set $S=\set{\l_0,p}$.
For each $m\geq 0, $
consider the point $$y_{m,0} =\begin{pmatrix} p^m & 0\\0 & 1\end{pmatrix}\in GL_2 (\Q_p)\subset \mathscr C^S;$$
note that $y_{m,0}$ is defined over $K[p^m]$.
If $T_{\l_0} $
is the usual adelic Hecke operator, then set $y_m = (T_{\l_0} - \l_0 -1 )\Tr_{K[p^m]/K_m} y_{m,0} \in \Z[\mathscr C ^S]^0$,
where $K_m$ is the $m$th layer of the anticyclotomic $\Z_p$-extension. 

\subsubsection{}Now suppose given any $\set {\QT,\epsilon_\QT} $, and let $\set {Q_n,\epsilon_{Q_n}}$ be a representative sequence of weakly admissible pairs as in Remark \ref{level raising remark}.
Since $T_p\not\in\m$, Hensel's Lemma implies that  the Hecke algebras $\T_{Q_n,\m_{Q_n}}$  contain a (unique) element $u\not\in \m_{Q_n}$
such that $u ^ 2 - ua_p+p= 0. $  Suppose first that $|\mathsf Q|+\nu(N^-)$ is even.
By the  usual Heegner point  norm relations (cf. e.g.  \cite[Proposition 3.10]{darmon2004rational}), the classes $$d(y_m,Q_n)'\coloneqq u^{-m+1}d (y_m, Q_n) - u^{-m} \Res_{K_m/K_{m-1}} d_n(y_{m-1}, Q_n)$$
of Construction \ref{construction and reciprocity}
 are compatible under the corestriction maps $$H^1(K_m, T_{\m_{Q_n}}J^{Q_n})\to H^1(K_{m-1}, T_{\m_{Q_n}}J^{Q_n}).$$
If we replace $d (y_m, Q_n) $
by $d (y_m, Q_n)' $, then, for any $j $
and for $\mathfrak F $-many $n $, we obtain classes $\kappa_j(y_m,Q_n)'\in \H^1(K_m^{\Sigma \cup Q_n}/K_m, T_j)$
that are compatible under corestriction (as long the choices in the construction are made compatibly as  $m$
varies, which is clearly possible). 
We let $$\boldsymbol\kappa(\QT) \in \lim_{\substack{\longleftarrow\\{m,j}}}\H ^ 1(K_m ^ {\Sigma\cup\QT}/K_m, T_j)\simeq \H^1(K, T_f\otimes \Lambda(\Psi))$$
be the class represented by the family $\kappa_j(y_m, Q_n)'$.

\subsubsection{}
Similarly, if $|\QT | +\nu (N ^ -) $ is odd, the elements
$$\lambda_j(y_m,Q_n)' \coloneqq \alpha_p^{-m+1} \lambda (y_m,Q_n)- \alpha_p^{-m+1}\lambda(y_{m-1}, Q_n) \in \O_j[\Gal(K_m/K)]$$ are compatible under the natural projection maps $$\O_j[\Gal(K_m/K)]\to \O_j[\Gal(K_{m-1}/K)].$$ We then obtain an element $$\boldsymbol\lambda (\QT)\in\lim_{\substack{\leftarrow\\m,j}}\mathcal U\left(\set {\O_j[\Gal (K_m/K)]}_{n\in \N}\right)
\simeq
\O\llbracket\Gal (K_\infty/K)\rrbracket\simeq\Lambda. $$

Let $\ST\subset\mathsf M_K $
be the set of constant ultraprimes $\underline v $
such that $v | Np\infty $. We define a Selmer structure $(\mathcal F_\Lambda, \ST)$ for $\mathbf T \coloneqq T_f\otimes \Lambda$ in the usual way (see e.g. \cite{howard2004heegner,burungale2019proof}):
$$\H^1_{\mathcal F_\Lambda}(K_\v, \mathbf T)=\begin{cases}\im\left (H ^ 1 (K_v, \Fil_v ^ +\mathbf T)\to H ^ 1 (K_v,\mathbf T)\right), & \v=\underline v, v|p,\\
H^1(K_v, \mathbf T), &\v=\underline v,v\nmid p,\\
\H ^ 1_{\unr}(K_\v, \mathbf T), &\text{otherwise}.\end{cases}$$
Here $\Fil_v^+\mathbf T\subset\mathbf T$ is the unique free, rank-one direct summand on which $I_v$ acts by the cyclotomic character. By \cite[proposition 3.3.1]{howard2006bipartite}, $\mathcal{F}_\Lambda $ extends naturally to an interpolated self-dual Selmer structure $(\ST,\mathcal{ F}_\Lambda,\mathcal F_\P,\Sigma_\Lambda) $
for $\mathbf T $. 
\begin{prop}\label{instructed Euler system overline}
The pair $(\boldsymbol\kappa,\boldsymbol\lambda) $ is a nontrivial bipartite Euler system for the triple $(\mathbf T,\mathcal F_\Lambda,\ST) $.

\end{prop}
\begin{proof}
We first show that $\boldsymbol\kappa (\QT)$ lies in $\Sel_{\mathcal F_\Lambda(\QT)} (\mathbf T)$ for all
$\set{\QT,\epsilon_\QT}\in\mathsf N ^ {\nu (N ^ -)} $.
The only local conditions to verify are those at $v | p $; the local conditions for $\q\in\QT $ follow from Construction \ref{construction and reciprocity}(1), and the rest are trivial.
If $\QT $
is represented by the sequence $Q_n $, let $\Fil_v ^ + T_{\m_{Q_n}} J^{Q_n}$
 be the maximal $\T_{\m_{Q_n}}$ submodule on which $I_v $
 acts by the cyclotomic character (adopting the notation of  Construction \ref{construction and reciprocity} and if necessary restricting our attention to $\mathfrak F$-many $n$). As in \cite[Proposition 4.7]{chida2015anticyclotomic}, it suffices to show that, for all $m$ and $n$ and a fixed extension of $v$ to $K_\infty$, the image $d_{n,m}$ of the class $d(y_m, Q_n)'$  under the composite $$H^1(K_m, T_{\m_{Q_n}}J^{Q_n})\to H^1(K_{m,v}, T_{\m_{Q_n}}J^{Q_n}/\Fil^+T_{\m_{Q_n}}J^{Q_n})$$
 is trivial.
Since $d_n(y_m)'$
is a $\T_\m$-linear combination of Kummer images over $K_m $, by \cite[Example 3.11]{bloch2007functions} and \cite[Proposition 12.5.8]{nekovavr2006selmer} $d_{n, m} $
lies in the kernel of
$$H^1(K_{m,v}, T_{\m_{Q_n}}J^{Q_n}/\Fil^+T_{\m_{Q_n}}J^{Q_n})\to H ^ 1 (K_{m, v}, \Q_p\otimes T_{\m_{Q_n}}J^{Q_n}/\Fil^+T_{\m_{Q_n}}J^{Q_n}).$$
Since the classes $d_{n,m}$ are corestriction-compatible as $m$ varies, the argument of \cite[Proposition 2.4.5]{howard2007variation} shows that indeed $d_{n,m}= 0$ for all $n,m$.

The explicit reciprocity laws are a consequence of Construction \ref{construction and reciprocity}(2,3), and the nonvanishing of either $\boldsymbol\kappa (1) $ or $\boldsymbol\lambda (1)$ (according to the parity of $\nu (N ^ -) $) is due to the work of Cornut \cite{cornut2007nontriviality} and Vatsal \cite{vatsal2003special}.
\end{proof}
\subsection{Kolyvagin classes}

\subsubsection{}\label{gross subsection}
Before defining the Kolyvagin classes in patched cohomology, we begin by recalling a calculation explained in \cite{gross1991kolyvagin}.

Let $m $
be a squarefree product of primes $\l $
inert in $K $, and let $$y(m)_0\in\prod_{\l|m} GL_2 (\Q_{\l}) $$
be the element with $\l$th component $\begin{pmatrix}\l& 0\\ 0& 1\end{pmatrix} $. If $\l_0$ is the auxiliary prime of the previous subsection, and $S \set {\l_0}\cup\set{\l|m}$, then we define $$y(m) = (T_{\l_0} - (\l_0 + 1))y(m)_0\in \Z [\mathscr C ^ S] ^ 0. $$ 
\subsubsection{}
Note that $K [y (m)] = K [m]$ and that $\Gal(K[m]/K[1])\simeq \prod_{\l|m}\Gal (K [\l]/K [1]); $
each $\Gal (K [\l]/K [1]) $
is cyclic of order $\l +1.$
For any place $\lambda$ of $\overline\Q$ over $\l$, fixed for the time being, let $\Frob_\lambda\in G_\Q $
be a lift of absolute Frobenius, and  $\sigma_\lambda \in I_\lambda \subset G_K$ a generator of $\Gal (K [\l]/K [1]) $.
Recall the Kolyvagin derivative operators \cite{gross1991kolyvagin}:
$$D_\l = \sum_{i = 1}^\l i \sigma_\lambda^i \in \Z [\Gal (K [\l]/K [1])],\;\; D_m =\prod_{\l | m} D_\l.$$
Let $$P(m) =  D_m y(m)\in \Z [\mathscr C ^ S] ^ 0.$$
Finally, let $Q $
be a set of primes inert in $K $ that is disjoint from $S $, $p$, and the factors of $N$, and let $P (m, Q) =\CM_{N^+,N^-Q}(P (m)),$
$y(m, Q)  =\CM_{N^+,N^-Q}(y(m)).$
\begin{prop}[\cite{gross1991kolyvagin}, Proposition 3.7]\label{Gross facts}
For all $\l | m $, we have:
\begin{enumerate}
    \item $(\sigma_\l -1) P(m) = (\l +1)  D_{m/\l} y(m) - T_\l P (m/\l). $
    \item Suppose $\nu (N ^ - Q) $ is even. If  $\lambda$  lies over $\l$, then $$D_{m/\l}y (m, Q)\equiv\Frob_\lambda P({m/\l}, Q)\pmod {\lambda}.$$
\end{enumerate}
\end{prop}
\begin{prop}\label{ Kolyvaginswitcheroo}
Suppose $\m_Q\subset\T_Q =\T_{N ^ +, N ^ - Q} $
is a maximal ideal whose associated residual representation has no $G_{K [m]} $-fixed points, and let $I_m\subset\T_Q $ be the ideal generated by $\l +1 $ and $T_\l $ for all $\l | m $. Then if
$\nu (N ^ - Q) $ is even:
\begin{enumerate}
    \item Restriction induces an isomorphism $$\Res_m: H ^ 1 (K [1], T_{\m_Q} J ^ Q/I_m)\xrightarrow\sim H ^ 1 (K [m],T_{\m_Q} J ^ Q/I_m)^{\Gal (K [m]/K[1])}. $$
    \item The Kummer image $d'(m, Q)$ of $P(m, Q)$ in $H ^ 1 (K [m],T_{\m_Q} J ^ Q/I_m)$
    lies in the image of $\Res_m$.
    \item If $c(m, Q) = \cores_{K [1]/K}\Res_m^{-1}d'(m, Q)$, then for all $\l | m $ and any  choices of representatives, $$c (m,Q) (\sigma_\lambda) = \Frob_\lambda^{-1}d'(m/\l, Q) (\Frob_\lambda^2) \pmod I_{m}. $$
    \item The class $c (m, Q) $ is unramified at any place $v\nmid NpmQ\infty. $
\end{enumerate}
\end{prop}
\begin{proof}
(1) follows from the inflation restriction exact sequence as in \cite{gross1991kolyvagin}, and (2) is immediate from Proposition \ref{Gross facts}. Also (4) is clear from the construction. For (3), it suffices to check the corresponding statement for $c' (m, Q) =\Res_m ^ {-1} d' (m, Q) $. The proof is a modification of the argument in \cite{mccallum1991kolyvagin}. Fix division points $\frac {P (m, Q)} {\l+1}$ and $\frac{P(m/\l,Q)}{\l+1}$; one may verify that  $c'(m,Q) (\sigma_\l) $ is the unique element $A\in T_{\m_Q} J ^ Q/I_m$ such that, for all $g\in G_{K [m/\l]} $,
\begin{align*}
    (g -1) A & \equiv(g-1) (\widetilde\sigma_\l -1)\frac {P (m, Q)} {\l+1} \in T_{\m_Q} J ^ Q/I_m.
\end{align*}
By 
Proposition \ref{Gross facts}(1), $A $
is also the image of the (unique)
point $T\in J^Q[\l+1]$ such that $$T\equiv D_{m/\l} \CM_{N^+,N^-Q}(P(m, Q)) -  T_\l\frac {P (m/\l, Q)} {\l+1}  \pmod \lambda.$$
But by Proposition \ref{Gross facts}(2), this is equivalent to 
\begin{equation}
    \label{torsion cong}T\equiv \Frob_\lambda^ {-1} P (m/\l, Q) -  T_\l\frac {P (m/\l, Q)} {\l+1}  \pmod\lambda.
\end{equation}
By the Eichler-Shimura relation, and the fact that $\l$ splits completely in $K [m/\l] $, the image of $T$ in $ T_{\m_Q} J ^ Q/I_m $
is precisely $$\Frob_\lambda ^ {-1} d' (m/\l, Q) (\Frob_\lambda^2).$$ 
\end{proof}
\begin{definition}\label{Kolyvagin definition}
For a squarefree product $m $ of primes inert in $K $, let $I_m (f)\subset\O$
be the ideal generated by $a_\l (f) $ and $\l +1 $
for all $\l | m. $
Suppose $\set {Q,\epsilon_Q}\in\mathcal N_j $ is $j $-admissible and $j\geq v_\p (I_m (f)).$ 
If $\nu (N ^ - Q) $ is even,  then the Kolyvagin class \begin{equation}\overline c (m, Q)\in H ^ 1 (K ^ {\Sigma\cup Q\cup m}/K, T_f/I_m (f)) \end{equation}
is defined to be the image of $c (m, Q) $. If $\nu (N ^ - Q) $ is odd, then the reduction of $\lambda (P (m), Q) $ in $$(\O/I_{m}) [\Gal (K [m]/K)] $$ is constant on cosets of $\Gal (K [m]/K [1]) $ by Proposition \ref{Gross facts}(1) and therefore descends to \begin{equation}\lambda' (m, Q)\in (\O/I_m) [\Gal (K [1]/K)].\end {equation}
The Kolyvagin element is then defined as:
\begin{equation}
   \lambda (m, Q) =\tr_{K [1]/K}\lambda' (P (m), Q)\in\O/I_m.
\end{equation}
\end{definition}
\begin{rmk}
When $Q = 1 $ and $\nu (N ^ -) $ is even, this agrees with Kolyvagin's construction
 \cite{kolyvagin1991sha}.
\end{rmk}

For applications to the parity conjecture for $f $, we will require the following:
\begin{prop}\label{ parity proposition}
If $\nu (N ^ -) $ is even, then $\overline c(m, 1) $
lies in the ${\epsilon_f}\cdot (-1)^{\nu (m)+1}$-eigenspace for  the action of $\tau$. If $\nu (N^-)$ is odd and $\lambda(m, 1) \neq 0,$ then $\epsilon_f = (-1)^{\nu (m)}$.
\end{prop}
\begin{proof}
The maps $T_\m J^{N^+,N^-} \to T_f/\pi^j$ or $\Z[X_{N^+,N^-}]^0\to \O/\pi^j$ used in Construction \ref{construction and reciprocity} are equivariant for the action of the Atkin-Lehner involution because of the uniqueness property derived from Lemma \ref{approximation Lemma}. (Note this is not necessarily true of the corresponding maps at level $N^+N^-Q$, which are not necessarily reductions of genuine modular parameterizations.) Since the Atkin-Lehner eigenvalue of $f$ is $-\epsilon_f$, the proposition follows exactly as in
 \cite[Proposition 5.4]{gross1991kolyvagin}.
\end{proof}
\begin{definition}
An ultraprime $\mathsf l $
is called Kolyvagin-admissible if $$\Frob_{\mathsf l}\in\Gal (K (T_f)/\Q) $$
is a complex conjugation.
A Kolyvagin-admissible set is a finite set of Kolyvagin-admissible ultraprimes, and the collection of all Kolyvagin-admissible sets is denoted $\mathsf K $.
\end{definition}

\subsubsection{}
If $\mathsf l $
is Kolyvagin-admissible, then the local cohomology $$\H ^ 1 (K_{\mathsf l}, T_f) $$
is free of rank four over $\O $, and carries a natural action of the complex conjugation $\tau\in\Gal (K/\Q) $.
It has a canonical splitting of the finite-singular exact sequence: $$\H^1(K_{\mathsf l}, T_f)=\H^1_{\unr}(K_{\mathsf l}, T_f)\oplus \H^1_{\tr} (K_{\mathsf l}, T_f),$$
defined as follows. If the sequence $\l_n$ represents $\mathsf l$, then for any $j $ and for $\mathfrak F $-many $n $, $\Frob_{\l_n}$ acts as complex  conjugation on $T_f/\pi^j$, and $$H ^ 1_{\tr}(K_{\l_n}, T_f/\pi^j) = \ker \left( H ^ 1 (K_{\l_n}, T_f/\pi ^ j)\to H ^ 1 (K [\l_n]_{\lambda_n}, T_f/\pi ^ j)\right) $$
is isomorphic to $H ^ 1 (I_{\l_n},T_f/\pi^j)^{\Frob_\l ^ 2 = 1},$
where $\lambda_n $
is a prime of $K [\l_n] $ over $\l_n $.
Then $$H^1_{\tr} (K_{\mathsf l}, T_f) = \lim_{\leftarrow} \mathcal U \left (\set {H ^ 1_{\tr} (K_{\l_n}, T_f/\pi^j)}_{n\in \N}\right) \subset \H^1(K_{\mathsf l}, T_f)$$
is our transverse subspace.
We denote by $\loc_\mathsf l^\pm$ and $\partial_\mathsf l ^\pm $
the composites $\H^1(K, T_f) \to \H^1_{\unr} (K_\mathsf l, T_f)^\pm$ and $\H^1(K, T_f) \to \H^1_{\tr} (K_\mathsf l, T_f)^\pm$, respectively, where $\pm $
is the Frobenius eigenvalue. The codomain of each is free of rank one over $\O$.

Let $\ST \subset\mathsf M_K$  be the set of constant ultraprimes $\underline v $ such that $v\nmid Np\infty$. We will consider the Kolyvagin-transverse Selmer structure $(\mathcal{F}(\mathsf m), \ST \cup \mathsf m)$ on $T_f$, for any $\mathsf m\in \mathsf K $:
\begin{equation}
    \H ^ 1_{\mathcal F(\mathsf m)}(K_\v, T_f) = \begin{cases} \im\left(A(K_v) \otimes_{\O_f} \O \to H ^ 1 (K_v, T_f)\right), &\v=\underline v,\\
    \H ^ 1_{\tr} (K_\mathsf l, T_f), &\v=\mathsf l \in \mathsf m, \\
    \H^1_{\unr} (K_\v, T_f), & \text{otherwise}.
    \end{cases}
\end{equation}
Here $A$ is an (optimal) abelian variety with real multiplication by $\O_f$ associated to $f $.
If $\set {\QT,\epsilon_\QT}\in\mathsf N_{\mathsf m }$, then we denote by $(\mathcal F(\mathsf m, \mathsf Q), \ST\cup \mathsf m \cup \mathsf  Q) $
the modified Selmer structure of (\ref{where admissible set defined}).
\subsubsection{}\label{construction Kolyvagin}
Let  $\set{\mathsf Q,\epsilon_\QT}\in\mathsf N_{\mathsf m} ^ {\nu (N ^ -)}$, and fix representatives $Q_n$ and $m_n$, which we may assume to be disjoint. 
Our patched Kolyvagin class is  the  element
$$\kappa(\mathsf m, \QT)\in\H ^ 1 (K^ {\ST\cup \mathsf m\cup \QT}/K, T_f)$$
whose image in $T_j $
is represented by the sequence of images of the classes $\overline c(m_n, Q_n)$, well-defined for $\mathfrak F $-many $n $. 

If $\set {\mathsf Q,\epsilon_\QT}\in\mathsf N ^ {\nu (N ^ -) +1},$ 
then we similarly set $$\lambda (\mathsf m,\mathsf Q)\in\O\simeq\lim_{\longleftarrow}\mathcal U\left (\set {\O/\pi ^ j}\right) $$
to be the element whose image in $\O/\pi ^ j $
is represented by the sequence $\lambda (m_n, Q_n) $.
\begin{prop}\label{ Selmerstaff for Kolyvagin classes}
For any $\mathsf m\in\mathsf K $
and $\set {\QT,\epsilon_\QT}\in\mathsf N_{\mathsf m} ^ {\nu (N ^ -)}$,
$$(\kappa (\mathsf m,\QT))\subset\Sel_{\mathcal F (\mathsf m,\QT)} (T_f). $$
Moreover:
\begin{enumerate}
    \item For all $\mathsf l\in \mathsf m $, $$(\loc_{\mathsf l}^\pm(\kappa(\mathsf m/\mathsf l,\QT))) = (\partial_{\mathsf l}^{\mp}(\kappa(\mathsf m, \QT))$$
    as submodules of $\O .$
   \item For all $\q\in\QT $,
   $$(\partial_\q (\kappa (\mathsf m,\mathsf Q))) = (\lambda (\mathsf m,\mathsf Q/\q)) $$
   as submodules of $\O. $
   \item For all $\q\not\in \QT$, admissible with sign $\epsilon_\q $,
   $$(\loc_\q(\kappa (\mathsf m,\mathsf Q))) = (\lambda (\mathsf m,\QT\q)) $$
    as submodules of $\O $.
\end{enumerate}
In particular, for any fixed $\mathsf m$, $(\kappa (\mathsf m,\cdot),\lambda (\mathsf m,\cdot)) $ forms a bipartite Euler system with sign $\nu (N ^ -) $ for the triple $(T_f,\mathcal F (\mathsf m),\ST\union \mathsf m). $
\end{prop}
\begin{proof}
We verify the local conditions for each $\v\in\ST\union\mathsf m\union\QT $. If $\v =\underline v $
for a prime $v | N\infty $, then the local condition is  all of $H ^ 1 (K_v, T_f) $, so there is nothing to show. If $v | p $, then we show that, for all $j$, the image $c_j$ of $\Res_v \kappa (\mathsf m,\QT) $
in $H ^ 1 (K_v, T_f/\pi ^ j) $
is a Kummer image.
Recalling the notation used to construct $\kappa (\mathsf m,\QT) $,
the proof of \cite[Lemma 7]{gross2012local} implies that: \begin{align*}
    \delta_v(A(K_v))&=H^1_{\text{fl}}(K_v,T_f/\pi ^ j)\subset H^1(K_v,T_f/\pi ^ j)\\
    \kappa_v(J^{Q_n}(K[m_n]_v))&=H^1_{\text{fl}}(K[m_n]_v,J^{Q_n}[p^M])\subset H^1(K[m_n]_v,J^{Q_n}[p^M])
\end{align*}
where we have extended $v$ to a place of $K[n]$. 
For $\mathfrak F $-many $n $, the restriction of $c_j$ to $K[m_n]_v$ is the image of a Kummer class in $H^1(K[m_n]_v, J^{Q_n}[p^M])$ by a map of Galois representations $J^{Q_n}[p^M]\to T_f/\pi^j \simeq A[\pi^j]$, which extends to a map of finite flat group schemes.
As a consequence, the image of $c_j$ in $H^1(K_v,A)$ is inflated from a class in $H^1(K[m_n]_v/K_v,A)$, which is trivial by \cite[Proposition I.3.8]{milne2006arithmetic}.  (This argument is essentially \cite[Proposition 6.2(1)]{gross1991kolyvagin}.)

If $\mathsf v=\mathsf l|\mathsf m$, then, adopting as well the notation of (\ref{gross subsection}), the class $c(P(m_n), Q_n)$ is zero when restricted to $K[m_n]_{\lambda_n}$ because $D_{\l_n} = \l_n(\l_n+1)$ on $\mathbb F_{\lambda_n}$; hence $\Res_\v\kappa(\mathsf m, \QT) \in \H^1_{\tr} (K_\mathsf l, T_f).$ The local conditions at $\q\in \QT$ are satisfied because every factor of $Q_n$ splits completely in $K[m_n]$; for the same reason, (2, 3) follow from Construction \ref{construction and reciprocity}(2, 3).  (Note that the projection step in Definition \ref{Kolyvagin definition} makes no difference to these identities.) Moreover (1) is clear from Proposition \ref{ Kolyvaginswitcheroo}(3).
\end{proof}
\begin{rmk}\label{ specialization rmk}
The Euler system $(\kappa (1,\cdot),\lambda(1,\cdot))$ may be viewed as a specialization of $(\boldsymbol{\kappa},\boldsymbol{\lambda})$.
Indeed, by the usual Heegner point norm relations \cite[Proposition 3.10]{darmon2004rational},  if $p $ splits in $K $, $\mathbbm 1(\boldsymbol{\lambda}(\QT)) = (\alpha_p-1)^2(\lambda(1, \QT))$ and $\mathbbm 1 (\boldsymbol \kappa(\QT)) = (\alpha_p-1)^2(\kappa(1, \QT))$ when $\set{\QT,\epsilon_\QT} \in \mathsf N^{\nu(N^-)+1}$ and $\set{\QT,\epsilon_\QT} \in \mathsf N^{\nu(N^-)}$, respectively. (Here $\mathbbm 1:\Lambda\to \O$ is specialization at the trivial character.)
\end{rmk}

\section{Deformation theory}\label{deform sec}
Theorem \ref{level raising} allows us to produce weak eigenforms (i.e. ring maps) $\T_{N^+,N^-Q}\to \O/\pi^j$ for arbitrarily large $j $, simply by requiring sufficiently deep congruence conditions on all $q | Q $. However, in general these maps do not lift to characteristic zero. To prove the main results, we also need to be able to $\p$-adically approximate $f$ by genuine level-raised newforms. In this section, we provide this input via the relative deformation theory of Fakhruddin-Khare-Patrikis \cite{fakhruddin2019relative}.
\subsection{Patched adjoint Selmer groups}
\subsubsection{}
Continuing to fix notation as in (\ref{modular forms notation}) and (\ref{construction notation}), we now
assume moreover that $f$ is non-CM.
Consider the (irreducible) adjoint representation $$L =\ad ^ 0 T_f$$ and its $\O $-dual, $L ^\dagger\simeq L (1),$ and let $\overline L$ and  $\overline L^\ast\simeq L ^\dagger/\pi$ be the associated residual representations. For all $v | Np\infty, $
a choice of framing for $T_f$ defines a smooth point of the generic fiber of an appropriate framed universal deformation ring (of fixed determinant, and fixed Hodge type if $v|p$) by \cite[Theorem D]{allen2016deformations}. Taking this smooth point as the input, the construction of \cite[Proposition 4.7]{fakhruddin2019relative} yields, for all $j $,
certain orthogonal local conditions $$ H ^ 1_\mathcal S (\Q_v, L/\pi ^ j)\subset H ^ 1 (\Q_v, L/\pi ^ j),\;\; H^1_{\mathcal S^\ast}(\Q_v, L^\ast[\pi ^ j])\subset H ^ 1(\Q_v, L ^\ast[\pi ^ j]).$$
By \cite[Lemma 6.1]{fakhruddin2019relative}, taking inverse limits yields dual local conditions $$ H ^ 1_\mathcal S (\Q_v, L)\subset H ^ 1 (\Q_v, L),\;\; H^1_{\mathcal S^\dagger}(\Q_v, L^\dagger)\subset H^1(\Q_v, L^\dagger). $$
We use these to define generalized Selmer structures $(\mathcal S,\ST) $ and $(\mathcal S^\dagger, \ST)$ for $L $
and $L^\dagger$, where $\ST\subset \mathsf M_\Q$ is the set of constant ultraprimes $\underline v$ for $v|Np\infty $.

\subsubsection{}
Now suppose that $\q $ is an admissible ultraprime with sign $\epsilon_\q $. Using the exact sequence of $\O[G_\q]$-modules in Definition \ref{admissible def},
$$0\to \Fil ^ +_\q T_f\to T_f\to T_f/\Fil ^ +_\q T_f\to 0, $$
we define $$\Fil ^ +_\q L =\Hom(T_f/\Fil ^ +_\q T_f, \Fil^+_\q T_f) \subset L$$
and $$\H ^ 1_{\ord} (\Q_\q, L) =\im\left (\H ^ 1 (\Q_\q,\Fil ^ +_\q L)\to \H ^ 1 (\Q_\q, L)\right). $$
We also define $\H ^ 1_{\ord} (\Q_\q, L^\dagger)$
as the orthogonal complement of $\H ^ 1_{\ord} (\Q_\q, L)$ under the local Tate pairing; note that, since $\H ^ 1 (\Q_\q, L)$ is torsion-free,  $\H ^ 1_{\ord} (\Q_\q, L^\dagger)$ and   $\H ^ 1_{\ord} (\Q_\q, L)$ are exact annihilators. 
We will require the restriction maps
\begin{equation}
    \begin{split}
        \loc_\q: \H^1(\Q, L) &\to \frac{\H^1(\Q_\q, L)}{\H^1_{\ord} (\Q_\q, L) \intersection \H^1_{\unr} (\Q_\q, L)},\\
        \loc_\q^\dagger: \H^1(\Q, L^\dagger) &\to \frac{\H^1(\Q_\q, L^\dagger)}{\H^1_{\ord} (\Q_\q, L^\dagger) \intersection \H^1_{\unr} (\Q_\q, L^\dagger)}.
    \end{split}
\end{equation}
Analogously, if $q\in M_\Q$ is $j$-admissible with sign $\epsilon_q$, then we may define $H^1_{\ord} (\Q_q, L/\pi^j)$, $H^1_{\ord}(\Q_q, L^\dagger/\pi^j)$, and the localization maps $\loc_q,$ $\loc_q^\dagger$.
\subsubsection{}
For any $\set {\mathsf P\QT\mathsf R,\epsilon_{\mathsf P\QT\mathsf R}}\in\mathsf N $, define the modified Selmer structure $(\mathcal S(\QT), \ST\union\QT)$
for $L$:
\begin{equation}\label{Define all her structures}
    \H ^ 1_{\mathcal S ^\mathsf P_\mathsf R(\QT)}(\Q_\v, L) =\begin {cases} \H ^ 1_\mathcal S (\Q_\v, L), &\v\not\in\mathsf P\QT\mathsf R\\ \H ^ 1_{\ord} (\Q_\q, L), &\v =\q\in\QT,\\
\H ^ 1_{\ord} (\Q_\q, L) +\H ^ 1_{\unr} (\Q_\q, L), &\v =\q\in\mathsf P,\\
\H ^ 1_{\ord} (\Q_\q, L) \intersection\H ^ 1_{\unr} (\Q_\q, L), &\v =\q\in\mathsf R
.\end {cases}
\end{equation} 
The corresponding dual Selmer structure for $L ^\dagger$ will be written $\mathcal S^{\mathsf R,\dagger}_{\mathsf P}(\QT)$. 
Finally, define, for any finite set of places
$\Sigma $ containing all $v|Np\infty $:
\begin{equation}
    \Sha^1_{\Sigma} (\overline L ^\ast) =\ker\left (H ^ 1 (\Q ^ {\Sigma}/\Q,\overline L ^\ast)\to\prod_{v\in \Sigma} H ^ 1 (\Q_v,\overline L ^\ast)\right).
\end{equation}
\begin{prop}\label{balanced}
For all $\set {\QT,\epsilon_\QT}\in\mathsf N $,
$$d_\QT\coloneqq \rank_\O \Sel_{\mathcal S (\QT)} (L) = \rank_\O\Sel_{\mathcal S^\dagger (\QT)} (L ^\dagger). $$
\end{prop}
\begin{proof}
It follows from the construction \cite[Proposition 4.7]{fakhruddin2019relative} that $$\rank_\O H ^ 1_\mathcal S (\Q_v, L) =\rank_\O H ^ 1_{\unr} (\Q_v, L)$$
 for all $v|N$ and $$\rank_\O H ^ 1_\mathcal S (\Q_v, L) =\rank_\O H ^ 0 (\Q_v, L) +2 $$ if $v|p$. Since $$\rank_\O \H ^ 1_{\ord} (\Q_\q, L) =
 \rank_\O \H ^ 0 (\Q_\q, L) = 1 $$
 for  all $\q\in \QT$,
 the claim results from Proposition \ref{Greenberg Wiles over DVR}.
\end{proof}
The ``relative deformation theory'' developed in \cite{fakhruddin2019relative} may be summarized (for our context) as follows.
\begin{thm}[Fakhruddin-Khare-Patrikis, Kisin]\label{relative def}
Suppose given $\set {\QT,\epsilon_\QT}\in\mathsf N $
such that $d_\QT = 0,$
and a finite set of places $\Sigma$
containing all $v|Np\infty$ such that $\Sha^1_{\Sigma} (\Q,\overline L ^\ast) = 0 $. Fix a sequence $\set{Q_n,\epsilon_{Q_n}}$ of weakly admissible pairs representing $\set {\QT,\epsilon_\QT}$ and an integer $j\geq 0. $ Then there is a sequence (defined for $\mathfrak F$-many $n$) of newforms $g_n$  of weight two, level $NQ_n$, and trivial nebentypus, with a prime $\p_{g_n}$ of the ring of integers of its coefficient field $\O_{g_n}$, such that:
\begin{itemize}
    \item The completion $\O_{g_n,\p_{g_n}}$ is isomorphic to $\O $.
    \item The associated Galois representations satisfy $T_f|_{I_\l}\simeq T_{g_n} |_{I_\l} $ for all $\l\nmid Q_n$, and $\rho_{g_n}|_{I_{q_n}}$
    is a Steinberg representation twisted by the unramified character $\Frob_{q_n}\mapsto \epsilon_{Q_n}(q_n)$ for all $q_n|Q_n$.
    \item For any fixed $j$, there is a congruence of Galois representations (in some basis) $$T_f\equiv T_{g_n}\pmod {\pi^j}$$ for $\mathfrak F$-many $n$. In particular, the maps $$\T_{Q_n}=\T_{N^+,N^-Q_n}\to \O/\pi^j$$
of Remark \ref{level raising remark}
admit $\O$-valued lifts for $\mathfrak F$-many $n$.
\end{itemize}
\end{thm}
\begin{proof}
Since $d_\QT = 0,$ Proposition \ref{Selmer proposition investment} implies that there exists some $k\geq 0$ such that the natural maps
\begin{equation}\label{attached relative}
    \Sel_{\mathcal S (\QT)} (L/\pi ^ k)\to\Sel_{\mathcal S  (\QT)} (\overline L),\;\;\Sel_{\mathcal S (\QT)}(L ^\dagger/\pi ^ k)\to\Sel_{\mathcal S(\QT)}(\overline L^\ast)
\end{equation}
are identically zero.

Now, for any $j$ and for $\mathfrak F $-many $n $, we have the weakly admissible pair $\set{Q_n,\epsilon_{Q_n}} \in N_j$
of Remark \ref{level raising remark}.
Consider the (non-patched) Selmer groups
$$\Sel_{Q_n}(L/\pi^k) = \ker\left (H ^ 1 (\Q^{S\cup Q_n}, L/\pi ^ k)\to\prod_{v|Np\infty}\frac{H ^ 1 (\Q_v, L/\pi ^ k)} { H^1_\mathcal S(\Q_v, L/\pi ^ k)}\times\prod_{q_n | Q_n}\frac {H^1(\Q_{q_n}, L/\pi ^ k)} {H ^ 1_{\ord}(\Q_{q_n},L/\pi^k)}\right),$$
where $S$ is the set of places dividing $NpQ_n\infty$.
For each $q_n | Q_n $,  $H^1_{\ord}(\Q_{q_n}, L/\pi^k)$ is exactly the local condition obtained from the smooth, Steinberg-with-sign-$\epsilon_{Q_n}(q_n)$ quotient of the framed local deformation ring at $q_n$ via \cite[Proposition 4.7]{fakhruddin2019relative} (as long as $q_n$ is $k$-admissible). There are also dual Selmer groups $\Sel_{Q_n}(L ^\dagger/\pi ^ k) =\Sel_{Q_n}(L^\ast[\pi^k]) $ defined in the same way (see Proposition \ref{two different duals} for the equality).
By (\ref{attached relative}), for $\mathfrak F $-many $n$, the maps 
$$\Sel_{Q_n} (L/\pi^k) \to \Sel_{Q_n}(\overline L),\;\; \Sel_{Q_n} (L ^\ast[\p ^ k])\to\Sel_{Q_n} (\overline L ^\ast)$$
are identically zero. Moreover, for all $n$, $$\Sha_{\Sigma\cup Q_n}^1 \subset \Sha_\Sigma^1 = 0.$$
The proof of \cite[Claim 6.12]{fakhruddin2019relative} now implies that, for $j\geq 3k$ large enough in a manner depending on the local representations $\rho_f|_{G_{\Q_\l}}$ for $\l|Np$ and for $\mathfrak F$-many $n$, there exists a representation $$\tau_n: G_\Q \to GL_2 (\O) $$
such that:
\begin{itemize}
    \item $\tau_n\equiv \rho_f \pmod {\pi^j}$ (for some choice of basis of $T_f$);
    \item $\det\tau_n = \chi$;
    \item the local representations $\tau_n|_{G_{\Q_\l}}$ lie on the same irreducible component of the framed local deformation ring as $\rho_f|_{G_{\Q_\l}}$ if $\l\nmid Q_n$, and are Steinberg representations
    twisted by the unramified character $\Frob_{q_n}\mapsto \epsilon_{Q_n}(q_n)$ for all $q_n|Q_n$.
\end{itemize}
It remains to apply a modularity lifting theorem to conclude that $\tau_n$ arises from a suitable modular form $g_n$. We claim that  $\overline T_f |_{G_{\Q(\zeta_p)}}$ 
is absolutely irreducible: suppose otherwise for contradiction. Then the image $G$ of the action of $G_\Q $ on $\overline T_f $ fixes a pair of lines. Since the image of the inertia group $I_p$ contains a matrix whose square has distinct eigenvalues by (\ref{Taylor Wiles}), the pair of lines must be the eigenspaces of $I_p$. But $I_p$ surjects onto $\Gal(\Q(\zeta_p)/\Q)$, so by assumption no element of $G_\Q$ interchanges the two eigenspaces of $I_p$, and $\overline T_f$ is reducible -- a contradiction.
Thus the Taylor-Wiles hypothesis in Kisin's result \cite{kisin2009moduli} is satisfied, and $\tau_n$ indeed arises from a newform $g_n$.
\end{proof}
\subsection{Annihilating two Selmer groups}
\subsubsection{}
In order to apply Theorem \ref{relative def}, we must make a suitable choice of $\set {\QT,\epsilon_\QT}\in\mathsf N $. In this  subsection, we show that such a choice as possible.

\begin{prop}\label{killing sha}
There exists a finite set of places $\Sigma$, containing all $v|Np\infty$,
such that $$\Sha_{\Sigma}^1 = 0.$$
\end{prop}
\begin{proof}
It suffices to show that \begin{equation}\label{Inflation goal}H^1(\Q(\overline L ^\ast)/\Q,\overline L ^\ast)= 0,\end{equation}
for if so,
 for any $0\neq c\in H ^ 1 (\Q,\overline L ^\ast)$, $c$ restricts to a nonzero homomorphism $G_{\Q(\overline L^\ast)} \to \overline L ^\ast, $
 and there exist primes which are totally split in $\Q(\overline L^\ast)$ but not in the extension cut out by the restriction of $c$.
 
 We now show (\ref{Inflation goal}). If $\mu_p\not\subset \Q(\overline L) $,
 then the center of $\Gal (\Q (\overline L ^\ast)/\Q) $ contains elements that act by nontrivial scalars on $\overline L ^\ast, $ and (\ref{Inflation goal}) follows from inflation-restriction. So suppose that $\mu_p\subset \Q(\overline L)$; then the projective image $\overline G = \Gal(\Q (\overline L)/\Q) $ of the (irreducible) Galois action on $\overline T_f$ has a cyclic quotient of order $p -1, $
 and a classical result of Dickson implies that $p = 3$ and $\overline G$ is either a dihedral group, or $S_4.$ In the former case, the order of $\Gal(\Q(\overline L^\ast)/\Q)$ is prime to $p$, so (\ref{Inflation goal}) still holds. We are left to consider the case $\overline G = S_4$ and $p = 3.$ Let $G= \Gal(\Q(\overline T_f)/\Q) $
 be the image of the Galois action; since we have assumed that $\det: G\to\F_3 ^\times $
 factors through $\overline G $, a complex conjugation $c$ in $G$
 projects to a transposition in $\overline G$. Let $\overline H\subset \overline G$ be a copy of $S_3$ containing the image of $c$, and $H$ the normalizer of its preimage in $G$, which is contained in a unique Borel subgroup  $B$. Let $N$ be the unipotent radical of $B\intersection G$.
  To prove (\ref{Inflation goal}), it suffices to check that $$H ^ 1 (H, \overline L ^\ast) = H^1(N, \overline L^\ast)^{H/N}= 0. $$
  This holds because
 $\im  (N-1)$ is isomorphic to the subgroup $\begin{pmatrix} \ast & \ast \\ 0 & \ast\end{pmatrix}\subset \overline L^\ast$, while $c\in H$ acts on $N$ by $-1$ and on $\overline L^\ast / \im  (N-1)$ by $1$.
 \end{proof}
We now require a more elaborate version of Theorem \ref{localization nontrivial}; the proof is inspired by \cite{castella2020anticyclotomic}, and begins with a series of lemmas.
 \begin{lemma}\label{first constant}
 There exists an integer $j $  that, for all $n\geq 0, $
$$ \pi ^ j H ^ 1 (K (T_f)/\Q, L/\pi ^ n) =\pi ^ j H ^ 1 (K(T_f)/\Q, L ^\dagger/\pi ^ n) = 0. $$
\end{lemma}
\begin {proof}
Let $E= \Q (\mu_{p ^\infty})\subset K(T_f)$, and  note that $L$ and $L^\dagger$ are isomorphic $G_E$-modules.
Since $f$ is non-CM, $(L\otimes \Q_p)^{G_E} = 0$, and so  $(L/\pi^n)^{G_E}$ is uniformly bounded in $n$. 

The pro-$p$-Sylow subgroup of $\Gal(K(T_f)/E)$ is a compact $p$-adic Lie group with semisimple Lie algebra; hence, by \cite[Lemma B.1]{fakhruddin2019relative}, $H^1(K(T_f)/E, L/\pi^n)$ 
is uniformly bounded in $n$. 

Now, by inflation-restriction,  we have exact sequences
\begin{equation}
\begin{split}
   & 0\to H^1(E/\Q, (L/\pi^n)^{G_E})\to H^1(K(T_f)/\Q, L/\pi^n) \to H^1(K(T_f)/E, L/\pi^n), \\
   &  0\to H^1(E/\Q, (L^\dagger/\pi^n)^{G_E})\to H^1(K(T_f)/\Q, L^\dagger/\pi^n) \to H^1(K(T_f)/E, L^\dagger/\pi^n),
    \end{split}
\end{equation}
where the outer terms are isomorphic and uniformly bounded in $n$; the lemma follows.
\end{proof}
For the next lemma, we abbreviate $L_m\coloneqq L/\pi^m$, $L^\dagger_m\coloneqq L^\dagger/\pi^m$, and $T_m \coloneqq T_f/\pi^m$. Moreover, if $y\in M$ for any torsion $\O$-module $M$, let $\ord(y)$ be the smallest integer $t\geq 0$ such that $\pi^t y = 0 $. 
\begin{lemma}\label{ugly lemma}
There is a global constant $C$, depending on $T_f$, with the following property.
  Given cocycles $\xphi\in H^1(\Q, L_m)$, $\psi\in H^1(\Q,L^\dagger_m),$ and $c_1,c_2\in H^1(K, T_m)^\delta$ for some $\delta = \pm 1,$ there exist infinitely many primes $q\nmid Np$ such that all the cocycles are unramified at $q$ and:
  \begin{itemize}
      \item The Frobenius of $q$ in $\Gal(K(T_m)/\Q)$ is a complex conjugation; in particular, $q$ is $m$-admissible with sign $\delta$.
      \item $\ord \loc_q \xphi \geq \ord \xphi - C.$
      \item $\ord \loc_q^\dagger \psi\geq \ord \psi - C$, or $\loc_q^\dagger \psi = 0$, as desired.
      \item $\ord \loc_q c_i \geq \ord c_i - C$ for $i = 1,2.$
  \end{itemize}
\end{lemma}
\begin{proof}
Let us first fix a complex conjugation $c\in G_\Q$ and choose a basis for $T_m$ in which $c$ acts as $\begin{pmatrix} -\delta & 0 \\ 0 &\delta\end{pmatrix}$.

The restriction of the cocycles $\xphi, \psi, c_i$ to $G_{K(T_m)}$ may be considered as a homomorphism $$h: G_{K(T_m)} \to L_m\oplus L^\dagger_m \oplus (T_m)^2$$ compatible with the action of $G_K$; let $H$ be the image of this homomorphism. Since there exists an element of $g_z\in G_K$ that acts by a scalar $z\neq \pm 1$ on $T_f$, we have:\begin{align*}
 H \supset &(g_z - z)(g_z - z^2) H + (g_z - z)(g_z -1) H + (g_z - z^2)(g_z - 1) H\\
    & \supset (z - 1)(z^2 - 1)(z^2 - z) \left(\pi_{L_m}(H) \oplus \pi_{L_m^\dagger} (H) \oplus \pi_{T_m^2}(H),\right)
\end{align*}
where $\pi_\bullet$ are the projection operators. Now, since $L$ and $L^\dagger$ are absolutely irreducible, the natural maps $\Q_p[G_K]\to \End(L\otimes \Q_p)$ and $\Q_p[G_K]\to \End(L^\dagger \otimes \Q_p)$ are surjective. Combining these observations with Lemma \ref{first constant}, we see that, for some constant $C$ depending only on $T_f$, there exists $\gamma\in G_{K(T_f)}$ satisfying:
\begin{itemize}
    \item The $\begin{pmatrix} 1 & 0 \\ 0 & -1 \end{pmatrix}$ component of $\xphi(\gamma)$ has order at least ${\ord \xphi - C}$.
    \item The $\begin{pmatrix} 0 & 0 \\ 1 & 0 \end{pmatrix}$ component of $\psi(\gamma_\psi)$  has order at least ${\ord \psi - C}$, 
    or is 0, as desired.
    \item The components of $c_i(\gamma)$ and $c_2(\gamma)$ in the $\delta$ eigenspace have order at least $\ord c_i - C$, where $i=1,2$.
\end{itemize}
For the final item, we are using the elementary fact that a group cannot be the union of two trivial subgroups, as well as the irreducibility of $T_f$. 

Since $\xphi(c^2) = c\xphi(c) + \xphi(c) = 0,$ $\xphi(c)$ lies in the $-1$ eigenspace for complex conjugation, whereas $\begin{pmatrix} 1 & 0 \\ 0 & -1 \end{pmatrix}$ has eigenvalue 1;
hence the $\begin{pmatrix} 1 & 0 \\ 0 & -1 \end{pmatrix}$ component of $\xphi(c\gamma)$ has order at least ${\ord \xphi - C}$. Similarly, the $\begin{pmatrix} 0 & 0 \\ 1 & 0 \end{pmatrix}$ component of $\psi(c\gamma_\psi)$ has order at least ${\ord \psi - C}$, or is 0, as desired.
    
    Any prime with Frobenius $c\gamma$ in $\ker h$ satisfies the conclusion of the lemma; cf. the proof of Lemma \ref{just one localization} for the assertions about $c_i$.
\end{proof}
\begin{cor}
 Suppose given a finite set of ultraprimes $\mathsf T$ and non-torsion cocycles:
 \begin{itemize}
     \item $\xphi\in\H ^ 1 (\Q^\mathsf T/\Q, L) $;
     \item $\psi\in\H ^ 1 (\Q ^\mathsf T/\Q, L ^\dagger) $;
     \item $c_1, c_2\in\H ^ 1 (K ^\mathsf T/K, T_f)^\delta $
     for $\delta = \pm 1.$
 \end{itemize}
 Then there exist infinitely many admissible ultraprimes $\q\not\in \mathsf T$ with sign $\delta$ such that:
 \begin{itemize}
     \item $\loc_\q \xphi\neq 0 $.
     \item Either $\loc^\dagger_\q\psi \neq 0 $ or $\loc_\q^\dagger\psi =0$, as desired.
     \item $\loc_\q c_i \neq 0.$
 \end{itemize} 
\end{cor}
\begin{proof}
Choose  a sequence $T_n$ representing $\mathsf T$ and sequences $\xphi_n,\psi_n,c_n^1, c_n^2$  representing the respective cocycles in $H^1(\Q^{T_n}/\Q, L/\pi^n)$, etc. For each $n$, apply Lemma \ref{ugly lemma} with $m = n$ and the appropriate desideratum for $\psi_n$; by definition,  any resulting admissible ultraprime $\q$, represented by a sequence $q_n$, satisfies the desired conclusion.
\end{proof}
\begin{prop}\label{killing Selmer}
Suppose given a self-dual Selmer structure $(\mathcal F ,\mathsf T)$
for $T_f $. Then there exists $\set{\QT,\epsilon_{\QT}}\in \mathsf N_\mathsf T $
such that $$r_{\QT} = d_{\QT} = 0. $$
(Recall that $r_{\QT} =\rank_\O\Sel_{\mathcal F (\QT)} (T_f).$)
\end{prop}
\begin{proof}
Without loss of generality, by Corollary \ref{the Selmer annihilation corollary} we may assume that $r_1 = 0$; for if not, choose any $\set{\QT,\epsilon_{\QT}}\in\mathsf N_\mathsf T $ with $r_ \QT = 0, $ and then relabel
 $\mathcal F(\QT)$ as $\mathcal F$.

We will show that, if $d_1 > 0$,  we may find $\set{\QT,\epsilon_{\QT}}$ such that 
$r_\QT= 0 $ and $d_\QT < d_1$; this clearly suffices by induction. 
By Proposition \ref{balanced}, there exist non-torsion elements $\xphi\in \Sel_{\mathcal S} (L),\;\psi\in\Sel_{\mathcal S ^\dagger}(L ^\dagger).$ We choose any admissible $\q\not\in \mathsf T$ with sign $\epsilon_\q$ such that $\loc_\q \xphi \neq 0,$ $\loc_\q^\dagger \psi \neq 0.$ 
Then by Proposition \ref{Greenberg Wiles over DVR},
$$\rank_\O \Sel_{\mathcal S^{\dagger, \q}}(L^\dagger) + \rank_\O \Sel_{\mathcal S^{\q}}(L) = 2 + \rank_\O \Sel_{\mathcal S^\dagger_\q}(L^\dagger) + \rank_\O \Sel_{\mathcal S_\q} (L)$$
(in the notation of (\ref{Define all her structures})).
The images of the localization maps
$$\loc_\q: \frac{\Sel_{\mathcal S^{\q}}(L)}{\Sel_{\mathcal S_\q} (L)} \hookrightarrow \frac{\H^1_{\ord} (\Q_\q, L)}{\H^1_{\ord \intersection \unr} (\Q_\q, L)} \oplus \frac{ \H^1_{\unr}(\Q_\q, L)}{\H^1_{\ord \intersection\unr} (\Q_\q, L)}$$
and 
$$\loc^\dagger_\q: \frac{\Sel_{\mathcal S^{^\dagger,\q}}(L^\dagger)}{\Sel_{\mathcal S^\dagger_\q} (L^\dagger)} \hookrightarrow \frac{\H^1_{\ord} (\Q_\q, L^\dagger)}{\H^1_{\ord \intersection \unr} (\Q_\q, L^\dagger)} \oplus \frac{ \H^1_{\unr}(\Q_\q, L^\dagger)}{\H^1_{\ord \intersection\unr} (\Q_\q, L^\dagger)}$$
have total rank two and annihilate each other under the induced Tate pairing by Proposition \ref{ dualityover DVR}. Hence the image in the ordinary part is zero for both maps, and $d_\q < d_1.$
However, by adding $\q$, we have made $r_\q = 1.$  Let $c\in \Sel_{\mathcal F (\q)} (T_f) $
be a generator; since $\partial_\q c \neq 0$ by Proposition \ref{Ray duality}, $c$ has nonzero component in the $\epsilon_\q$ eigenspace for $\tau$.

Now consider the set $\mathsf P$ of admissible ultraprimes $\mathsf s$ with sign $\epsilon_\s = \epsilon_\q$ such that $\loc_\s c \neq 0.$ If, for any $\s\in \mathsf P$, $d_{\q\s} \leq d_\q$, then we may take $\QT = \q\s$ and complete our induction step. For example, this will occur provided $d_\q > 0,$ by the argument above; so without loss of generality, $d_\q = 0$ and $d_{\q\s} =1 $ for all $\s \in \mathsf P$.
By definition, we therefore have non-torsion elements $\xphi(\s)\in \Sel_{\mathcal S (\q\s)} (L) $
and $\psi (\s) \in \Sel_{\mathcal S ^\dagger (\q\s)} (L ^\dagger) $
such that $\loc_\s \xphi(s)$ and $\loc_\s \psi(s)$ do not lie in the unramified subspace of the ordinary cohomology. 

Choose any $\s_1\in\mathsf P$, and then choose $\s_2\in \mathsf P$ such that $\loc_{\s_2}\xphi(\s_1) \neq 0$ but $\loc_{\s_2}\psi(\s_1) = 0.$ By another application of Proposition \ref{Ray duality}, $r_{\q\s_1\s_2} = 1,$ and a generator $c'$ of $\Sel_{\mathcal F(\q\s_1\s_2)}(T_f)$ again has nonzero component in the $\epsilon_\q$ eigenspace.
We now choose $\s_3\in \mathsf P$ such that $\loc_{\s_3} c'\neq 0$, $\loc_{\s_3} \xphi(\s_2)\neq 0,$ and $\loc_{\s_3}\psi(\s_1) \neq 0.$ Clearly $r_{\q\s_1\s_2\s_3} = 0.$
Note that
$\rank_\O \Sel_{\mathcal S^{\s_1\s_2\s_3}(\q)} = 3;$ up to torsion, $\xphi(\s_i)$ are generators. So to show that $d_{\q\s_1\s_2\s_3} = d_\q$, it suffices to show that the images of $\xphi(\s_i)$ form a rank-three subspace of $$S\coloneqq \bigoplus_{i = 1}^3 \frac{\H^1_{\unr+\ord}(\Q_{\s_i}, L)}{\H^1_{\ord}(\Q_{\s_i}, L)}$$ under the localization $$\loc:\frac{ \Sel_{\mathcal S^{\s_1\s_2\s_3}}(\q)(L)}{\Sel_{\mathcal S_{\s_1\s_2\s_3}(\q)}(L)}\hookrightarrow S.$$ By pairing $\xphi(\s_i) $
and $\psi (\s_j) $
for $i\neq j $
and applying Proposition \ref{ dualityover DVR} once more,
we see that  $\loc_{\s_i}\xphi(\s_j)\neq 0$ if and only if $\loc_{\s_j} \psi(\s_i)\neq 0.$ Hence, the images of $\xphi(\s_i)$ in $S$ are of the form:\begin{align*}
    \loc(\xphi(\s_1)) &= (0, \ast,\cdot) \\
    \loc(\psi(\s_2)) &= (0, 0, \ast),\\
    \loc(\psi(\s_3)) &= (\ast, 0, 0),
\end{align*}
where $\ast$ is nonzero and $\cdot$ may or may not be zero. This completes the inductive step since $d_{\q\s_1\s_2\s_3} = d_\q < d_1.$
\end{proof}
 \section{Proof of main results}\label{proof section}
For this section, let $f $
be a modular form of weight two, level $N $, and trivial character, with ring of integers $\O_f$, and let $\p\subset\O_f$ be an ordinary prime lying over $p\nmid 2N.$ Denote by $\O$ the completion. 

\subsection{A result of Skinner-Urban}

The following result is a corollary to the proof of the main conjecture \cite{skinner2014iwasawa}.
\begin{thm}[Skinner-Urban]\label{skinner urban}
Let $K $
be an imaginary quadratic field of discriminant prime to $Np$ in which $p$ splits. Assume that $\p$ is ordinary for $f$ and that:
\begin{itemize}
    \item the mod $\p$ representation  $\overline T_f$ is absolutely irreducible; 
    \item $N = N_1N_2,$ where every factor of $N_1 $ is split in $K$ and $N_2$ is the squarefree product of an odd number primes inert in $K $.
    \end{itemize}
        If $\Sel_{\mathcal F_\Lambda}(\mathbf W_f)$ is $\Lambda$-cotorsion, then $$\Char_\Lambda\Sel_{\mathcal F_\Lambda} (\mathbf W_f)^\vee \subset (\boldsymbol\lambda(1))^2$$
as ideals of $\Lambda $,
where  $\boldsymbol{\lambda}(1) \in \Lambda$ is the element constructed in (\ref{p adic interpolation}).

\end{thm}
\begin{proof}
We must explain some details and notations of \cite{skinner2014iwasawa}, in which it is assumed that $\overline T_f$ is ramified at every $\l|N_2.$ 
As in \cite{skinner2014iwasawa}, we let $O_L$ be the ring of integers of a suitable finite extension of $\mathbb{Q}_p$ and consider $f$ as a specialization of a suitable Hida family $\mathbf{f}$.  This family is parametrized by   $\mathbb{I}, $ which is a normal domain and a finite integral extension of $O_L\llbracket W\rrbracket$.  We write $\Gamma_{{K}}=\Gamma_{{K}}^+\times\Gamma_{{K}}^-$ for the Galois group of the maximal $\Z_p$-extension of ${K}$ and its decomposition into cyclotomic/anticyclotomic components.  For a sufficiently large finite set of primes $\Sigma$, there is \cite[Theorem 12.3.1]{skinner2014iwasawa} a three-variable $p$-adic $L$-function $\mathcal{L}_{\mathbf{f},{K}}^\Sigma\in \mathbb{I} \llbracket\Gamma_{{K}}\rrbracket.$ (Here the superscript $\Sigma $ refers to removing Euler factors at primes in $\Sigma$, or relaxing local conditions for a Selmer group.)
Letting $\gamma^-$ be a topological generator of $\Gamma_{{K}}^-$, we may expand: 
\begin{equation}\label{expand}
 \mathcal{L}_{\mathbf{f},{K}}^\Sigma=a_0+a_1(\gamma^--1)+a_2(\gamma^--1)^2+\ldots
\end{equation}
 where $a_i\in \mathbb{I} \llbracket\Gamma_{{K}}^+\rrbracket$.  Let $Ch_{{K}_\infty}^\Sigma(\mathbf{f})\subset\mathbb{I}\llbracket\Gamma_{{K}}\rrbracket$ be the characteristic ideal of the three-variable Selmer group as considered in \cite{skinner2014iwasawa}. Skinner and Urban deduce 
 \begin{equation}\label{imcdiv}
 Ch_{K_\infty}^\Sigma(\mathbf{f})\subset
 (\mathcal{L}_{\mathbf{f},K}^\Sigma)
 \end{equation}
 by proving (see \cite[Theorem 6.5.4, Proposition 12.3.6, Proposition 13.4.1]{skinner2014iwasawa}):
\begin{enumerate}
    \item If $P\subset\mathbb{I}\llbracket\Gamma_{K}\rrbracket$ is a height one prime which is not of the form $P_+\mathbb{I}\llbracket\Gamma_{K}\rrbracket$ for some $P_+\subset\mathbb{I}\llbracket\Gamma_{K}^+\rrbracket$, then $$\ord_{P}Ch_{K_\infty}^\Sigma(\mathbf{f})\geq\ord_P(\mathcal{L}_{\mathbf{f},K}^\Sigma).$$
    \item If $\overline T_f $
    is ramified at every $\l | N_2, $ then $\ord_P(\mathcal{L}_{\mathbf{f},K}^\Sigma)=0$ for all height one primes $P$ of the form $P_+\mathbb{I}\llbracket\Gamma_{K}\rrbracket$ for some $P_+\subset\mathbb{I}\llbracket\Gamma_{K}^+\rrbracket$.
\end{enumerate}
Although (2) does not apply, we claim that we may replace (\ref{imcdiv}) by the weaker inclusion: 
\begin{equation}\label{weakimcdiv}
    Ch_{K_\infty}^\Sigma(\mathbf{f})\cdot (a_i)\supset(\mathcal{L}_{\mathbf{f},K}^\Sigma)
\end{equation} where $a_i$ is any of the terms in (\ref{expand}). Indeed,
because both sides of (\ref{weakimcdiv})
are divisorial, it suffices to check that $\ord_{P} (\mathcal{L}_{\mathbf{f},K}^\Sigma)\leq\ord_P (a_i) $
for all $P $
as in (2). But this is clear: if $(\mathcal{L}_{\mathbf{f},K}^\Sigma) $
is zero modulo $P ^ k $
for such a prime $P $, then $a_i $
is as well.
By \cite{vatsal2003special} $a_i$ may be chosen so that its image under the specialization map $\mathbbm{1}: \mathbb{I}\llbracket\Gamma_{K}^+\rrbracket\to O_L$ is nonzero. Fix such a choice $\widetilde\alpha$.

The divisibility (\ref{weakimcdiv})  also (trivially) implies a divisibility for the Fitting ideal of the 3-variable Selmer group:
\begin{equation}\label{Sitting divisible}
   (\widetilde\alpha) Fitt_{K_\infty}^\Sigma(\mathbf g)\subset (\mathcal L_{\mathbf g, K} ^\Sigma).
\end{equation}


 Specializing (\ref{Sitting divisible}) to the anticyclotomic variable, we obtain $$\Char_\Lambda \Sel^\Sigma(K_\infty, f)\subset L_p^\Sigma(K_\infty, f)\text{ in }\Lambda\otimes\Q_p, $$
 where $L_p^\Sigma(K_\infty, f)$ is a certain $\Sigma$-primitive anticyclotomic $L$-function, and $\Sel^\Sigma(K_\infty, f)$ is the $\Sigma$-primitive Selmer group.  Replacing \cite[Proposition 3.3.19]{skinner2014iwasawa} by  \cite[Proposition A.2]{pollack2011anticyclotomic} (and using the hypothesis that the Selmer group is $\Lambda $-cotorsion), we may convert this to an imprimitive divisibility
 \begin{equation}\label{with p inverted}
     \Char_\Lambda \Sel_{\mathcal F_\Lambda}(\mathbf W_f)^\vee \subset (\boldsymbol \lambda (1)) ^ 2\text{ in }\Lambda\otimes\Q_p. \end{equation}
NB: The anticyclotomic $p$-adic $L$-function appearing 
 in \cite[\S12.3.5]{skinner2014iwasawa}, which emerges naturally from the specialization of the three-variable $p$-adic $L$-function, is normalized using Hida's canonical period, whereas $\boldsymbol{\lambda}(1)^2$ is the $L$-function constructed in \cite{bertolini2005iwasawa}, normalized using Gross's period. However, these $L$-functions differ only by a power of $p$. Similarly, the local cohomology groups $H^1(K_\l, \mathbf W_f)$ for $\l|N_2$ have characteristic ideal a power of $(\p)$, so the choice of local condition at primes $\l|N_2$ does not change the characteristic ideal in $\Lambda\otimes\Q_p $. See the appendix, and \cite{pollack2011anticyclotomic} for a detailed discussion.
 
 To upgrade (\ref{with p inverted}) to a divisibility in $\Lambda $, we simply note that  the $\mu$-invariant of $\boldsymbol{\lambda}(1)$ is 0 by \cite{vatsal2003special}. 
\end{proof}

\subsection{The Heegner point main conjecture}
In this subsection, we prove the following main theorem.
\begin{thm}\label{Heegner point main conjecture}
Let $f$ be a modular form of weight two, level $N$, and trivial character, with an ordinary prime $\p$ of its ring of integers $\O_f$, and let $K$ be an imaginary quadratic field. Assume:
\begin{itemize}
    \item $N = N^+N^-$, where every factor of $N^+$ is split in $K $, and $N ^ - $
is a squarefree product of 
 primes inert in $K $.
 \item The residue characteristic $p$ of $\p $ does not divide $2D_KN$, and $p $ splits in $K $.
 \item The modulo $\p$ representation $\overline T_f $
associated to $f $ 
is absolutely irreducible; if $p = 3, $ assume that $\overline T_f $ is not induced from a character of  $G_{\Q\sqrt{-3}}$.
\end{itemize}
Then, for all $\set{\QT,\epsilon_\QT}\in\mathsf N^{\nu(N^-)}$ such that $(\boldsymbol{\kappa}(\QT,\epsilon_\QT))\neq 0,$
we have  $$\rank_\Lambda \Sel_{\mathcal F_\Lambda (\QT)} (\mathbf T) = \corank_\Lambda \Sel_{\mathcal F_\Lambda (\QT)} (\mathbf W_f) =1 $$ and
$$\Char_\Lambda \left(\left(\Sel_{\mathcal F_\Lambda (\QT)} (\mathbf W_f)^\vee\right)_{\tors}\right)  =\Char_\Lambda \left(\frac{\Sel_{\mathcal F_\Lambda (\QT)} (\mathbf T)} {(\boldsymbol\kappa (\QT))}\right)^2\text{   in }\Lambda\otimes \Q_p. $$
 For all $\set{\QT,\epsilon_\QT}\in\mathsf N ^ {\nu(N^-) +1}$ such that $\boldsymbol\lambda (\QT)\neq 0, $  $$\rank_\Lambda \Sel_{\mathcal F_\Lambda (\QT)} (\mathbf T) = \corank_\Lambda \Sel_{\mathcal F_\Lambda (\QT)} (\mathbf W_f) =0 $$
 and $$\Char_\Lambda \left(\Sel_{\mathcal F_\Lambda (\QT)} (\mathbf W_f)^\vee\right) = (\boldsymbol\lambda (\QT))^2 \text{   in } \Lambda\otimes \Q_p. $$
If moreover the image of the $G_\Q$ action on $\overline T_f $
contains a nontrivial scalar, then the equalities hold in $\Lambda $.
 \end{thm} 
 \begin{proof}
 Given $f $, apply Proposition \ref{killing Selmer} to the standard Selmer structure $(\mathcal F,\ST)$ on $T_f$ (with local conditions the image of the Kummer map at all $\underline v$ such that $v|Np$).
 Let $\set{Q_n,\epsilon_{Q_n}}$
 be a sequence of weakly admissible pairs representing the resulting 
 pair $\set{\QT, \epsilon_\QT}\in \mathsf N$. Let $g_n$ be the resulting sequence of newforms of level $NQ_n$ obtained from Theorem \ref{relative def} (and Proposition \ref{killing sha}); $g_n$ may only be defined for $\mathfrak F$-many $n$. 
 \setcounter{step}{0}
 \begin{step}
 $\set{\QT,\epsilon_\QT}\in \mathsf N^{\nu(N^-) + 1}$.
 \end{step}
 \begin{proof}
 The prime $(T)$, corresponding to the trivial character, does not lie in the exceptional set $\Sigma $ for $\mathcal F_\Lambda$ (see the proofs of \cite[Lemma 5.3.13]{mazur2004kolyvagin} and \cite[Lemma 2.2.7]{howard2004heegner}). Hence $\Sel_{\mathcal F_\Lambda(\QT)} (\mathbf T)= 0,$ which by Theorem \ref{ lambda Euler system bound} and the nontriviality of $(\boldsymbol\kappa,\boldsymbol\lambda) $
 implies the claim.
 \end{proof}

 \begin{step}
 For any fixed $j$, $$(\boldsymbol\lambda (\QT)) \equiv (\boldsymbol\lambda_{g_n}(1))\pmod{\p^j,T^j}$$
    for $\mathfrak F$-many $n$.
 \end{step}
 \begin{proof}
 Recall notations of (\ref{ level raising section}) and (\ref{p adic interpolation}).
 By definition, the image of $\boldsymbol\lambda (\QT)$ modulo $(\p^j,T^j)$ 
is a map $\Gal(K_j/K)\to \O$ obtained, for $\mathfrak F$-many $n$,
by evaluating a surjective map $F_n: M_{Q_n}\otimes_{\T^{Q_n}} \O(f)\to \O(f)/\p^j$ of $\T^{Q_n}$-modules at certain CM points, where $\O(f)$ is $\O$ with $\T^{Q_n}$-action by $f$. Recall that the map is chosen to  factor through multiplicatation by $\O(f)/\pi^{j+C}$ for the constant $C$ of Lemma \ref{approximation Lemma}, and that $\pi^C(M_{Q_n}\otimes_{\T^{Q_n}}\O(f))$ is principal. When $g_n $
has a sufficiently deep congruence to $f $,
$\O(g_n)/\pi^{j+C} = \O(f)/\pi^{j+C}$ as $\T ^ {Q_n} $-modules, and
the composite $G_n: M_{Q_n}\to\O(g_n)\to\O (g_n)/\pi^j$ induces a unit multiple of $F_n$, where $M_{Q_n}\to \O(g_n)$ is the quaternionic modular form associated to $g_n $. But $G_n$ is the very map whose evaluation at CM points  is used to define $\boldsymbol\lambda_{g_n}(1)$, and the claim follows.
 \end{proof}
  \begin{step}
 For any fixed $j$,
    $$\Fitt_\Lambda \Sel_{\mathcal F_\Lambda (\QT)} (\mathbf W_f) ^\vee\equiv    \Fitt_\Lambda\Sel_{\mathcal F_{g_n,\Lambda}} (\mathbf W_{g_n}) ^\vee\pmod{\p^j,T^j}$$
    for $\mathfrak F$-many $n$.   
 \end{step}
 \begin{proof}
 Since fitting ideals are stable under base change and $\overline T_f $
 has no $G_K $-fixed points, it suffices to show that
 \begin{equation}
    \Fitt_\Lambda\left ( \Sel_{\mathcal F_\Lambda (\QT)} (\mathbf W_f[\pi^j,T^j])\right)=\Fitt_\Lambda\left (\Sel_{\mathcal F_{g_n,\Lambda}} (\mathbf W_{g_n}[\pi ^ j, T ^ j])\right)
 \end{equation}
 for $\mathfrak F$-many $n$. Note that, for $\mathfrak F $-many $n $, $\mathbf W_f [\pi ^ j, T^j] =\mathbf W_{g_n}[\pi^j, T^j]$ as finite Galois modules, and $$\Sel_{\mathcal F_\Lambda (\QT)} (\mathbf W_f[\pi^j,T^j])$$ is isomorphic to a submodule of $H^1(K^{\Sigma\union Q_n}/K, \mathbf W_f[\pi ^ j, T ^ j]) $
 defined by certain local conditions. We will show that these local conditions coincide with the ones defining $\Sel_{\mathcal F_{g_n,\Lambda}} (\mathbf W_{g_n}[\pi ^ j, T ^ j]) $. At a $j$-admissible prime $q_n|Q_n$, this is clear (cf. the proof of Construction \ref{construction and reciprocity}(1)). 
 At $v|N$, the local conditions are simply the kernels $$\ker \left(H ^ 1 (K_v,\mathbf W_f[\pi ^ j, T^j])\to H ^ 1 (K_v,\mathbf W_f)\right), \;\;\ker\left (H^ 1 (K_v,\mathbf W_{g_n} [\pi ^ j, T ^ j])\to H ^ 1 (K_v,\mathbf W_{g_n})\right). $$
If $v | N $ is a prime of multiplicative reduction for $f $, then $\mathbf W_f =\mathbf W_{g_n}$ as $G_{K_v}$ modules for $\mathfrak F $-many $n$, so the local conditions clearly coincide. For other places of bad reduction,
the inertia co-invariants $ T_{f, I_v}$ and $T_{g_n, I_v}$ are finite, and we may assume they are isomorphic, say with exponent bounded by $\pi^{M-j}$ for some $M\geq 0$. It follows
that the local conditions are also given by the kernels of $$\ker \left(H ^ 1 (K_v,\mathbf W_f[\pi ^ j, T^j])\to H ^ 1 (K_v,\mathbf W_f[\pi^M])\right), \;\;\ker\left (H^ 1 (K_v,\mathbf W_{g_n} [\pi ^ j, T^j])\to H ^ 1 (K_v,\mathbf W_{g_n}[\pi^M])\right),$$
which  also agree for $\mathfrak F $-many $n $.

For primes $v|p$, it suffices to compare the kernels  $$\left(H ^ 1 (K_v,\gr \mathbf W_f[\pi ^ j, T^j])\to H ^ 1 (K_v,\gr\mathbf W_f)\right), \;\;\ker\left (H^ 1 (K_v,\gr\mathbf W_{g_n} [\pi ^ j, T^j])\to H ^ 1 (K_v,\gr\mathbf W_{g_n})\right).$$
 A similar argument as above applies provided that $H^0(K_v, \gr\mathbf W_f)= \prod_{w|v}H^0(K_{\infty, w}, \gr W_f) $ is finite, which it  is because $a_p$ cannot be a root of unity.
 \end{proof}
 \begin{step}
 Conclusion of the proof.
 \end{step}
Step 1 shows that $N ^ - Q_n$
  is the squarefree product of an odd number primes inert in $K $ for $\mathfrak F$-many $n$, and by Theorem \ref{ lambda Euler system bound} applied to the Euler system $(\boldsymbol \kappa_{g_n}, \boldsymbol\lambda_{g_n})$ for $T_{g_n}$, $\Sel_{\mathcal F_{g_n,\Lambda}} (\mathbf W_{g_n})$ is then $\Lambda $-cotorsion.  
 By Theorem \ref{skinner urban}, for such  $n$ we have: 
 \begin{equation}\label{divisibility work planet levels}
     \Fitt_\Lambda\Sel_{\mathcal F_{g_n,\Lambda}} (\mathbf W_{g_n}) ^\vee\subset (\boldsymbol\lambda_{g_n} (1)) ^ 2\subset\Lambda.
 \end{equation}
 On the other hand, by Theorem \ref{ lambda Euler system bound}, the theorem would follow from:
 \begin{equation}\label{ divisibilityact in fallible}
     \Fitt_\Lambda\Sel_{\mathcal F_\Lambda (\QT)} (\mathbf W_f) ^\vee\subset (\boldsymbol\lambda (\QT)) ^ 2\subset \Lambda.
 \end{equation}
 For the passage between characteristic ideal and fitting ideals, recall that the characteristic ideal of any $\Lambda$-module is the smallest divisorial ideal containing the Fitting ideal; cf. \cite[Corollary 3.2.9]{skinner2014iwasawa}. Steps 2 and 3 allow us to pass from (\ref{divisibility work planet levels}) to (\ref{ divisibilityact in fallible}).
  \end{proof}
  \begin{cor}
  Under the hypotheses of Theorem
  \ref{Heegner point main conjecture}, if additionally $\nu(N^-)$ is even,  then the Heegner point main conjecture holds for $f $ in $\Lambda\otimes \Q_p$; that is,
there is a pseudo-isomorphism of $\Lambda $-modules: $$\Sel(K_\infty, A_f[\p^\infty])^\vee\approx \Lambda \oplus M\oplus M$$
 for some torsion $\Lambda $-module $M $, and
$$\Char_\Lambda \left(\frac{\Sel(K, T_f\otimes \Lambda)}{\Lambda\boldsymbol{\kappa}(1)}\right)= \Char_\Lambda(M)$$
 as ideals of $\Lambda\otimes\Q_p$. If additionally the image of the Galois action on $\overline T_f $ contains a nontrivial scalar, then the equality is true in $\Lambda. $
  \end{cor}
  \begin{cor}\label{specialization corollary}
  Under the hypotheses of Theorem \ref{Heegner point main conjecture}, the bipartite Euler system $(\kappa (1,\cdot), \lambda (1,\cdot))$ of (\ref{construction Kolyvagin}) is nontrivial. 
  \end{cor}
  \begin{proof}
 Let $\set{\QT,\epsilon_\QT}\in \mathsf N$ be such that $\Sel_{\mathcal F(\QT)}(T_f) = 0,$ where again $\mathcal F$ is the standard Selmer structure for $T_f$. As noted in Step 1 of the proof of Theorem \ref{Heegner point main conjecture}, we have $\mathbbm 1(\Sel_{\mathcal F_\Lambda (\QT)}(\mathbf{W}_f))\neq 0$, so $\mathbbm 1(\boldsymbol\lambda(\QT)) \neq 0$; this implies $\lambda (1, \QT)\neq 0 $ by Remark \ref{ specialization rmk}.
  \end{proof}
  Corollary \ref{specialization corollary} is generalized by Theorem \ref{ appendix theorem} of the appendix.
  \subsection{Kolyvagin's conjecture}
Let $f$, $\p$, $K$, and $N^+N^-$ be as in (\ref{modular forms notation}) and (\ref{construction notation}). 
\subsubsection{}
For any $\mathsf m\in\mathsf K $, define the $\mathsf m $-transverse Selmer ranks
\begin{equation}
    r_\mathsf m ^\pm = \rank_\O\Sel_{\mathcal F (\mathsf m)} (T_f) ^\pm,
\end{equation}
 where $\pm$ refers to the $\tau $ eigenvalue; note that this is well-defined because the local conditions defining $\mathcal F (\mathsf m) $
 are all $\tau $-stable. When $\mathsf m = 1, $
 the $r_1 ^\pm $ are the classical Selmer ranks of $f $.
 \begin{prop}\label{Kolyvagin duality}
 For all $\mathsf m\mathsf l\in\mathsf K $, and for each $\delta\in\set{\pm}$,
 either:
 \begin{itemize}
     \item $r_{\mathsf m\mathsf l} ^\delta =r_\mathsf m ^\delta -1,$ $\loc_\mathsf l ^\delta(\Sel_{\mathcal{F} (\mathsf m)} (T_f))^\delta\neq 0,$ and $\partial_\mathsf  l ^\delta (\Sel_{\mathcal F (\mathsf m\mathsf l)} (T_f))^\delta = 0.$
     \item $r_{\mathsf m\mathsf l} ^\delta =r_\mathsf m ^\delta +1,$ $\loc_\mathsf l ^\delta(\Sel_{\mathcal{F} (\mathsf m)} (T_f))^\delta= 0,$ and $\partial_\mathsf  l ^\delta (\Sel_{\mathcal F (\mathsf m\mathsf l)} (T_f))^\delta \neq 0.$
 \end{itemize}
 \end{prop}
 \begin{proof}
 If $\mathcal F ^\mathsf l (\mathsf m,\QT) = \mathcal F (\mathsf m\mathsf l,\QT) +\mathcal F (\mathsf m,\QT) $ and $\mathcal F_\mathsf l (\mathsf m,\QT) =\mathcal F (\mathsf m\mathsf l,\QT)\intersection\mathcal F (\mathsf m,\QT) $, then we have a $\tau$-equivariant exact sequence
$$0\to \Sel_{\mathcal F _\mathsf l (\mathsf m,\QT)}(T_f)\to \Sel_{\mathcal F ^\mathsf l (\mathsf m,\QT)}(T_f)\to \H^1( K_\mathsf l, T_f), $$
where the image of the final arrow has rank two and is self-annihilating under the local Tate pairing by Propositions \ref{ dualityover DVR} and \ref{Greenberg Wiles over DVR}. Since the Tate pairing of two classes with opposite $\tau $
eigenvalues is necessarily zero, the proposition follows.
 \end{proof}
\begin{lemma}\label{Kolyvagin nontrivial loc}
Suppose given elements $c^\pm\in \H^1(K,T_f) ^\pm$, where $\pm$ is the $\tau $ eigenvalue. Then there exists a Kolyvagin-admissible ultraprime $\mathsf l$ such that $$c^\pm \neq 0 \implies \loc_\mathsf l ^\pm c^\pm.$$ If (\ref{contains scalar}) holds for $T_f $, then the same is true for elements $c ^\pm\in\H ^ 1 (K, T_f/\pi ^ j) $.
\end{lemma}
\begin{proof}
The proof of Theorem \ref{localization nontrivial} applies almost verbatim, except that  in the proof of Lemma \ref{just one localization} we will have two homomorphisms $\xphi^\pm \in \Hom_{G_K}(G_L, \overline T_f)^\pm$, and we must choose $g\in G_L$ so that $\xphi^\epsilon(g)$ has nonzero component in the $\epsilon$ eigenspace of $\tau$ for both signs $\epsilon$ (unless $\xphi^\epsilon$ is itself 0); for each $\epsilon$, this condition is satisfied outside a proper subgroup of $G_L$, so indeed there exists $g\in G_L$ such that both conditions are satisfied. With this modification, the rest of the proof applies unchanged. 
\end{proof}

\begin{lemma}\label{ KolyvaginTheorem lemma}
Suppose that the bipartite Euler system $(\kappa (1,\cdot),\lambda(1,\cdot)) $ of (\ref{construction Kolyvagin})  is nontrivial. Then, for all $\mathsf m\in\mathsf K $, $(\kappa (\mathsf m,\cdot),\lambda (\mathsf m,\cdot))$ is nontrivial.

In particular, for all $\mathsf m\in\mathsf K $ and  $\set {\QT,\epsilon_\QT}\in\mathsf N_\mathsf m $: $$ \rank_\O\Sel_{\mathcal F (\mathsf m,\QT)} (T_f) \leq 1\iff\begin {cases}\kappa (\mathsf m,\QT)\neq 0, &\nu(N^- ) + |\QT | \text{ even}\\\lambda (\mathsf m,\QT)\neq 0, &\nu(N^- ) + |\QT | \text{ odd}.\end{cases} $$ 
\end{lemma}
\begin{proof}
Proposition \ref{ Selmerstaff for Kolyvagin classes} implies that, for fixed $\mathsf m $, the pair $(\kappa (\mathsf m,\cdot), \lambda (\mathsf m,\cdot)) $
forms a bipartite Euler system with sign $\nu(N^-)$ for the self-dual Selmer structure $(\mathcal F(\mathsf m),\ST\union\mathsf m) $ on $T_f$. We will prove that, for any $\mathsf m\mathsf l\in \mathsf K$, if
$(\kappa(\mathsf m, \cdot), \lambda(\mathsf m, \cdot))$ is nontrivial then so is $(\kappa (\mathsf m\mathsf l,\cdot),\lambda (\mathsf m\mathsf l,\cdot)). $

Choose $\set{\QT,\epsilon_\QT}\in\mathsf N_\mathsf m ^{\nu(N^-)+1}$
such that $\Sel_{\mathcal F(\mathsf m,\QT)} (T_f) = 0$ and $\mathsf l \not\in \QT$; this is possible by Corollary \ref{the Selmer annihilation corollary}.
By Proposition \ref{Kolyvagin duality}, we may choose a nonzero $$d\in \Sel_{\mathcal F (\mathsf m\mathsf l, \QT)} (T_f).$$
Applying Theorem \ref{localization nontrivial} to $d$, let $\q $
be admissible with sign $\epsilon_\q $
such that $\q\not\in\QT\mathsf m\mathsf l$ and $\loc_\q d \neq 0.$ 
By Proposition \ref{Ray duality} for the Selmer structures $\mathcal F (\mathsf m,\QT\q)$ and $\mathcal F (\mathsf m,\QT)$, \begin{equation}
    \rank_\O \Sel_{\mathcal F (\mathsf m,\QT\q)} (T_f) = 1. 
\end{equation}
Hence, by hypothesis, $ \kappa(\mathsf m, \QT\q)$ generates $\Sel_{\mathcal F (\mathsf m,\QT\q)} (T_f)$ up to finite index, and in particular $\partial_\q\kappa(\mathsf m, \QT\q) \neq 0$. Now, taking the sum of local pairings and using Proposition \ref{ dualityover DVR},
\begin{equation}\label{pairing}0=\sum_\v \langle d, \kappa(\mathsf m, \QT\q)\rangle_\v = \langle d, \kappa(\mathsf m, \QT\q)_{\mathsf l} + \langle d, \kappa(\mathsf m, \QT\q)\rangle_\q.
\end{equation}
Since the latter pairing is nonzero by construction, the former is as well, and so, by Proposition \ref{ Selmerstaff for Kolyvagin classes}(1), $$\Res_\mathsf l  \kappa(\mathsf m, \QT\q)) \neq 0 \implies \kappa(\mathsf m\mathsf l, \QT\q) \neq 0.$$
\end{proof}
\subsubsection{}For any $\mathsf m\in\mathsf K$, define the vanishing order of the Kolyvagin system at $\mathsf m$:
\begin{equation}
    \nu_\mathsf m = \begin{cases} \min\set{|\mathsf n|\,:\,\mathsf n \in \mathsf K,\,\lambda (\mathsf n\union\mathsf m, 1)\neq 0}, &\nu(N ^ -)\text { odd,}\\
    \min\set{|\mathsf n |\,:\,\mathsf n \in \mathsf K,\,\kappa (\mathsf n\union \mathsf m, 1)\neq 0}, &\nu (N ^ -)\text { even}.\end {cases}
\end{equation}

\begin{cor}\label{ Selmer rank corollary}
If $(\kappa(1,\cdot), \lambda(1, \cdot)) $ is nontrivial, and in particular under the hypotheses of Theorem \ref{ appendix theorem}, 
we have for all $\mathsf m\in\mathsf K $:
\begin{itemize}
    \item If $\nu (N ^ -) $
    is odd, then $\nu_\mathsf m =\max\set {r_{\mathsf m} ^ +, r_{\mathsf m} ^ -} $ and $r_\mathsf m^\pm\equiv \frac{\epsilon_f -1} {2}\pmod 2.$
    \item If $\nu (N ^ -) $
    is even, then $\nu_\mathsf m =\max\set {r_{\mathsf m} ^ +, r_{\mathsf m} ^ -}-1$ and $\epsilon_f\cdot (-1)^{|\mathsf m| + \nu_\mathsf m + 1}$ is the larger $\tau$ eigenspace.
\end{itemize}
In particular, if $\rank_\O\Sel (K, T_f) = 1, $
then $L' (f/K, 1)\neq 0. $
\end{cor}
\begin{proof}
Note that the parity statement in Theorem \ref{DVR Euler system bound} implies $r_\mathsf m^++r_\mathsf m^- \equiv \nu(N^-) + 1\pmod 2$ for all $\mathsf m\in\mathsf K$
by Proposition \ref{Kolyvagin duality}.
So if we have $\rank_\O\Sel_{\mathcal F (\mathsf m\mathsf n)} (T_f)\leq 1, $
for some $\mathsf n\in\mathsf K$ such that $\mathsf m\intersection\mathsf n =\emptyset $, then the $\tau$-equivariant localization map $$\Sel_{\mathcal F (\mathsf m)}(T_f)\to\oplus_{\mathsf l\in\mathsf n} H ^ 1_{\unr} (K_\mathsf l, T_f) $$
has kernel of rank either zero (if $\nu(N^-)$ is odd), or at most one (if $\nu(N^-)$ is even). It follows that $\nu_\mathsf m \geq \max\set{r_\mathsf m^+,r_\mathsf m^-}$ in the former case and $\nu_\mathsf m \geq \max\set{r_\mathsf m^+,r_\mathsf m^-}-1 $ in the latter. The opposite inequality follows readily
 from repeated applications of Proposition \ref{Kolyvagin duality} and Lemma \ref{Kolyvagin nontrivial loc}. The additional claims about parity and $\tau$ eigenvalues follow from Proposition \ref{ parity proposition}.
 
 For the final statement, take $\mathsf m = 1. $
 The corollary implies $\nu_1 = 0, $ so $\kappa (1, 1)\neq 0. $ Since $\kappa (1, 1) $ is the Kummer image of the classical Heegner point $y_K\in E (K), $
 the result follows from the Gross-Zagier theorem of \cite{zhang2001gross}.
\end{proof}

It remains to relate the vanishing of the patched Kolyvagin classes to the classical vanishing order 
\begin{equation}
    \nu_{\text {classical}}\coloneqq\begin{cases}
    \min\set {\nu (m)\,:\,\overline\lambda (m, 1)\neq 0}, &\nu (N ^ -)\text { odd,}\\
    \min\set {\nu (m)\,:\,\overline c (m, 1)\neq 0}, &\nu (N ^ -)\text { even.}\end{cases}
\end{equation}
\begin{cor}
\label{classical vanishing order}
If $(\kappa (1,\cdot),\lambda (1,\cdot)) $ is nontrivial, and in particular under the hypotheses of Theorem \ref{ appendix theorem}, $\nu_{\text {classical}} $
is finite. If (\ref{contains scalar}) holds for $f $, then $\nu_{\text {classical}} =\nu_1, $ and in particular:
\begin{itemize}
    \item If $\nu (N ^ -) $ is odd, then $\nu_{\text {classical}} = \max\set {r_1 ^ +, r_1 ^ -} $ and $r_1 ^\pm\equiv\frac {\epsilon_f -1} {2}\pmod 2. $
    \item If $\nu (N ^ -) $ is even, then $\nu_{\text {classical}} =\max\set {r_1 ^ +, r_1 ^ -} -1 $ and $\epsilon_f\cdot (-1) ^ {1+\nu_{\text {classical}}} $
    is the larger $\tau $ eigenspace.
\end{itemize}
\end{cor}
\begin{proof}
The finiteness of the classical vanishing order is clear by construction: if a patched Kolyvagin class or element is nontrivial, then infinitely many of the classical Kolyvagin classes or elements defining it are nontrivial. This also shows $\nu_{\text {classical}}\leq\nu_1. $
We will check that equality holds under the condition (\ref{contains scalar}). Suppose first $\nu (N ^ -) $ is even. We abbreviate by $c_j (m, 1)\in H ^ 1 (K, T_j) $
the image of $\overline c (m, 1) $ when $v_\p (I_m)\geq j. $ Given some nonzero $c_j (m, 1) $,
one may show as in \cite[p. 309]{mccallum1991kolyvagin}  that there exist classes $c_j (m_n, 1) \neq 0$
 with $v_\p (I_{m_n}) \to\infty $ and $\nu (m_n) =\nu (m).$ (In \cite{mccallum1991kolyvagin}, additional hypotheses are put on the image of the Galois action, but the argument goes through by invoking Lemma \ref{Kolyvagin nontrivial loc}.) In particular, the sequence $m_n $ defines a nonzero $\kappa (\mathsf m, 1)$ witnessing $\nu_1\leq\nu_{\text {classical}} $.
 
 Now suppose that $\nu (N ^ -) $ is odd, and that $\lambda_j (m, 1)\neq 0 $ where $\nu (m) =\nu_{\text {classical}}$. 
 We choose an auxiliary $\set{\mathsf q,\epsilon_{\mathsf q}}\in\mathsf N $ with the following properties:
\begin{itemize}
    \item $\epsilon_{\mathsf q} $ is the sign of the larger $\tau $ eigenspace in $\Sel_{\mathcal F} (T_f). $
    \item The localization map $\loc_{\mathsf q} $
    is trivial on $\Sel_{\mathcal F} (T_f). $
    
\end{itemize}
To ensure the second condition, we may choose $\Frob_{\mathsf q}\in G_\Q $ to be a complex conjugation. Let $\set{q_n,\epsilon_{q_n}} $
represent $\set {\mathsf q,\epsilon_{\mathsf q}} $
as in Remark \ref{level raising remark}.
 Once again, the argument of \cite[p. 309]{mccallum1991kolyvagin} implies that, for each $n $, there exists $m_n $ with $c_j (m_n, q_n)\neq 0 $
 and $v_\p (I_{m_n})\to\infty. $ We therefore obtain a nonzero patched class $\kappa(\mathsf m,\mathsf q) $ with $ |\mathsf m | =\nu_{\text {classical}}. $ Repeating the argument of Lemma \ref{ KolyvaginTheorem lemma}, it follows that $\nu_1\leq\nu_{\text {classical}} +1. $
 For contradiction, we assume that $$\nu_{\text {classical}}=\nu_1-1 = r ^ {\epsilon_{\mathsf q}} -1. $$
 This implies $\partial_\q\kappa (\mathsf m,\q) = 0, $
 so by Lemma \ref{ KolyvaginTheorem lemma} and Proposition \ref{Ray duality}, we conclude $$\rank_\O\Sel_{\mathcal F (\mathsf m, 1)} (T_f)= 2. $$
 However, by Proposition \ref{Kolyvagin duality} and the assumption $ |\mathsf m | = \nu_1-1 = r ^ {\epsilon_{\mathsf q}}-1, $ $\rank_\O\Sel_{\mathcal F (\mathsf m, 1)} (T_f) ^ {\epsilon_\q} $ is odd, hence equal to one. Proposition \ref{Kolyvagin duality} then implies $$\Sel_{\mathcal F (\mathsf m, 1)} (T_f) ^ {\epsilon_\q}\subset\Sel_{\mathcal F} (T_f) ^ {\epsilon_\q}, $$
 so $$\loc_\q (\Sel_{\mathcal F (\mathsf m, 1)} (T_f)) = 0$$
 by the choice of $\q $.
 However, this contradicts Proposition \ref{Ray duality}, so we must have $\nu_1 =\nu_{classical}. $
 
 
\end{proof}

\appendix
\section{Kolyvagin's conjecture for inert or non-ordinary $p$}\label{ appendix section}
\subsection{The main result}
In this appendix, we shall prove the following:
\begin{thm}\label{ appendix theorem}
Let $f$ be a non-CM modular form of weight two, level $N$, and trivial character, with a prime $\p$ of its ring of integers $\O_f$, and let $K$ be an imaginary quadratic field. Assume:
\begin{itemize}
    \item $N = N^+N^-$, where every factor of $N^+$ is split in $K $, and $N ^ - $
is a squarefree product of an even number of 
 primes inert in $K $.
 \item The residue characteristic $p$ of $\p $ does not divide $2D_KN$.
 \item The modulo $\p$ representation $\overline T_f $
associated to $f $ 
is absolutely irreducible; if $p = 3, $ assume that $\overline T_f $ is not induced from a character of  $G_{\Q\sqrt{-3}}$.
\item If $p$ is inert in $K$ or  or $a_p$ is not a $\p$-adic unit, then there exists some prime $\l_0||N$.
\item If $a_p$ is not a $\p$-adic unit, then either $\l_0$ may be chosen above so that $A_f$ has non-split toric reduction at $\l_0$, or the image of the Galois action on $T_f$ contains a conjugate of $SL_2 (\Z_p) $.
\end{itemize} Then  $(\kappa(1,\cdot),\lambda (1,\cdot))$ is nontrivial.
\end{thm}

\subsubsection{}
If $p $ is split in $K $, then this is simply Corollary \ref{specialization corollary}. 
If $p$ is non-ordinary or inert in $K $, then the anticyclotomic main conjecture is currently not known in full generality; however, since all we are interested in specialization the trivial character, we will show that
the result may instead be obtained, more circuitously, by combining  main conjectures for quadratic twists of $f $. The proof applies equally well to the split ordinary case.

\subsection{Comparing periods}
\subsubsection{}Let $f $ be a modular form of weight two, level $N$, and trivial character,  with ring of integers $\O_f$ of its coefficient field, and let $\p\subset\O_f$ be an ordinary prime lying over $p\nmid 2N$, with associated completion $\O$; we assume that $\overline T_f $ is absolutely irreducible.  There are two ways to normalize the anticyclotomic $p $-adic $L $-function, as explained in \cite{vatsal2003special, pollack2011anticyclotomic}.
For any factorization $N = N_1N_2,$ where $N_1$ and $N_2$ are coprime, the congruence ideal $\eta_f(N_1,N_2)\subset \O$
is defined as
\begin{equation}
    \pi_f (\Ann_{\T_{N_1, N_2}} (\ker\pi_f))\cdot,
\end{equation}
where $\pi_f:\T_{N_1, N_2}\to\O $ is the projection giving the Hecke eigenvalues of $f $.
Hida's canonical period \cite{hida1988modules} is defined (up to $p$-adic units) by:
\begin {equation}
\Omega_f ^ {can} =\frac {(f,f)}{\eta_f(N, 1)},
\end{equation}
where $(f, f) $
is the Peterson inner product.

On the other hand, if $N=N_1N_2$ where $N_2$ is squarefree with an odd number of prime factors, then the $f$-isotypic part of the Hecke module $ \Z_p[X_{N_1,N_2}]$ is free of rank one, generated by an element $\phi_{f,N_2}$. Gross's period is defined as:
\begin{equation}
\Omega_{f,N_2} =\frac {(f,f)}{\langle\phi_{f, N_2},\phi_{f, N_2}\rangle},
\end{equation}
where $\langle \cdot, \cdot \rangle$
is the canonical intersection pairing on $\Z_p [X_{N_1,N_2}]$.  This period occurs naturally in anticyclotomic Iwasawa theory due to the well-known special value formula of Gross.
\begin{prop}
Let $K $
be an imaginary quadratic field of discriminant prime to $Np $, and suppose that $N = N ^ + N ^ - $
where all factors of $N ^ + $
are split in $K $, and $N ^ - $
is a squarefree product of an odd number primes inert in $K $. Then $L(f/K, 1) \in \overline\Q\cdot\Omega_{f,N^-}$ and the element $\lambda (1)\in\Lambda$ constructed in (\ref{construction Kolyvagin}) satisfies:
\begin{equation}
  \frac{L(f/K, 1)}{\Omega_{f, N ^ -}} = \lambda (1)^ 2
\end{equation}
up to $p$-adic units.
\end{prop}

\subsubsection{}
 Consider the N\'{e}ron model $\mathcal A_f$ of $A_f$ over $\Z_\l$. The
 Tamagawa numbers over $K $ are: 
\begin{equation}
    t_f(\l)=\lg_{\O}\Phi(\mathcal A_f)(\O_K/\l)_{\p}
\end{equation}
We also require a variant:
\begin{equation}
    c_f(\l)=\lg_{\O}\Phi(\mathcal A_f)(\overline{\F}_\l)_\p
\end{equation}
The number $c_f(\l)$ is the maximal exponent $e$ such that $A_f[\p^e_f]$ is unramified at $\l$. 
If $\l$ is inert in $K$, then  $t_f(\l)=c_f(\l).$ 
The following theorem  generalizes  \cite{ribet1997parametrizations,khare2003isomorphisms} and \cite[Theorem 6.8]{pollack2011anticyclotomic}.  

\begin{prop}\label{rtstufthm}
Suppose $N = N_1N_2 $ where $N_2 $
is squarefree with an odd number of prime factors and coprime to $N_1$.
For any $\l_0||N$, we have: $$v_{\p}\eta_f(N,1)-v_{\p}\langle\phi_{f,N_2},\phi_{f,N_2}\rangle\geq\sum_{\l|N_2}c_f(\l)-\sigma(N_2)c_f(\l_0).$$
\end{prop}
\begin{proof}
For a decomposition $N=N_1'N_2'$ with $N_2'$ the squarefree product of an even number of primes that do not divide $N_1'$, one defines $\delta(N_1',N_2')\subset\O$ to be the ``degree'' of an optimal modular parametrization $J_{N_1',N_2'}\to A_f$ as explained in \cite{pollack2011anticyclotomic, khare2003isomorphisms}.
 By \cite[Proposition 6.6]{pollack2011anticyclotomic}, we have: 
\begin{equation}\label{pollackprop}
    v_{\p} \delta(N,1)=c_f(\l_0)+v_{\p}\langle\phi_{f,\l_0},\phi_{f,\l_0}\rangle.
\end{equation}
On the other hand, \cite[Lemma 4.17]{darmon1995fermat} implies that:
\begin{equation}
    v_{\p}\eta_f(N/\l_0,\l_0)\geq v_{\p}\langle\phi_{f,\l_0},\phi_{f,\l_0}\rangle.
\end{equation}
Because  $\T_{N/\l_0,\l_0}$ is a quotient of $\T_{N,1}$, we conclude that:
\begin{equation}\label{etatodelta}
      v_{\p}\eta_f(N,1)\geq v_{\p}\delta(N,1)-c_f(\l_0).
\end{equation}
We apply \cite[Proposition 6.6]{pollack2011anticyclotomic} again, this time to the decomposition $N=N_1N_2$ and any  $r|N_2.$ This yields: 
\begin{equation}\label{moveoneout}
v_{\p}\delta(N_1r,N_2/r)=c_f(r)+v_{\p}\langle \phi_{f,N_2},\phi_{f, N_2}\rangle.
\end{equation} 
If $N_2$ is prime, this is sufficient to conclude.  If not, we may choose $r\neq\l_0.$

The results of Ribet-Takahashi and Khare \cite{ribet1997parametrizations, khare2003isomorphisms} imply that:
\begin{equation}\label{rt}
    v_{\p}\delta(N,1)\geq v_{\p}\delta(N_1r,N_2/r) +\sum_{\l|N_2/r}c_f(\l)-\sigma(N_2/r)\cdot c_f(\l_0).
\end{equation}
(If $\l_0|N_1,$ we are using the fact that both $r$ and $\l_0$ exactly divide $N_1r.$) Combining (\ref{etatodelta}), (\ref{moveoneout}), and (\ref{rt}) completes the proof.

\end{proof}
\begin{rmk}
 If $\l_0$ is residually ramified, the inequality is an equality. In \cite[Theorem 6.4]{zhang2014selmer}, more restrictive conditions are given under which this result holds.

\end{rmk}
\subsubsection{}We will also require the following related result:
\begin{prop}\label{muupper}
Let $f$ and $\p$ be as above and suppose $\l_0||N$.  Then, in the notation of the proof of Theorem \ref{rtstufthm},  $$v_{\p}\eta_f(N,1)-v_{\p}\langle\phi_{f,\l_0},\phi_{f,\l_0}\rangle\leq c_f(\l_0).$$
\end{prop}
\begin{proof}
Let $J=J_0(N)$ be the modular Jacobian and let $\T$ be the full Hecke algebra of level $N$. Write $\pi:\T_\m\to \O_{f,\p}$ for the projection associated to $g$ and let $I$ be its kernel.  The claim will follow from (\ref{pollackprop}) once we establish \begin{equation}\label{etadelta}
v_{\p}\eta_f(N,1)\leq v_{\p}\delta(N,1).
\end{equation}
Indeed, if $J\to A$ is an optimal parametrization, then the dual map $A^\vee\to J^\vee $ is an inclusion. The composition $$\xphi:J \to A\to A^\vee\to J^\vee\xrightarrow{w_N}J$$ is a Hecke-equivariant endomorphism; by (\ref{all endomorphisms}), its image in $\End(J)_{\m}$ may be identified with some $y\in \T_\m$. Because $\im\xphi\subset J[I],$ we have $y\in \Ann(I)$.  By the definition of $\delta(N,1)$, $$(\pi(y))=\delta(N,1)\subset \O.$$  This implies (\ref{etadelta}).
\end{proof}

\subsection{Cyclotomic Iwasawa theory: ordinary case}
\subsubsection{}
Let $\Lambda_{\Q_\infty} =\O\llbracket \Gal (\Q_\infty/\Q)\rrbracket$
be the cyclotomic Iwasawa algebra. We denote by $\mathbbm 1:\Lambda_{\Q_\infty}\to\O $
and $\mathbbm 1:\Lambda\to \O$
the specializations at the trivial character. 
If $\p$ is ordinary and $\Sigma $ is a finite set of rational primes, we consider the $\Sigma $-ordinary  cyclotomic Selmer group
$$\Sel (\Q_\infty, W_f) =\ker\left (H ^ 1 (\Q,\mathbf W_f)\to\prod_{v\not\in\Sigma\union\set {p}} H ^ 1 (I_v,\mathbf W_f)\times\frac {H ^ 1 (\Q_p,\mathbf W_f)} {H ^ 1_{\ord} (\Q_p,\mathbf W_f)}\right), $$
and
denote by $Ch_{\Q_\infty, f}^\Sigma\subset\Lambda_{\Q_\infty} $
the characteristic ideal of its Pontryagin dual.


From the work of Skinner and Urban, we deduce the following result.

\begin{thm}[Skinner-Urban]\label{skinner urban prop}
Let $K $
be an imaginary quadratic field of discriminant prime to $Np$ in which $p$ splits. Assume that $\p$ is good ordinary for $f$ and that:
\begin{itemize}
    \item the mod $\p$ representation  $\overline T_f$ is absolutely irreducible;
    \item $N = N_1N_2,$ where every factor of $N_1 $ is split in $K$ and $N_2$ is the squarefree product of an odd number primes inert in $K $.
    \end{itemize}
    Then there exists an element $\alpha\in \Lambda_{\Q_\infty} $ such that $\mathbbm 1 (\alpha)$ divides $$\frac{\Omega_{f,N_2}}{\Omega_f^{can}}\sim \frac{\eta_f(N,1)}{\langle \phi_{f,N_2},\phi_{f, N_2}\rangle}$$ in $\O$ and $$(\alpha)Ch_{\Q_\infty, f}Ch_{\Q_\infty, f\otimes \chi_K} \subset (L _p (\Q_\infty, f))( L _p (\Q_\infty, f\otimes\chi_K)).$$
\end{thm}
\begin{proof}
Recall the divisibility established in the course of the proof of Theorem \ref{skinner urban} for the
Fitting ideal of the 3-variable Selmer group:
\begin{equation}\label{Sitting divisible 2}
   (\widetilde\alpha) Fitt_{K_\infty}^\Sigma(\mathbf g)\subset (\mathcal L_{\mathbf f, K} ^\Sigma),
\end{equation}
 where $\widetilde\alpha \in \mathbb I [\Gamma_{K}^ +]$ may be chosen such that $\widetilde\alpha$ specializes to a unit multiple of $\Omega_{f, N_2}/\Omega_f ^ {can} $ at the trivial character (by \cite{vatsal2003special}).
 By Lemma 3.2.5, Corollary 3.2.9(i), and Corollary 3.2.20(iii) of  \cite{skinner2014iwasawa}, specializing to the cyclotomic variable yields a divisibility
\begin{equation}
    (\alpha )Ch^\Sigma_{\Q_\infty, f}Ch^\Sigma_{\Q_\infty, f\otimes \chi_K} \subset( L ^\Sigma_p (\Q_\infty, f))( L ^\Sigma_p (\Q_\infty, f\otimes\chi_K)),
\end{equation}
where $\alpha $ is the image of $\widetilde\alpha$. The desired divisibility for the imprimitive $L$-functions and Selmer groups  follows by \cite[Proposition 3.2.18]{skinner2014iwasawa}. 
\end{proof}

 \subsection{Cyclotomic Iwasawa theory: general case}
 \subsubsection{}
Kato gave a formulation \cite{kato2004p} of the cyclotomic main conjectures which also applies to non-ordinary primes. For any finite set of rational primes $\Sigma $, consider the strict Selmer group:
$$\Sel_{\str}(\Q_\infty, W_f) =\ker\left (H ^ 1  (\Q,\mathbf W_f)\to\prod_{v\not\in\Sigma\union\set {p}} H ^ 1 (I_v,\mathbf W_f)\times H ^ 1 (\Q_p,\mathbf W_f)\right), $$
and
denote by $Ch_{\Kato, f}^\Sigma\subset\Lambda_{\Q_\infty} $
the characteristic ideal of its Pontryagin dual.

The Iwasawa cohomology $H ^ 1 (\Q ^\Sigma/\Q,\mathbf T) $ (for $\Sigma $ the set of primes dividing $Np $) is free of rank one when $\overline T_f $  is absolutely irreducible, and under that hypothesis Kato defined an element $z_{\Kato}\in H ^ 1 (\Q ^\Sigma/\Q,\mathbf T) $ which is closely related to $L $-values of $f $. 
\begin{rmk}
In \cite[Theorem 12.4]{kato2004p}, it is only asserted that $z_{\Kato}$ is integral when the image of the Galois action on $T_f $
contains a conjugate of $SL_2 (\Z_p) $. However, the proof (13.14 of \textit{loc. cit.}) only requires the fact that any two $\O$-lattices in $T_f\otimes \Q_p$ are homothetic, which holds whenever $\overline T_f$ is absolutely irreducible. Note that  $z_{\Kato} = \mathbf z_\gamma^{(p)}$ in the notation of \cite{kato2004p}, where $\gamma\in T_f$ is the sum of any generators of the two eigenspaces for complex conjugation.
\end{rmk}
\subsubsection{}Kato's main conjecture \cite[Conjecture 12.10]{kato2004p} then asserts that:
\begin{equation}\label{Canto conjecture}
    \Char_{\Lambda_{\Q_\infty}}\left(\frac{H ^ 1 (\Q ^\Sigma/\Q,\mathbf T)}{\Lambda_{\Q_\infty}\cdot z_{\Kato}}\right)=Ch_{\Kato,f}.
\end{equation}
If $\p$ is ordinary, then this conjecture is equivalent to the usual cyclotomic main conjecture for $f $ by \cite[17.13]{kato2004p}.
We abbreviate by $Z_{\Kato,f}\subset \Lambda_{\Q_\infty}$
the left side of (\ref{Canto conjecture}).

 The following is the analogue of Theorem \ref{skinner urban prop} in the non-ordinary case, due to Wan.
 \begin{thm}[Wan]\label{wan cyclotomic}
 Let $K$ be an imaginary quadratic field in which $p$ splits. Assume that $\p$ is a good, non-ordinary prime for $f$ and that:
 \begin{itemize}
     \item For all primes $\l|N$ ramified in $K$, that $T_f|_{G_{\Q_\l}}$ is Steinberg with sign $-1$. At least one such prime exists.
     \item For all $\l|N$, $\l$ is either split or ramified in $K$.
 \end{itemize}
 Then $$Ch_{\Kato, f}\cdot Ch_{\Kato, f\otimes \chi_K}\subset Z_{\Kato,f}\cdot Z_{\Kato, f\otimes \chi_K} \text{ in }\Lambda_{\Q_\infty}.$$
 \end{thm}
 \begin{proof}
 This is proven  in \cite[p. 29]{wan2016iwasawa}; compare to the proof of Corollary 3.32, where it is assumed that $Z_{\Kato,f \otimes \chi_K}\subset \Ch_{\Kato, f\otimes\chi_K}$.
 \end{proof}
 Additionally, Kato \cite{kato2004p} has proven one direction of his conjecture in our setting:
\begin{thm}[Kato]\label{kato theorem}
 Let $f$ be a modular forms of weight two, level $N $, and trivial character, and $\p\subset \O_f $
 a prime of good reduction with odd residue characteristic. Then $Ch_{\Kato,f}\neq 0$  and $$Z_{\Kato, f}\subset Ch_{\Kato, f}$$ in $\Lambda_{\Q_\infty}\otimes\Q_p.$
 In particular, if $a_p$ is a $\p$-adic unit, then $\Sel (\Q_\infty, W_f)$ is $\Lambda_{\Q_\infty}$-cotorsion  and $$L_p (\Q_\infty,f)\subset Ch_{\Q_\infty, f}$$ in $\Lambda_{\Q_\infty}\otimes\Q_p.$
If the image of the Galois action on $T_f $ contains $SL_2 (\Z_p) $, then all of the inclusions hold  in $\Lambda_{\Q_\infty}$.
 \end{thm}
 \subsubsection{}
 Denote by $\mu (f) $ the $\mu $-invariant of $Ch_{\Q_\infty, f}$ or $Ch_{\Kato, f}$
 in the ordinary or non-ordinary case, respectively.
 To control the powers of $\p $ in Theorem \ref{kato theorem}, we will use the following.
 \begin{lemma}\label{control invariant}
 Let $f $ and $g $ be modular forms of weight two and trivial character such that $\overline T_f $ is absolutely irreducible. Suppose that $f $ and $g $ have a congruence modulo $\p^j$, i.e. there is a common completion $\O$ of $\O_f $ and $\O_g $ and, in some basis, a congruence of $\O$-valued associated Galois representations $$T_f\equiv T_g\pmod {\p ^ j}. $$ If $\mu (f) <j $, then $\mu (g) =\mu (f). $
 \end{lemma}
 \begin{proof}
 For the sake of notation, assume $\p$ is ordinary; this makes no difference to the proof.
 By \cite{greenberg2000iwasawa}, $\mu(f)$ is also the $\mu $-invariant of $Ch_{\Q_\infty, f}^\Sigma$  for any finite set of primes $\Sigma $, and likewise for $g $. If $\Sigma $ contains all primes dividing the level of either $f $ or $g $, then we have:
 \begin{equation}\Sel ^\Sigma (\Q_\infty, W_f) [\p ^ j]\simeq\Sel ^\Sigma (\Q_\infty, W_g) [\p ^ j] \end{equation} as $\Lambda_{\Q_\infty} $-modules.
 Let $M_f$ and $M_g$  be the Pontryagin duals of
 $\Sel ^\Sigma (\Q_\infty, W_f)$ and $\Sel ^\Sigma (\Q_\infty, W_g) [\p ^ j] $, respectively, and let $\P=(\p)\subset\Lambda_{\Q_\infty} $. Then we have a congruence $$M_f\otimes \Lambda_{\Q_\infty}/\P ^ j\simeq M_g\otimes\Lambda_{\Q_\infty}/\P ^ j. $$
 Since $\mu(f) = \lg M_{f,(\P)} < j,$ where $(\P)$ denotes the localization, 
 \begin{equation}
     M_{f,(\P)} \otimes \Lambda_{\Q_\infty}/\P ^ j =  M_{f,(\P)} \otimes \Lambda_{\Q_\infty}/\P ^ {j-1},
 \end{equation}
 which implies the same for $g$. Therefore $M_{g,(\P)}\otimes\Lambda_{\Q_\infty}/\P ^ j = M_{g,(\P)}$ and the result follows.
 \end{proof}
 
 \begin{proof}[Proof of Theorem \ref{ appendix theorem}]
 Let us suppose first that $A_f$   has non-split toric reduction at $\l_0$ if $\p$ is non-ordinary.
  Fix once and for all an auxiliary quadratic imaginary field $\mathcal K$, not contained in the fixed field $K(T_f)$, such that:
  \begin{itemize}
      \item If $\p$ is ordinary, then $\l_0$ is inert in $\mathcal{K}$ and every other factor of $Np $ is split in $\mathcal K $.
      \item If $\p $ is non-ordinary, then $\l_0 $ is ramified in $\mathcal K $ and every other factor of $Np $ is split in $\mathcal K $.
  \end{itemize}
 As in the proof of Theorem \ref{Heegner point main conjecture}, begin by
 applying Proposition \ref{killing Selmer} and Theorem \ref{relative def} to obtain some $\set {\QT, \epsilon_\QT}\in\mathsf N $, represented by $Q_n $, and a resulting sequence of newforms $g_n$ of $NQ_n$; we make sure to choose each $\q\in \QT$ such that $\Frob_\q$ has trivial image in $\Gal (\mathcal K/\Q) $, which is clearly possible.  
 \begin{claim}
 There exists a constant $C $, depending only on $f $, such that
  \begin{equation}\label{Uniform L value}
     v_\p\left(\frac {L (g_n/K,1)} {\Omega_{g_n}^{can}} \right)\leq\lg_\O\Sel(K,W_{g_n})+\sum_{\l|NQ_n}t_{g_n}(\l) + C
 \end{equation}
 for $\mathfrak F$-many $n$.
 \end{claim}
\begin{proof}[Proof of claim]
Consider first the ordinary case. By Lemma \ref{control invariant} and Theorem \ref{kato theorem}, 
\begin{equation}
    \p^{\mu(f\otimes\chi_\mathcal K)}\cdot (L_p(\Q_\infty, g_n\otimes\chi_\mathcal K))\subset Ch_{\Q_\infty, g_n\otimes\chi_\mathcal K} \text{ in } \Lambda_{\Q_\infty} 
\end{equation}
for $\mathfrak F $-many $n $. 
By Theorem \ref{skinner urban prop} for $g_n$, for $\mathfrak F $-many $n$ we have
\begin{equation}
      (\alpha)\cdot \p^{\mu(f\otimes\chi_K)} \cdot Ch_{\Q_\infty, g_n}\cdot Ch_{\Q_\infty,g_n\otimes \chi_\mathcal K}\subset (L_p (\Q_{\infty}, g_n))\cdot Ch_{\Q_\infty,g_n\otimes \chi_\mathcal K},
\end{equation}
 where, by Proposition \ref{muupper}, $\mathbbm 1 (\alpha)$ divides $\p^{c_f(\l_0)}$ in $\O$.
 Since $Ch_{\Q_\infty, g_n\otimes\chi_\mathcal K}\neq 0$, and since characteristic ideals are divisorial, we conclude that 
 \begin{equation}
     (\alpha)\cdot \p^{\mu(f\otimes\chi_K)} \cdot Ch_{\Q_\infty, g_n}\subset (L_p (\Q_{\infty}, g_n)).
 \end{equation}
  Applying the same argument to $g_n\otimes\chi_K $, we have:
 \begin{equation}\label{Silly equation}
      (\alpha)^2\cdot \p^{{\mu(f\otimes\chi_K\otimes\chi_\mathcal K)}+\mu( f\otimes\chi_\mathcal K)} \cdot Ch_{\Q_\infty, g_n}\cdot Ch_{\Q_\infty, g_n\otimes\chi_K}\subset (L_p (\Q_{\infty}, g_n)\cdot (L_p(\Q_\infty, g_n\otimes\chi_K)).
 \end{equation}
 The result now follows from standard interpolation properties of both sides of (\ref{Silly equation}), cf. e.g. \cite[Theorem 3.6.11]{skinner2014iwasawa}.
 
 The non-ordinary case is similar: combining Theorem \ref{wan cyclotomic} and Theorem \ref{kato theorem}, we have for $\mathfrak F $-many $g_n $
 $$\p^{\mu(f\otimes\chi_\mathcal K)} \cdot Ch_{\Kato, g_n}\subset Z_{\Kato,g_n},$$
 and likewise for the twist $g_n\otimes\chi_K $.
 Since $Z_{\Kato,g_n}$ is principal, the result follows as in \cite[Corollary 3.35]{wan2016iwasawa}. Note that the $p$-adic Tamagawa factor appearing there is trivial by \cite[Proposition II.2]{berger2002tamagawa}.
 \end{proof}
 As in Step 3 of the proof of Theorem \ref{Heegner point main conjecture}, $\#\Sel(K,W_{g_n})=\#\Sel_{\mathcal{F}(\QT)}(W_f)<\infty$ for $\mathfrak F $-many $n $ (the local conditions at $v|N$ may be compared in the same way, and at $v|p$ we use \cite[Lemma 7]{gross2012local}). 
 
 Now, by combining the claim above with Proposition \ref{rtstufthm}, we have for $\mathfrak F $-many $n $:   \begin{equation}\label{Uniform L value 2}
     v_\p\left(\frac {L (g_n/K,1)} {\Omega_{g_n,N^-Q_n}} \right)\leq\lg_\O\Sel(K,W_{g_n})+\sum_{\l|N^+}t_{g_n}(\l) + C'
 \end{equation}
  for a constant $C' $
  that does not depend on $n$.  In particular, for $\mathfrak F $-many $n$, $L(g_n/K,1)\neq 0$,  which by parity considerations implies that $\nu(N^-)+|\QT |$ is odd.
  Exactly as in Step 2 of the proof of Theorem \ref{Heegner point main conjecture}, we then conclude from (\ref{Uniform L value 2}) that $\lambda (1,\QT)\neq 0. $

  If we assume instead that $\p$ is non-ordinary but the image of the Galois action on $T_f$ contains a conjugate of $SL_2(\Z_p)$, then, rather than choosing $\mathcal K $ at the beginning, we choose $\set{\QT,\epsilon_\QT}\in \mathsf N $ as in the proof of Theorem \ref{Heegner point main conjecture}, but we also ensure that $\epsilon_\QT(\q) = -1$ for at least one $\q$.
  If $\set{Q_n, \epsilon_{Q_n}}$ represents $\set{\QT,\epsilon_\QT}$ and $q_n|Q_n$ represents $\q$, then for each $n$, we choose an auxiliary imaginary quadratic field $\mathcal K_n$ such that $q_n$ is ramified in $\mathcal K_n$ and all other factors of $NQ_np$ are split. Note that, by the proof of \cite[Lemma 6.15]{fakhruddin2019relative}, the image of the Galois action on $T_{g_n}$ contains a conjugate of $SL_2(\Z_p)$ for $\mathfrak F$-many $n$. The argument is then the same as above; $\mu(f\otimes \chi_{\mathcal K_n})$ and $\mu(f\otimes \chi_{\mathcal K_n}\otimes \chi_K)$ may not be uniformly bounded in $n$, but by Kato's result these error terms are not needed under the large-image hypothesis. 
 \end{proof}
  \begin{rmk}
In all cases, we crucially use the prime $\l_0 $
  to make the error terms uniform in $n$.
  \end{rmk}

\bibliography{mybib.bib}
\bibliographystyle{plain}
\section*{Declarations}
\textbf{Funding.}
This research was funded by NSF grant \#DGE1745303.

\textbf{Conflicts of interest.}
None.

\textbf{Data transparency.}
Not applicable.

\textbf{Code availability.}
Not applicable.
\end{document}